\documentclass[11pt,reqno]{amsart}
\usepackage{amsmath,amssymb,amsthm,mathtools,booktabs}
\usepackage[margin=1.05in]{geometry}
\usepackage[hidelinks]{hyperref}
\hypersetup{
  pdftitle={Real equiangular lines in dimensions 18--20 and sharp fixed-angle bounds},
  pdfauthor={Valjakas Djaljapayan, Zhen-Qi Liao, and Wei-Hsuan Yu},
  pdfsubject={Equiangular line maxima in dimensions 18, 19, and 20 and at prescribed angles}
}
\allowdisplaybreaks[2]
\renewcommand{\paragraph}[1]{\par\addvspace{4pt}\noindent\textbf{#1}\enspace\ignorespaces}
\newtheorem{theorem}{Theorem}[section]
\newtheorem{lemma}[theorem]{Lemma}
\newtheorem{proposition}[theorem]{Proposition}

\newtheorem{corollary}[theorem]{Corollary}
\theoremstyle{remark}

\newcommand{\R}{\mathbb R}
\newcommand{\Z}{\mathbb Z}
\newcommand{\Q}{\mathbb Q}
\newcommand{\F}{\mathbb F}
\DeclareMathOperator{\tr}{tr}
\DeclareMathOperator{\rank}{rank}
\newcommand{\rk}{\rank}
\DeclareMathOperator{\supp}{supp}
\DeclareMathOperator{\wt}{wt}

\newcommand{\ip}[2]{\langle #1,#2\rangle}
\newcommand{\sq}[1]{\lVert #1\rVert^2}
\newcommand{\eps}{\varepsilon}

\title[Equiangular lines and sharp fixed-angle bounds]
      {Real equiangular lines in dimensions 18--20 and sharp fixed-angle bounds}
\author[V. Djaljapayan]{Valjakas Djaljapayan}
\address[All authors]{Department of Mathematics, National Central University, No. 300, Zhongda Road, Zhongli District, Taoyuan City 320317, Taiwan}
\email[Valjakas Djaljapayan]{115281004@cc.ncu.edu.tw}
\author[Z.-Q. Liao]{Zhen-Qi Liao}
\email[Zhen-Qi Liao]{l0911487082@gmail.com}
\author[W.-H. Yu]{Wei-Hsuan Yu}
\email[Wei-Hsuan Yu]{whyu@math.ncu.edu.tw}
\thanks{All authors are at the Department of Mathematics, National Central University, Taoyuan City, Taiwan.
The code, certificates and replay instructions are available at
\url{https://github.com/Wideep886/Equiangular_Lines_Maxima}.}
\date{25 November 2026}

\begin{document}
\begin{abstract}
In collaboration with OpenAI's GPT-6, we prove that the maximum numbers of real equiangular lines in dimensions 18, 19, and 20 are $N(18)=57$, $N(19)=72$, and $N(20)=90$, respectively. We reduce the upper bounds to paired lattice points in sections of even unimodular lattices of rank twenty-four. Root-system estimates, exact certificates, and Golay-code arguments bound the sections. In the Leech lattice, generalized Hamming weights of a binary evaluation code give sharp bounds attained by nested sections. At prescribed cosine $1/q$, we also show that $N_{1/q}(d)=d+7$ for explicit intervals of consecutive dimensions for each odd $q\ge319$, and $N_{1/q}(d)=d+23$ on such intervals for each odd $q\ge5599$. Regular graphs in repeated Seidel blocks vary the cosine while preserving nullity, and a coherence bound proves optimality. Symmetric conference matrices give a further sharp family at quadratic irrational cosines.
\end{abstract}
\maketitle

\section{Introduction}
A set of lines through the origin in \(\R^d\) is equiangular if the angle between every two distinct lines is the same. Write \(N(d)\) for the maximum cardinality of such a set over all choices of the common angle, and \(N_\alpha(d)\) for the corresponding maximum with common cosine \(\alpha\).

\begin{theorem}[Main results]\label{thm:combined}
The maximum numbers of real equiangular lines in dimensions eighteen, nineteen, and twenty are
\[
N(18)=57,\qquad N(19)=72,\qquad N(20)=90.
\]
Every set of \(57\) equiangular lines in \(\R^{18}\) spans \(\R^{18}\) and has common angle \(\arccos(1/5)\).
\end{theorem}

Sections~\ref{part:18}, \ref{part:19}, and~\ref{part:20} prove the three equalities using mathematical reductions and exact finite calculations. In dimension nineteen, the dual unit-vector and marked-section exclusions in Theorems~\ref{d19:thm:A} and~\ref{d19:thm:B} rule out a system of 73 lines.

The semidefinite computations of Barg and Yu~\cite{BY}, with subsequent analytic bounds of Yu~\cite{Yu2017}, give \(N(d)=276\) for \(24\le d\le41\) and \(N(43)=344\). The table of Greaves and Syatriadi~\cite[Table~1]{GS} left precisely dimensions eighteen, nineteen, and twenty undetermined among \(d\le41\). The three equalities in Theorem~\ref{thm:combined} determine the remaining entries through dimension forty-one.

The previously published upper bounds in dimensions \(18,19,20\) were \(59,74,94\), respectively \cite{GS,GSY1}. Xiyue Zhou reports the improved bound \(N(19)\le73\) in \emph{On the maximum number of real equiangular lines in dimension nineteen}, preprint, 8 September 2026, \href{https://doi.org/10.13140/RG.2.2.25857.62560}{Theorem~1.1}. The constructions giving the lower bounds \(57,72,90\) are described or cited at the start of the corresponding dimension arguments.

The shared reduction realizes each line as a pair of points in an even integral lattice. It applies to any odd reciprocal cosine, and the rank inequality in Lemma~\ref{lem:paired-sections} gives the available root markings after embedding into rank twenty-four. The final exclusions use the resulting lattice structure. In dimension twenty the Leech branch relates the number of axes in a section to generalized Hamming weights and gives a single chain of sharp sections of the 276 axes in dimensions twenty-two, twenty-one, and twenty.

The same 276-line configuration also serves as a starting point for the fixed-angle results in Section~\ref{fa:section}. Repeating the Seidel blocks and adding a regular graph within each block varies the odd denominator while preserving the dimension of the kernel. A common projection permits the addition of one line and one dimension at a time. Theorem~\ref{fa:intervals} combines this construction with a coherence bound to determine \(N_{1/q}(d)\) on explicit intervals. For instance, at cosine \(1/19\) the upper bound \(d+7\) holds and is attained in every dimension from \(49\) through \(51\). Theorem~\ref{fa:irrational} gives exact values at a family of quadratic irrational cosines using symmetric conference matrices.

Separate manuscripts dated 25 November 2026 are available for
\href{https://github.com/Wideep886/Equiangular_Lines_Maxima/blob/d7ffd626aca5398c83e840dd670a58bf339643b1/dimension-18/paper/el18_max_20261125.pdf}{dimension eighteen},
\href{https://github.com/Wideep886/Equiangular_Lines_Maxima/blob/80b8514f964c9c729c5b53085e41da22aa5f81cc/dimension-19/paper/el19_max_20261125.pdf}{dimension nineteen},
and \href{https://github.com/Wideep886/Equiangular_Lines_Maxima/blob/d7ffd626aca5398c83e840dd670a58bf339643b1/dimension-20/paper/el20_max_20261125.pdf}{dimension twenty}.
Section~\ref{sec:reproduction} gives the dimension-specific guides, the combined replay command, and the version of the computational supplement used here.

\subsection{Even unimodular extension}\label{sec:completion}
\begin{lemma}\label{lem:completion}
Every positive definite even integral lattice \(E\) of rank at most \(21\) embeds isometrically in a positive definite even unimodular lattice of rank \(24\).
\end{lemma}
\begin{proof}
Start with the even unimodular lattice \(U=E_8^{\oplus3}\). Diagonalize the rational quadratic space of \(E\) and embed its positive diagonal directions successively into \(U\otimes\Q\). At each step the orthogonal complement has dimension at least four. For a required squared norm \(a>0\), adjoining the one-dimensional form \(-at^2\) gives an indefinite rational form of dimension at least five. Meyer's theorem \cite{Cassels} gives a nonzero isotropic vector, whose last coordinate is nonzero by positive definiteness of the complement. Dividing by that coordinate gives a vector of the required norm. Proceed by taking orthogonal complements, as in the rational embedding result of Kapovich, Potyagailo and Vinberg \cite[Proposition~2.3]{KPV}. We may therefore regard \(E\) as a sublattice of \(U\otimes\Q\).

Set \(D=U+E\), a full lattice since the generators of \(E\) have a common denominator in a basis of \(U\). Its dual is
\[
A=D^*=\{u\in U:(u,E)\subseteq\Z\},
\]
since \(U^*=U\). Put \(P=A+E\). Both summands are even and pair integrally with each other, so \(P\) is even integral. If \(y\in P^*\), then \(y\in A^*=D\), so \(y=u+e\) for \(u\in U\) and \(e\in E\). Since \(y\) and \(e\) pair integrally with \(E\), we have \(u\in A\), and hence \(y\in P\). Integrality gives the reverse inclusion, so \(P=P^*\). As \(P\) contains the full lattice \(A\), it has rank \(24\). The embedding need not be primitive.
\end{proof}

\subsection{Paired lattice sections}\label{sec:paired-sections}
\begin{lemma}\label{lem:paired-sections}
Let an equiangular line system span dimension \(r\) and have common cosine \(1/q\), where \(q\ge3\) is odd. For any positive definite even integral lattice \(T\), including the zero lattice, there is an even integral lattice of rank \(r+1+\rk T\) containing \(T\) and a vector \(w\perp T\) with \(w^2=2(q-2)\). Each line yields a pair \(x,w-x\) perpendicular to \(T\) with
\[
x^2=q-1,\qquad (x,w)=q-2.
\]
The pairs lie in a common affine section of dimension \(r\). If \(r+1+\rk T\le21\), the entire construction embeds in an even unimodular lattice of rank twenty-four.
\end{lemma}
\begin{proof}
Choose representatives \(v_i\) with squared length \(q\) and pairings \(\pm1\). Choose \(w\) of the indicated squared length and a copy of \(T\) so that they are orthogonal to each other and to the span of the \(v_i\). Set \(x_i^\pm=w/2\pm v_i/\sqrt2\). These points have the stated norm and pairing with \(w\), and \(x_i^++x_i^-=w\). Vectors from distinct pairs have pairings \((q-2\pm1)/2\). The integral span of \(w\), the \(x_i^+\), and a lattice basis of \(T\) therefore has integral pairings and even squares. It is discrete because its generators pair integrally with a linearly independent subset spanning the same real space. Its rank is \(r+1+\rk T\), and all paired points lie in \(w/2+\operatorname{span}_{\R}\{v_i\}\). The last assertion follows from Lemma~\ref{lem:completion}.
\end{proof}

For \(q=5\) the section has \(w^2=6\), \(x^2=4\), and \((x,w)=3\). In ambient dimensions eighteen, nineteen, and twenty, the rank inequality allows the choices \(T=A_2\), \(A_1\), and \(0\), respectively. The positive-rank root markings in dimensions eighteen and nineteen force the rank-twenty-four lattice to have roots. The zero marking in dimension twenty leaves the Leech lattice among the possibilities.

\begingroup
\small
\setlength{\parskip}{.8pt plus .15pt}
\tableofcontents
\endgroup

\section{Dimension eighteen}\label{part:18}
\subsection{Statement and angle reduction}
A set of lines through the origin in $\R^d$ is \emph{equiangular} if the angle between every two distinct lines is the same. Write $N(d)$ for the maximum size of such a set, with the common angle allowed to vary. In dimension eighteen, Greaves, Syatriadi and Yatsyna constructed $57$ lines~\cite[Theorem~1.1]{GSY2}, and Greaves and Syatriadi proved that there are at most $59$~\cite{GS}. Further constructions of $57$ lines are given in~\cite{LMTY} and in the preprint of Munemasa, Sz\"oll\H{o}si and Yoshino, \emph{Sets of equiangular lines in dimension $18$ constructed from $A_5^3\oplus A_1^4$}, \href{https://arxiv.org/abs/2606.10421}{arXiv:2606.10421} (2026). We determine the exact maximum.

\begin{theorem}\label{d18:t:main}
$N(18)=57$. Every set of $57$ equiangular lines in $\R^{18}$ spans $\R^{18}$ and has common angle $\arccos(1/5)$.
\end{theorem}

To exclude $58$ lines we keep both representatives of each line and translate all $116$ vectors by the same amount. They become lattice points in a single affine subspace. The key is to prove that its direction space contains no root. Sections~\ref{d18:r:section}--\ref{d18:sec:rootcases} establish the required statement about short vectors. Section~\ref{d18:s:height} then constructs the ambient lattice and a linear functional that vanishes on the selected points. Sections~\ref{d18:rb:section} and~\ref{d18:hc:section} bound their number.

\begin{lemma}\label{d18:a:angle}
Any set of at least $57$ equiangular lines in $\R^{18}$ spans the space and has common cosine $1/5$.
\end{lemma}
\begin{proof}
Choose $57$ lines and one unit vector on each. Since $57>2\cdot18$, the common cosine is $1/q$ for an odd integer $q\ge3$~\cite[Theorem~3.4]{LS}. At cosine $1/3$ there are at most $34$ lines in dimension $18$~\cite[Theorem~4.5]{LS}. If $G$ is the Gram matrix of the chosen vectors and $r=\rank G$, then Cauchy's inequality on the nonzero eigenvalues gives
\[
  57^2=(\tr G)^2\le r\tr G^2=r\cdot57(1+56/q^2).
\]
As $r\le18$, we have $q^2\le336/13<49$. Thus $q=5$, and the same inequality yields $r\ge475/27>17$.
\end{proof}

\subsubsection{A configuration attaining the bound}\label{d18:s:lower}
We include a construction using the blocks of Lin, Munemasa, Taniguchi and Yoshino~\cite[Definition~5.2]{LMTY}. It splits the $57$ vectors into $45$ indexed by the
edges of $K_{10}$ and $12$ indexed by the directed edges of $K_4$.
Centered indicators express all three types of inner product as
intersection numbers.

\begin{proposition}\label{d18:l:construct}
There are $57$ lines in $\R^{18}$ with common angle $\arccos(1/5)$.
\end{proposition}
\begin{proof}
For $B\subseteq\{0,\ldots,9\}$ put
\[
 p(B)=\sqrt2\bigl(\mathbf1_B-\tfrac{|B|}{10}\mathbf1\bigr)
 \in W:=\mathbf1^\perp\subset\R^{10}.
 \qquad (p(B),p(C))=2|B\cap C|-|B||C|/5 .
\]
Partition the edges of $K_{10}$ into nine perfect matchings.
Explicitly, for each $i\in\Z/9\Z$, take $\mathcal E_i$ to consist
of $\{i,9\}$ and $\{i-j,i+j\}$ for $1\le j\le4$, with arithmetic
on the other nine labels modulo nine. The midpoint of an edge
avoiding vertex $9$ determines its matching. In $W\oplus W$ take
the $45$ vectors $(p(\{i\}),p(e))$ with $e\in\mathcal E_i$.
For two distinct such vectors their inner product is
$2(\delta_{ii'}+|e\cap e'|)-1\in\{\pm1\}$: edges in one matching
are disjoint, and distinct edges in different matchings meet at
most once.

Identify the ten coordinates of the first factor with the four
vertices of $K_4$ and its six directed Hamilton cycles. For each
arc $i\to j$ use the five-subset of~\cite[Definition~5.2]{LMTY},
\[
 B_{ij}=\bigl(V(K_4)\setminus\{i\}\bigr)
       \sqcup\{C : i\to j\text{ belongs to }C\}.
\]
An arc occurs in two cycles, since the two remaining vertices may be
visited in either order. For two arcs with the same tail their blocks
share three vertices and no cycle. For distinct tails they share two
vertices and at most one cycle: if both arcs can lie in one cycle,
contracting them leaves two directed paths (counting an isolated vertex
as a path), which can be joined in a unique cyclic order. Distinct
blocks thus intersect in two or three points.

For each of the twelve blocks $B$ add $(p(B),0)$. Its inner product
with a matched-edge vector is $2\mathbf1_{\{i\in B\}}-1$, and two
different blocks have inner product $2|B\cap C|-5\in\{\pm1\}$.
All $57$ vectors have square five. Their Gram matrix has rank at least
$(\tr G)^2/\tr(G^2)=285^2/[57(25+56)]>17$. Since $\dim(W\oplus W)=18$,
the common angle is attained in full rank.
\end{proof}

\subsection{From a dual unit vector to a finite root system}\label{d18:r:section}
Suppose for the remainder of the upper-bound argument that $58$ lines exist. By Lemma~\ref{d18:a:angle} choose representatives $v_i$ spanning $V=\R^{18}$ with squared length $5$ and pairings $\pm1$ between distinct vectors. We use $z^2=(z,z)$ throughout. Put $L=\sum_i\Z v_i$ and $L^*=\{u\in V:(u,L)\subseteq\Z\}$. The generators of $L$ pair integrally with a basis extracted from them, so they lie in the dual of the lattice spanned by that basis. Thus $L$ itself is a full lattice.

\begin{theorem}\label{d18:r:unit}
The dual lattice $L^*$ contains no vector of squared length one.
\end{theorem}
We first show what a dual unit vector would force. The classification by root rank in Sections~\ref{d18:sec:puregeom} and~\ref{d18:sec:rootcases} will then yield a contradiction. Let $F=\sum_i v_iv_i^{\mathsf T}$. Directly from the Gram matrix,
\begin{equation}\label{d18:eq:frame}
\tr F=290,\qquad \tr F^2=58(25+57)=4756,\qquad \|F-16I\|_{\mathrm F}^2=84.
\end{equation}
\begin{lemma}\label{d18:r:saturate}
If $L^*$ contains a vector of squared norm $1$, there is an integral lattice $M$ containing $L$ with the following properties. Let $Q$ denote the sublattice generated by all vectors of squared norm $4$ in $M$. Then $Q\ne0$, $M=L+Q$, and
\[
 s^2=4,\quad (s,L+Q)\subseteq2\Z\qquad\Longrightarrow\qquad s\in Q.
\]
Writing $\tau(x)=(x,v_1)\bmod 2\in\{0,1\}$, one has
\[
 (x,y)\equiv\tau(x)\tau(y)\pmod 2,\qquad
 x^2\equiv\tau(x)\pmod 4\qquad(x,y\in M).
\]
The norm-one shell of $M$ is empty.
\end{lemma}

\begin{proof}
Start with $L$ and adjoin $2u$, where $u\in L^*$ and $u^2=1$. Continue by adjoining vectors of squared norm $4$ pairing evenly with the current lattice. Each lattice lies between $L$ and $L^*$. The quotient $L^*/L$ is finite since $L$ is integral and has full rank. Choose an endpoint of maximal index over $L$. The asserted saturation and generation properties follow directly.

For $x=\sum_i c_iv_i+\sum_j k_jr_j$, each adjoined $r_j$ pairs evenly with every preceding generator. All inner products of original generators are odd, and $v_i^2\equiv 1\pmod 4$. Consequently $x^2\equiv (\sum_i c_i)^2\pmod4$ and $(x,v_1)\equiv\sum_i c_i\pmod2$. Expanding $(x,y)$ gives the stated parity formula. If $x^2=1$, each $(x,v_i)$ is odd, whence $x^{\mathsf T}Fx\ge58$. But~\eqref{d18:eq:frame} gives $\lambda_{\max}(F)\le16+\sqrt{84}<26$.
\end{proof}

\subsubsection{Weights and their residual vectors}

Write $W=\R Q$, $r=\dim W$ and $c=18-r$. If $s\in M$ has $s^2=4$, then $s$ pairs evenly with $M$. Its inner product with an original vector belongs to $\{0,\pm2\}$, since $|(s,v)|=4$ would give $(v\mp s)^2=1$. Thus reflection in $s^\perp$ preserves $M$, and the norm-four vectors divided by $\sqrt2$ form a simply laced root system~\cite[Chapter~VI, \S4 and the plates]{Bourbaki}. Write $Q=\sqrt2\,Q_0$, where $Q_0$ is its root lattice in the convention that roots have squared length $2$. For $v=p+z$ with $p\in W$ and $z\in W^\perp$, $p/\sqrt2$ pairs with every root of $Q_0$ in $\{0,\pm1\}$. On each irreducible factor it is zero or a minuscule weight. The Weyl group acts on $M$, so these orbit calculations apply even when the original line system has trivial symmetry group.

Let $a$ be the class of $p/\sqrt2$ modulo $Q_0$ in its weight lattice. Write $w(a)$ for twice the minimum squared length of a representative. For a direct sum the function $w$ is additive. The component data, including zero, are
\[
\begin{array}{c|c|l}
 & Q_0^*/Q_0 & w \\ \hline
 A_d&\Z/(d+1)&2j(d+1-j)/(d+1),\quad 0\le j\le d\\
 D_d,\ d\text{ odd}&\Z/4&0,\ d/2,\ 2,\ d/2\\
 D_d,\ d\text{ even}&(\Z/2)^2&0,\ 2,\ d/2,\ d/2\\
 E_6&\Z/3&0,\ 8/3,\ 8/3\\
 E_7&\Z/2&0,\ 3\\
 E_8&0&0 .
\end{array}
\]
Indeed, a shortest representative pairs with every root in $\{0,\pm1\}$ by subtracting a root whenever this shortens the representative. In the dominant chamber the nonzero possibilities are the fundamental weights at nodes of coefficient one in the highest root. Their squared lengths are the corresponding diagonal entries of the inverse Cartan matrix. Each dominant minuscule weight represents a distinct nonzero class in the component table. After Weyl conjugation this also shows that every weight pairing in $\{0,\pm1\}$ with all roots is shortest in its class. In particular every original vector satisfies $w(a)+z^2=5$.

\begin{lemma}\label{d18:r:address}
Consider two original vectors with data $(a,z)$ and $(b,z')$ modulo simultaneous negation of a class and its residual vector. If the two classes with residuals are distinct up to this sign, both extrema
\[
 \frac{w(a)+w(b)-w(a-b)}2+(z,z'),\qquad
 \frac{w(a+b)-w(a)-w(b)}2+(z,z')
\]
belong to $\{-1,1\}$. Moreover, for any integral linear combination of realized classes with residuals $(d,t)$, the number $w(d)+t^2$ is an integer congruent to $0$ or $1$ modulo $4$. For $(d,t)\ne(0,0)$ it is at least $5$.
\end{lemma}

\begin{proof}
All Weyl images belong to $M$ and agree with the original vector modulo $Q$. For two such images an inner product of absolute value $5$ implies proportionality. An inner product of absolute value $3$ gives a sum or difference of squared norm $4$, hence in $Q$. Either would identify the two classes with residuals up to sign. Odd parity excludes the remaining values other than $\pm1$. For the upper extremum take dominant minuscule representatives of $a,b$. Their difference pairs with each positive root in $\{0,\pm1\}$, hence is shortest in its class. Applying this to $a,-b$ gives the lower extremum. Finally $M$ contains each translate of a realized vector by $Q$. Taking a shortest representative of its projected class realizes exactly $w(d)+t^2$ in $M$. Apply Lemma~\ref{d18:r:saturate} and the definition of $Q$.
\end{proof}

If two realized labels have the same residual after orientation, their
difference has zero residual and even parity. Thus
\begin{equation}\label{d18:r:separation}
 a\ne b\quad\Longrightarrow\quad w(a-b)\in4\Z,\qquad w(a-b)\ge8 .
\end{equation}
Indeed, its shortest representative lies in $M$, and every vector of
square four already belongs to $Q$.

\subsubsection{Bounds for the multiplicities}
For later use, two elementary multiplicity bounds are recorded here. In a $d$-dimensional subspace the number $n$ of original lines obeys
\begin{equation}\label{d18:r:bound}
 (5n)^2\le d(n^2+24n),\qquad n\le B(d):=\left\lfloor\frac{24d}{25-d}\right\rfloor .
\end{equation}
This is the trace inequality for their Gram matrix. If $n$ oriented original vectors have the same residual $z$ of squared length $h$, and their projections lie in a $d$-space, then
\[
 (h-1)n\le4,\qquad n\big((5-h)^2-d(1-h^2)\big)\le24d .
\]
Indeed $n^2h\le\|\sum_i v_i\|^2\le n^2+4n$. For the second inequality the projected Gram $P$ has trace $n(5-h)$ and $\operatorname{tr}(P^2)\le24n+n^2(1-h^2)$, obtained by expanding $P_{ij}=(v_i,v_j)-h$ and using the lower bound on the squared sum.

For every real symmetric operator $D$, centering $F$ and $D$ at their scalar parts yields
\begin{equation}\label{d18:r:cut}
 \big(18\operatorname{tr}(DF)-290\operatorname{tr}D\big)^2
 \le1508\big(18\operatorname{tr}(D^2)-(\operatorname{tr}D)^2\big).
\end{equation}
Here $\|F-(290/18)I\|_{\rm F}^2=754/9$. Choosing $D$ constant on specified blocks of $W$ bounds linear forms in the multiplicities of the weight classes. The same argument allows off-diagonal coefficients on the residual space.
\subsubsection{The finite search domain}
Fix one root type of rank $r$.
For each class $a\in Q_0^*/Q_0$ retain the finite set of shortest representatives in that class, scaled by $\sqrt2$. It is one orbit of the reflection group, and every member has square $w(a)$. Choose an integer $S>0$ clearing all values of $w$. For $58$ oriented lines assign classes $a_i$, projected representatives $p_i$, and residual squared lengths $h_i=5-w(a_i)\ge0$. If $H_{ij}=(z_i,z_j)$, the diagonal of $H$ is fixed and the possible off-diagonal entries belong to the following finite sets:
\[
  y_{ij}=2S H_{ij}\equiv-S\{w(a_i+a_j)+h_i+h_j\}\pmod{4S},
  \qquad y_{ij}^2\le 4S^2h_ih_j .
\]
Indeed the sum of the two realized vectors has even parity. Subtracting an element of $Q$ to minimize its projected length preserves the square modulo four and gives precisely $w(a_i+a_j)+h_i+h_j+2H_{ij}$. The inequality is Cauchy's inequality. A realized choice satisfies
\[
 H\succeq0,\qquad \operatorname{rank}H\le c,\qquad
 (p_i,p_j)+H_{ij}\in\{-1,1\}\quad(i\ne j).
\]
These conditions specify a finite domain containing every realization, including repeated classes, coincident residuals and $h_i=0$. The certificates discard families within this domain before all of its entries are assigned.

When $c>0$, the residuals span $W^\perp$, where $W=\R Q$, since the original vectors span $V$. Thus some $h_i>0$. Taking each possible occurring class with positive $h_i$ in turn as a reference covers every system. For $c=2$ its residual has a linearly independent partner. For that fixed reference, the positive numbers $4S^2h_i h_j-y_{ij}^2$ all belong to one rational square class: choose the anchor and an independent partner as a basis, with Gram matrix $G_0$. Every residual has rational coordinates in this basis, since $G_0$ and all its pairings with the basis are rational. For any independent pair with coordinate matrix $U$, its Gram determinant is $\det(G_0)\det(U)^2$. The common scale $4S^2$ is itself a square. Partners for which the determinant vanishes are included in every such branch. At $c=0$ all residuals vanish, and an occurring projected class alone serves as the reference. If the remaining projected and residual choices have multiplicities $m_j$, then $m_j\ge0$ and $\sum_jm_j=58$. Given necessary bounds $Cm\le b$, any $u\ge0$ with
\[
  58\min_j (u^{\mathsf T}C)_j>u^{\mathsf T}b
\]
excludes the entire domain, since the left side is at most $u^{\mathsf T}Cm$. All quantities in this final cut can be checked after clearing denominators.

\subsubsection{Coverage by root rank}
The input consists of all unordered sums of irreducible $A,D,E$ components of ranks at most $18$. Generate a type by adjoining a component no earlier than its predecessor in a fixed ordering until the prescribed rank is reached. Permutations of equal factors and diagram automorphisms are simultaneous changes of coordinates on the entire problem. They require no symmetry of the chosen lines. An \emph{anchor} specifies a label with its residual position that is assumed to occur. A system can occur for more than one anchor. Table~\ref{d18:table:root-coverage} summarizes the cases by rank.
\begin{table}[htbp]
\centering\small
\caption{Root-type coverage and the final necessary test by rank.}\label{d18:table:root-coverage}
\begin{tabular}{@{}rlp{0.56\linewidth}@{}}\toprule
Rank & ADE types & Last necessary test\\\midrule
$1$--$5$ & $1,2,3,6,9$ & Root incidence leaves $A_5,D_5$; bounds $39,35$.\\
$6$--$15$ & $1794$ & Pure and small types handled separately; the mixed residual test ends with five code cases on three types.\\
$16$ & $815$ & $21$ anchor records on nine types; the component-support bounds are at most $57$.\\
$17$ & $1137$ & Two terminal types with bounds $56,53$.\\
$18$ & $1599$ & Fourteen types; eleven admit integer coordinates and three have capacity at most $57$.\\
\bottomrule\end{tabular}
\end{table}

\subsection{Orthogonal roots and their support systems}
\label{d18:sec:puregeom}
Suppose that the complete set of vectors of square norm four in the ambient integral lattice $M$ is $\{\pm2e_i:1\le i\le r\}$, where the $e_i$ are orthonormal. The parity and reflection properties follow from Lemma~\ref{d18:r:saturate} and Section~\ref{d18:r:section}. Write $W=\langle e_1,\ldots,e_r\rangle_{\R}$ and write each of the 58 original vectors as $v=p+z$, with $p\in W$ and $z\perp W$. Then $p_i\in\{0,\pm1\}$. Put $I(v)=\{i:p_i\ne0\}$. Independent changes of signs of the $p_i$ preserve $M$.

\begin{lemma}\label{d18:lem:orth-support}
The possible support sizes are $0,1,2,3,5$. Distinct supports have at most one common point. Their residual inner product is zero if they meet and has absolute value one if they are disjoint. Each support of size five meets every other distinct support. A support of size two, three or five carries at most two, four or five original lines, respectively.
\end{lemma}
\begin{proof}
For $|I|=4$, the sum of $v$ and its image under reversal of all four coordinates is $2z$, a root outside $W$. For $I\ne J$, the inner product of vectors in the two reflection orbits is odd. Equality in Cauchy's inequality excludes absolute value five. Absolute value three gives, by adding or subtracting the two vectors, a root with a coordinate $\pm1$ on $I\mathbin{\triangle}J$. Thus all such inner products are $\pm1$, whereas their range has width $2|I\cap J|$ and midpoint $\langle z,z'\rangle$. A support of size five has $z=0$.

On the same support of size two or three, residuals on distinct axes would again forbid absolute values three and five, despite width at least four. Orient the lines to have the same $z$. Their sign words have mutual Hamming distances two or three. In length three a complementary pair allows no further word, and in its absence the words have the same parity. In length five orient the words to have even weight. They have mutual distance two. After translating by one of the words, the remaining supports are pairwise intersecting two-subsets of a five-set. They form a star or a triangle. Finally, on a one-point support orient all lines to have $p_i=1$. Distinct residuals cannot be equal or opposite, and sign changes then force their inner product to vanish.
\end{proof}

Let $F=\sum_vvv^{\mathsf T}$ and $H=F-16I$. On axis $i$ orient each of the $t_i$ active vectors as $e_i+u$. The residuals $u$ have norm squared four and mutual inner products $0,-2$. Let $b_i$ count the pairs attaining $-2$. Then
\[
 \|He_i\|^2=(t_i-16)^2+4t_i-4b_i=(t_i-14)^2+60-4b_i.
\]
Within this section write $B=\sum_i b_i$, $Q=\sum_i t_i$ and $Z=\sum_vzp^{\mathsf T}$. The sum over the remaining columns of $H$ contains $\|Z\|_F^2$ and the square norm of the complementary principal block. The latter has trace $2+16r-Q$. Hence
\begin{equation}\label{d18:eq:pureenergy}
84\ge 60r-4B+\sum_i(t_i-14)^2+
  \frac{(2+16r-Q)^2}{18-r}+\|Z\|_F^2.
\end{equation}
At $r=18$ one has $Q=290$ and the complementary-block term is omitted.

For a pair $v,w$, set $d=v-\langle v,w\rangle w$. Its square norm is eight. The pair contributes to $b_i$ precisely when $|d_i|=2$. In that case $d-\operatorname{sgn}(d_i)2e_i$ is a root orthogonal to $e_i$, hence equals $\pm2e_j$ for some $j\ne i$. The pair contributes at exactly two axes and lies on one support. Distinct supports have orthogonal contributions to $Z$ in the Frobenius inner product by Lemma~\ref{d18:lem:orth-support}. The following table gives the contribution to $\|Z\|_F^2-4B$ of a support with $a$ lines. For size one it is $4a$ for every allowed $a$. For sizes two and three it follows by evaluating $(5-|I|)\|\sum p\|^2$ and counting sign-word pairs at distance two. For size five all pairs contribute.
\[
\begin{array}{c|rrrrrr}
 |I|\backslash a&0&1&2&3&4&5\\\hline
 2&0&6&-8&&&\\
 3&0&6&0&-18&-48&\\
 5&0&0&-8&-24&-48&-80
\end{array}
\]

\begin{proposition}\label{d18:p:pure}
The complete root system of $M$ cannot have type $A_1^r$ for any $1\le r\le18$.
\end{proposition}
\begin{proof}
Let $q$ denote the number of distinct supports of size five. Any $m\le q$ such supports use at least $5m-\binom m2$ axes. Dropping $Z$ and completing the squares in~\eqref{d18:eq:pureenergy} gives $4B\ge60r-84+(2+2r)^2/18$. The contributing pairs form disjoint cliques of size at most five, and $B$ is their total degree. A clique of size five comes from a full five-point support. Thus $B\le232$, which rules out $r\ge16$. For $r=14$, at most four supports of size five give $B\le3\cdot58+4\cdot5=194<403/2$. For $r=15$, the relative bound permits at most $36$ lines in $W$, so $B\le3\cdot58+36=210<1964/9$.

For $r=13$ one needs $B\ge185$, which forces three full five-point supports. On the remaining 43 vertices only ten cliques of size four and one of size three supply the required degree of at least 125. Thus $B=186$. At least nine axes occur in just one of the three five-point supports, and at most three of these are used by the last triangle. At least six axes therefore have $t_i\equiv1\pmod4$. Completing the squares in~\eqref{d18:eq:pureenergy} gives $84\ge36+28^2/18+6(13/9)^2$, a contradiction. At the other extreme, $b_i\le2(r-1)$ because the negative pairs have distinct oriented sums $u+u'$, all roots orthogonal to $e_i$. Indeed, a sum pairs positively with its endpoints and nonpositively with all other $u$. For $r=1$ we have $b_i=0$, and the two copies of the off-diagonal block give $84\ge(t_i-12)^2+112$, a contradiction. For $2\le r\le7$ substitution gives $84\ge r(68-8r)+(2+2r)^2/18$, again impossible.

For $8\le r\le12$, write $E$ for the right-hand side of
\eqref{d18:eq:pureenergy}. The preceding contribution table gives a lower
bound $12-8h$ per line on each support of size $h=1,2,3$.
A size-two support carrying just one line exceeds this lower bound
by ten. We use the identity
\[
           (t-14)^2-8t=(t-18)^2-128 .
\]

\emph{No five-point support.}
For $r=12$ all lines contribute at least $-12$ each and $Q\le174$.
Hence $E\ge720-12\cdot58+20^2/6>84$.
For $r=8$ let $o,d,t$ be the numbers of lines on supports of sizes
one, two, three. Each of the eight axes belongs to at most three
distinct triples, so there are at most eight triples and $t\le32$.
Counting coordinate pairs gives $d/2+3t/4\le28$.
Cauchy's inequality bounds the sum of the squared deviations in
\eqref{d18:eq:pureenergy}, including its complementary block, below by
$(2+2r)^2/18=18$. Thus
\[
  84\ge480+4o-4d-12t+18\ge82+4o+6(32-t).
\]
It follows that $t=32$, $o=0$, and then that $d=8$.
But now $Q=112$ and the complementary block alone contributes
$18^2/10$, giving $E\ge64+18^2/10>84$.
For $r=9,10,11$ the binary residual Gram matrix restricts the
possible support partitions. Sections~\ref{d18:app:pure-nine},
\ref{d18:app:pure-ten} and~\ref{d18:app:pure-eleven} use these partitions,
the energy inequality and the residual dimension to exclude the
three cases.

\emph{One five-point support.} Let $a\le5$ be its multiplicity, $s=r-5$, $c=18-r$, and let $o,d,t$
be the other line counts of sizes one, two, three. Every such support
meets the five-set once. The activity totals on and off the five-set are
$58+4a$ and $U=d+2t$. The contribution table and the block trace
inequality yield
\begin{align*}
 E\ge L_r(a,U)={}&60r+232-4a^2-8U+(4a-12)^2/5\\
     &+(U-14s)^2/s+(16r-56-4a-U)^2/c .
\end{align*}
Here the non-quintet supports contribute at least
$4o-4d-12t=4(58-a)-8U$. Completing the square in the $s$-block
gives $E\ge872-68r-4a^2$. This is at least $92$ for $r\le10$.

At $r=12$, if there were at most twelve triples, then
$t\le48$ and $U\le106-a$. Use $a^2\le5a$, $x^2/5\ge2x-5$
and $y^2/6\ge4y-24$ on $x=4a-12$ and $y=136-4a-U$.
This gives $E\ge1443-28a-12U\ge171-16a\ge91$.
Thus there are at least thirteen triples.
At $r=11$, at most ten triples would give $t\le40$ and $U\le98-a$.
Using $a^2\le5a$ and the tangents to the three squared terms of $L_{11}$
at $(a,U)=(5,93)$ gives
\begin{equation}\label{d18:eq:q1tangent}
 E\ge813/10+(76/5)(5-a)+7(93-U).
\end{equation}
For $a\le4$ this is at least $813/10+(41/5)(5-a)>84$.
For $a=5$ it requires $U=93$, which forces $t=40$, $d=13$.
The odd value of $d$ requires a size-two support carrying just one
line, improving the bound by ten. Hence there are at least eleven
triples. They use distinct pairs among the six axes outside the
five-set, so there are at most fifteen.

For $r=11,q=1$, each triple support consists of a colour $c$ in
the quintet and an edge $e$ on the other six points. Each colour
is a matching and an edge occurs at most once. Choose one residual
of square two per triple. For distinct triples their residuals are
orthogonal if they share a colour or vertex, and otherwise pair to
$\pm1$. Write $E_j$ for the even-weight binary vectors of length $j$,
with the dot product. The binary Gram matrix of these residuals
is represented by the vectors
\[
 (\mathbf1+e_c,\mathbf1_e)\quad\hbox{in}\quad
 E_5\perp (E_6/\langle\mathbf1\rangle).
\]
Indeed their pairing is $1+[c=c']+|e\cap e'|\pmod2$, including
the diagonal. The displayed space is nondegenerate of dimension eight;
the Gram matrix has binary rank at most the real residual dimension
seven. Thus the displayed vectors have a nonzero common normal
$(\alpha,\mathbf1_A)$. Replacing $A$ by its complement, take
$|A|=0$ or $2$. Each selected edge satisfies
$\alpha_c=|e\cap A|\bmod2$.

There are at least eleven edges by the trace bound above. For $A$
empty at most three colours can be used, giving at most nine edges.
Hence $A=\{p,p'\}$; write $B$ for the other four vertices.
Only colours with $\alpha_c=1$ occur across this cut. There are
$0,2$ or $4$ such colours. Zero gives at most the seven internal
edges. Two allow at most four crossing edges and seven internal
edges. Four allow eight crossing edges and at most three edges of
the one internal colour. Thus there are exactly eleven edges, with
equality in every relevant bound.

We show that no size-two support can occur. It would consist of
a colour $c$ and a vertex unused by the matching of colour $c$.
Its residual $z$ would have squared norm three.

With two crossing colours all edges of the $K_4$ on $B$ occur.
Its three internal colours partition those six edges into three
perfect matchings. The six associated residuals $g_i$ are mutually
orthogonal: any two edges intersect or have the same colour.
Use $g_i/\sqrt2$ as the first six orthonormal coordinates.
A crossing residual has three pairings of absolute value one with
the $g_i$, and hence last coordinate $\pm1/\sqrt2$.
For the proposed $z$ the number $m$ of such pairings is four if its
colour is internal (its vertex must then lie in $A$), and three if
its colour is crossing (its vertex must lie in $B$).
Its last coordinate would be $\pm\sqrt{(6-m)/2}$. Its pairing with
a crossing residual would therefore have the form
$(j\pm\sqrt{6-m})/2$ with $j\in\Z$, contradicting integrality for $m=3,4$.

With four crossing colours all eight crossing edges occur.
The internal colour is a perfect matching and contains $pp'$.
The four crossing residuals at each vertex of $A$ are orthogonal.
Their combined Gram matrix has diagonal blocks $2I_4$ and
cross-block $C$ with two entries $\pm1$ in each row and column.
The other two entries are zero, from the coincident colour and the
coincident vertex in $B$. The residual for $pp'$ is orthogonal to
these eight residuals, so their Gram matrix has rank at most six.
It follows by a Schur complement that $C$ has singular value two
at least twice. The support graph of $C$ is a union of even cycles.
On each cycle, $\|Cx\|^2\le4\|x\|^2$ follows by
$(a+b)^2\le2(a^2+b^2)$ in every row; equality determines all
column entries of $x$ on that cycle up to one scalar.
Thus there are two four-cycles, and in each the signs can be switched
to make the two-by-two block all positive. Its four residuals can be
written in orthonormal coordinates as $f_1\pm f_2,f_1\pm f_3$.
A vector pairing integrally with them has all three coordinates
integral or all three in $\Z+1/2$.

The proposed size-two support must use a crossing colour and a
vertex in $B$, so it is disjoint from the triple for $pp'$.
The squared projection of $z$ on the latter axis is $1/2$.
The two three-dimensional blocks account for the remaining $5/2$.
To have fractional part $1/2$ both blocks must have half-integral
coordinates. But each such block has squared norm congruent to
$3/4$ modulo two, giving $3/2$ modulo two in total. This excludes $z$.

With no size-two support, at most $44$ lines use triples, so
$U\le88$. Equation~\eqref{d18:eq:q1tangent} gives
$E\ge813/10+35>84$, completing the rank-eleven case.

For $r=12,q=1$ take thirteen distinct triples as established above and let $T=\sum zz^{\mathsf T}$ for their residuals. The five axes of the quintet occur 13 times and the other seven occur 26 times. Consequently at least $11+36=47$ pairs of triples intersect. We have $\operatorname{tr}T=26$ and $26^2/6\le\operatorname{tr}T^2\le114$, with even integral middle member. Thus $\operatorname{tr}T^2=114$. Put $P=T-4I$ on the six-dimensional residual space. The residual lattice has full rank here, since $26^2>5\cdot114$, is even integral, and is preserved by $P$. We have $\tr P=\tr(P^2)=2$. If exactly $k$ eigenvalues $\mu_i$ are nonzero, their product is a nonzero integer up to sign, being a coefficient of the characteristic polynomial in a lattice basis. Hence $1\le\prod_i\mu_i^2\le(2/k)^k$. It follows that $k\le2$, and the trace identities force two eigenvalues equal to one. Thus $P$ is an orthogonal projection. Each residual $z$ splits into two orthogonal lattice vectors $Pz$ and $z-Pz$; one is zero since $z^2=2$ and both squares are even. As $\tr(PT)=10$, five residuals lie in its image. For any one of these, $10=(Tz,z)\le4+4=8$, since the other four pair with it in absolute value at most one. This contradiction excludes $r=12,q=1$.

\emph{Three or more five-point supports.} Four quintets require at least fourteen axes and three require
at least twelve. At $r=12$ three quintets have three distinct
pairwise intersections and nine points occurring once. Write
$n\le15$ for the number of lines on the three quintets. Every
other line uses either one pairwise intersection and one single
point, or three single points. If $d$ lines use the first pattern,
the activity totals on the two blocks are $2n+d$ and $174-2d$.
Since each quintet multiplicity $a$ satisfies $a^2\le5a$, we obtain
\[
 E\ge-120+12n+(2n+d-54)^2/3+(12-2d)^2/9 .
\]
The inequalities $x^2/3\ge-8x-48$ and $y^2/9\ge-4y-36$ give
$E\ge180-4n\ge120$. This completes the trace reductions.

\emph{Two five-point supports.} The supports meet at one point, and their two differences each have
size four. Denote their multiplicities by $a,b$, put $n=a+b\le10$,
and let $U$ count all lines using the common point. Thus $U\ge n\ge2$.
Every other line avoiding this point meets both differences once;
the activity totals on them are $4a+58-U$ and $4b+58-U$.
Since every line meets a quintet, the general per-line bound and the
quintet contributions give
\[
 E\ge696-68r+32n-4(a^2+b^2)+\sum_i(t_i-18)^2 .
\]
Keep just the common point and the two differences in this sum.
Cauchy's inequality reduces the right-hand side to at least
\[
 696-68r+(4-2U)n+(U-18)^2+(U+14)^2/2.
\]
Now $(U-18)^2\ge128-8U$ and $(U+14)^2/2\ge28U$. Consequently
\[
 E\ge824-68r+4n+2U(10-n)\ge864-68r .
\]
This exceeds $84$ for $r\le11$, leaving only $r=12$.

For $r=12$ denote the common point by $u$, the two four-point
differences by $A,B$, and the other three axes by $C$.
Besides the quintets the supports have types
$u,uC,uCC,AB,ABC$. Let $h$ be the total activity on $C$ and $m$
the number of lines of type $ABC$. The last two terms omitted in the
previous bound are the deviations on $C$ and the residual block:
\begin{align*}
E\ge {}&-120+(4-2U)n+(U-18)^2+(U+14)^2/2\\
       &+(h-54)^2/3+(78-4n+U-h)^2/6 .
\end{align*}
At most one $uCC$ support can occur, since any two would meet in
two axes. If there are at most nine distinct $ABC$ supports,
then $m\le36$ and $h\le m+(U-n)+4\le40+U-n$.
Bound the last two squares by their tangents at
$h-54=-8$ and $78-4n+U-h=8$, respectively. Substituting the upper
bound on $h$ and writing $x=U-16$, $k=10-n$ yields
\[
 E\ge86+\tfrac32x^2+\tfrac23x+(2U-\tfrac43)k\ge86.
\]
Here $x$ is integral, $k\ge0$ and $U\ge n\ge2$.

Thus at least ten distinct cells of $A\times B$ carry an $ABC$
support, coloured by its member of $C$. Each colour is a matching,
and each cell occurs only once. A colour with four cells exists;
label its matching as the diagonal. Its four residuals $g_i$
have square two and are mutually orthogonal. An off-diagonal
residual belonging to cell $(i,j)$ is orthogonal to $g_i,g_j$
and pairs to $\pm1$ with the other two. Its component $f$ in
the orthogonal complement of the four $g_i$ therefore has square
one. This complement has dimension at most two.

The six pairs $\{i,j\}$ split into the three perfect matchings
of $K_4$. In the same matching, two off-diagonal residuals
have four-dimensional projections with integral pairing:
the two nonzero coordinates have identical or disjoint support.
Between different matchings they have a half-odd-integral pairing,
since their supports meet once. The original residual pairings
are integral. Thus within one matching the $f$-axes are equal or
orthogonal; across two matchings their absolute pairing is $1/2$.
At least six off-diagonal cells occur and each matching accounts
for at most four. Two matchings therefore give unit vectors $f_1,f_2$
at angle $60^\circ$ after orientation. Every occupied matching has
just one $f$-axis: a unit vector from another matching cannot have
absolute pairing $1/2$ with two orthogonal unit vectors in this plane.

Consider now the residual $z$ of any remaining support and its
plane projection $z'$. The four pairings $(z,g_i)$ are integers,
and for an off-diagonal residual $x_j$ projecting to $f_j$ we have
\[
  2(z',f_j)=2(z,x_j)-\sum_{i=1}^4(z,g_i)(x_j,g_i)\in\Z .
\]
Writing these two integers
as $b_1,b_2$ and inverting the two-by-two Gram matrix gives
\[
                 3\|z'\|^2=b_1^2-b_1b_2+b_2^2 .
\]
The norm on the right is even only when both $b_i$ are even, in
which case it is divisible by four. In particular $\|z'\|^2=2$
is impossible.

A support $AB$ must be off the diagonal, hence exactly two of its
four pairings with the $g_i$ vanish, giving $\|z'\|^2=3-1=2$.
For a support $u,uC$ or $uCC$ the four pairings are either all
zero or all $\pm1$. Its four-dimensional projection then has
integral pairing with every off-diagonal projection, so both
$b_i$ are even and $3\|z'\|^2$ is divisible by four. The possible
plane squares are $2$ for $u$, $1$ or $3$ for $uC$, and $0$ or
$2$ for $uCC$. Only $uCC$ with square zero survives; its two
colours must avoid the chosen diagonal colour.
If there are ten cells the total line count is at most
$10+4\cdot10+4=54$. If there are eleven cells, the colour sizes are $4,4,3$. Repeating
the argument with either perfect colour as the diagonal rules out
every $uCC$. The bound is then
$10+4\cdot11=54$.

For twelve cells all three colours are perfect. Among the other
eleven cells, precisely $2+2+3=7$ share a row, column or colour
with a given cell. For one residual per cell put $T=\sum zz^{\mathsf T}$.
Then $\operatorname{tr}T=24$ and
$\operatorname{tr}T^2=12(4+4)=96$. Equality in the rank-six trace
bound gives $T=4I_6$. On the complementary plane this requires
$\sum ff^{\mathsf T}=4I_2$ from the eight off-diagonal cells.
The three possible projective directions there are separated by
$60^\circ$. Their traceless projectors $ff^{\mathsf T}-I_2/2$
sum to zero and any two are independent, as is seen in coordinates
at angles $0,60,120$ degrees. All three multiplicities must therefore
be equal, impossible since $3\nmid8$.
\end{proof}

\subsubsection{Nine orthogonal roots without a five-point support}
\label{d18:app:pure-nine}
Here $o,d,t$ count lines on one-, two- and three-point supports,
$e=58-o-d-t$ counts empty supports, and $m$ counts distinct triples.
Set $Q=o+2d+3t$. Averaging the activity squares over all nine axes
in~\eqref{d18:eq:pureenergy} gives
\[
 84\ge E\ge540+4o-4d-12t+\frac{200}{9}
                       +\frac{2}{9}(Q-136)^2 .
\]
A two-point support carrying a single line improves this by ten.
Completing the square gives $E\ge884/9-12e+2(Q-154)^2/9$, hence
$e\ge2$. Taking the tangent to the last square in the displayed
inequality at $Q=281/2$ gives $E\ge5143/18+6o-6t$,
so $t\ge34$ and $m\ge9$. We shall also use the tangent at $Q=132$:
\begin{equation}\label{d18:r9:double}
  9E\ge7204+20o-68d-156t .
\end{equation}
In particular, if $m=9$ then $d\ge13$, since $t\le36$ and
$d\le12$ would give $9E\ge772>9\cdot84$.

\paragraph{A partition of the nine coordinates.}
For every line with support $I$ the binary vector
$(1,\mathbf1_I)\in\F_2^{10}$ realizes its residual Gram pairings
modulo two, including its diagonal. The binary Gram rank is at
most the real residual dimension nine. As an empty support occurs,
these vectors include $(1,0)$. If the $\mathbf1_I$ spanned
$\F_2^9$, the displayed vectors would span the nondegenerate space
$\F_2^{10}$, giving Gram rank ten. Thus there is a nonempty
$A$ meeting every support evenly. Write $B$ for its complement and
$a=|A|$, $b=9-a$. A triple lies in $B$ or consists of a pair in
$A$ coloured by a point of $B$. Each colour is a matching, and no
pair is repeated. The number of coloured pairs is at most
$\min\{\binom a2,b\lfloor a/2\rfloor\}$. The triples within $B$
number at most $b\lfloor(b-1)/2\rfloor/3$, by counting incidences
at its points; on five points there are at most two since three
linear triples use at least six points. These bounds (and zero
when $b<3$) give, for $a=1,\ldots,9$, the respective upper bounds
\[
               m\le8,8,7,8,9,10,6,4,0.
\]
On four points the incidence bound is rounded down to one.
Every two-point support is an unused pair internal to $A$ or to $B$,
and singletons lie in $B$.
If $a=5$, attaining $m=9$ requires eight coloured pairs and one
triple in $B$. Only two pairs in $A$ and three in $B$ remain,
giving $d\le10$, contrary to~\eqref{d18:r9:double}.
Hence $a=6$. There are either eight coloured pairs with a triple
in $B$, or nine coloured pairs with or without that triple.

\paragraph{Eight coloured pairs.}
The matching sizes are $3,3,2$. The triple in $B$ uses all three
pairs of $B$, so $d\le14$. The displayed energy inequality already
gives $E>84$ for $t\le35$. At $d=13$ at least one two-point support
carries a single line, and the extra ten gives $E>84$ even at
$t=36$. By~\eqref{d18:r9:double} we are left with $t=36$ and $d=14$.
All nine triples have four lines and all seven unused pairs of $A$
have two. Four points of $A$ lie on three coloured edges, and two
on two, giving activities $16,16,16,16,14,14$. On $B$ the
activities are $16+i,16+j,12+k$ with $i+j+k=o$.
The exact support contributions in~\eqref{d18:eq:pureenergy} give
\[
 E\ge80+k^2+8(i+j)+(10-o)^2/9 .
\]
If $i+j\ge1$ this is at least $88$; otherwise it is at least
$80+k^2+(10-k)^2/9\ge90$.

\paragraph{Three perfect matchings.}
Choose one residual of square two from each coloured triple.
In the graph joining pairs of nonorthogonal residuals there are
no edges within a colour. Between two colours the graph is a
perfect matching: each pair on $A$ misses exactly one pair of the
other matching. Hence the graph is a union of cycles whose colours
repeat in order, with lengths divisible by three. Each cycle of
length $l$ has Gram diagonal two and edge entries $\pm1$. By changes
of signs its roots are realized as $f_i-f_{i+1}$ for $1\le i<l$,
and $f_l\pm f_1$ in orthonormal coordinates. The minus case has a
relation involving all $l$ roots with coefficients $\pm1$, so
$l$ must be even: an empty-support residual $y$ pairs to $\pm1$
with all of them.

If every cycle is nonsingular, their nine coordinates span the
residual space. In each cycle the coordinates of $y$ have successive
odd differences and the first and last have odd sum. Thus twice
each coordinate is congruent to $l$ modulo two. An odd cycle
contributes $l/4$ modulo one to $y^2$, and an even cycle contributes
an integer. The total length is nine, so $y^2$ cannot be the integer
five. The only remaining possibility is a singular six-cycle and a
nonsingular triangle. Their spans are orthogonal of ranks five and
three. The triangle Gram can be switched to $I+J$.

If there is a triple on $B$, its residual is orthogonal to all nine
and spans the last dimension. The projections of $y$ onto that line
and the triangle have squares $1/2$ and $3/4$ or $11/4$, respectively,
the latter by $(I+J)^{-1}=I-J/4$. The five-dimensional projection
would then have square $15/4$ or $7/4$. Five consecutive roots of
the six-cycle have the $A_5$ Cartan matrix of determinant six, so
every such projection with integral pairings has square in
$\tfrac16\Z$. This excludes the triple on $B$.

There are now just nine triples. Let $z$ be the residual of a
two-point support. The three triple supports of the triangle
partition the nine coordinates. Our support meets exactly two,
since it cannot reuse a pair contained in a triple. Hence the
triangle projection of $z$ has square $3/4$. In total the support
misses exactly three coloured triples: an $A$-pair misses one per
colour, and a $B$-pair misses all three in the third colour.
Thus $z$ pairs to $\pm1$ at exactly two positions of the six-cycle.
In its difference-root realization the two signs are opposite. If
the circular distance between them is $j=1,2,3$, the projection
has square $j(6-j)/6$: its six coordinates have mean zero and two
values differing by one with multiplicities $j,6-j$. The remaining
one-dimensional projection of $z$ therefore has square, respectively,
\[
                  17/12,\qquad 11/12,\qquad 3/4 .
\]
These three numbers are in distinct rational square classes. All
pairings of the one-dimensional projections are rational, since the
other two projections have rational Gram inverses on integral bases.
Thus the same $j$ occurs for every two-point support. For a $B$-pair
the missed positions have the same colour, hence are opposite on
the six-cycle ($j=3$). An $A$-pair cannot miss both positions of
one colour on the cycle: the only remaining pair of that perfect
matching is already used by a triple. Thus $A$-pairs and $B$-pairs
cannot both occur. There are at most $\max\{\binom62-9,\binom32\}=6$
available pairs, giving $d\le12$, contrary to~\eqref{d18:r9:double}.

\subsubsection{Ten orthogonal roots without a five-point support}
\label{d18:app:pure-ten}

Use $e,o,d,t$ for the numbers of lines on supports of sizes
$0,1,2,3$, and $m$ for the number of distinct triples. Put
$Q=o+2d+3t$. Two consequences of~\eqref{d18:eq:pureenergy} are
\[
\begin{split}
 E&\ge600+4o-4d-12t+22^2/18
   =368+22^2/18+4e+8o-8t,\\
 E&\ge600+4o-4d-12t+(162-Q)^2/8\ge1216-12d-24t .
\end{split}
\]
For the last inequality use $x^2/8\ge4x-32$. In particular
$d+2t\ge95$ and $t\ge39$. At $t=39$ some triple carries fewer
than four lines, which increases the first bound by at least 18.
Thus $t\ge40$ and $m\ge10$.

\paragraph{Parity checks on the supports.}
As in rank nine the augmented vectors $(1,\mathbf1_I)$ realize
the residual Gram matrix modulo two. Their span $U\subseteq\F_2^{11}$
has dimension at most nine: on a subspace of dimension $s$ the
dot product has rank at least $2s-11$, and here the rank is at most
eight. Dropping the first coordinate identifies the even-weight
subspace of $\F_2^{11}$ with $\F_2^{10}$. If every support is odd,
we obtain at least one parity check on the supports. If an even
support occurs, $U$ has an odd-weight vector and its even-weight
part has dimension at most eight. Hence there are two independent
parity checks on all odd supports and on the differences of even
supports.

\paragraph{Odd supports only.}
Here $e=d=0$ and the first energy bound gives $t\ge49$, so $m\ge13$.
Let $A$ be a nonempty check of weight $a$ and $B$ its complement.
Each triple is contained in $B$ or is a pair of $A$ coloured by
a point of $B$. Write $P(s)$ for an upper bound on the number of
linear triples on $s$ points. Counting at each point gives
$P(s)=\lfloor s\lfloor(s-1)/2\rfloor/3\rfloor$ for $s\ge1$,
with $P(0)=0$ and the improvement $P(5)=2$ (three triples use
at least six points). Thus
\[
 m\le P(10-a)+\min\{\tbinom a2,(10-a)\lfloor a/2\rfloor\}
 =12,9,10,10,12,13,10,8,4,0
\]
as $a$ runs from one to ten. Equality at $a=6$ forces four perfect
matchings on $A$ and one triple in $B$. The unused edges of $K_6$
form a fifth perfect matching.

To obtain a rank contradiction, use the ten perfect matchings
outside this factorization as vertices of a graph. Each of the
15 edges of $K_6$ belongs to two such matchings, which it joins.
The graph is the Petersen graph. Indeed, a matching outside the
factorization uses three distinct colours. Its label is the pair
of omitted colours. Two fixed colours cover a six-cycle whose
complement is a triangular prism. The prism has just four perfect
matchings (one uses all three edges between the triangles, the
others use one). Three are the remaining colours, giving exactly
one matching for each label. Two matchings sharing an edge of
colour $i$ cannot both omit a colour $j\ne i$: exactly one
$j$-edge avoids both endpoints, and precisely one of these two
matchings uses it. Thus adjacency means disjoint labels, as claimed.

Delete the three graph edges corresponding to the unused colour:
they form a matching. The remaining graph $H$ is connected. Indeed,
it has ten vertices, minimum degree two and girth at least five
(two intersecting labels have only one common neighbour). If
disconnected, it consists of two five-cycles. Each cycle is
chordless in the original graph, so five edges would have been
deleted. The twelve residuals of the coloured triples have
nonorthogonality graph the line graph of $H$: disjoint edges
of $K_6$ with different colours extend to a unique perfect
matching outside the factorization.

At each vertex of degree three in $H$ the three inner-product
signs around the corresponding triangle have product $+1$.
Otherwise the roots satisfy a dependence with all three
coefficients nonzero. An edge continuing at one of the three other
endpoints meets exactly one member of the star, contradicting the
dependence. Give each incidence at a vertex a sign so the pairwise
products give the prescribed inner products there. Such a choice
is possible at every star by the sign product just proved.
The resulting signed vertex-edge incidence matrix $D$ reproduces
the residual Gram matrix $D^{\mathsf T}D$. Since $H$ is connected,
a vector in $\ker D^{\mathsf T}$ is determined by one entry.
Thus $\operatorname{rank}D\ge9$, exceeding the residual dimension eight.

\paragraph{An even support occurs.}
Give each axis its two-bit label under the checks. Write
$n_0,n_1,n_2,n_3$ for the label multiplicities, taking
$n_1\le n_2\le n_3$; at least two nonzero labels occur.
Every odd support sums to zero, and every even support has the
same sum $k$. The five triple types are $000,011,022,033,123$.
With counts $l_0,l_1,l_2,l_3,h$ they satisfy
\[
 l_0\le P(n_0),\quad l_i\le\min\{\tbinom{n_i}{2},
 n_0\lfloor n_i/2\rfloor\},\quad
 h\le\min_{1\le i<j\le3}n_i n_j .
\]
For $n=(1,3,3,3)$ we can strengthen $h\le9$ to $h\le7$.
Eight or nine transversals form a Latin square of order three
with at most one cell missing (the deficient row and column
both lack the uniquely deficient symbol). A Latin square of order
three is, up to labels, addition modulo three. Its disjointness
graph consists of three triangles, one for each value of the
difference of row and column. Since $m\ge10$, a triple $0ii$
also occurs. Its residual $z$ pairs to $\pm1$ with precisely one
member of each complete triangle and to zero with the other two.
Thus a complete triangle is nonsingular, and switches to Gram
$I+J$. Nine cells give rank nine. At eight cells the two
triangles and the remaining edge already span eight dimensions.
Each triangle contributes $3/4$ to the square of the projection
of $z$, by $(I+J)^{-1}=I-J/4$. The remaining edge contributes
$0$ or $2/3$, by the inverse of its two-by-two Gram matrix.
Their sum cannot be $z^2=2$.

The following table counts the remaining unused pairs. Its entries
are upper bounds for $d$ in the order $k=0,1,2,3$. To reproduce it,
take $n_0+\cdots+n_3=10$, $n_2\ge1$, and retain the partitions
whose five bounds sum to $m\ge10$. For each column start with
twice the number of pairs of that colour and subtract twice the
pairs used by triples. Before doubling, the four values are
\[
 C_0=\sum_j\tbinom{n_j}{2}-3l_0-\sum_{i=1}^3l_i,\qquad
 C_i=n_0 n_i+n_j n_k-2l_i-h\quad(\{i,j,k\}=\{1,2,3\}).
\]
Maximize each expression subject to the five displayed upper
bounds and $l_0+l_1+l_2+l_3+h=m$. This only requires filling
the five types in increasing order of their subtraction coefficient.
\[
\begin{array}{c|rrrr|c}
(n_0,n_1,n_2,n_3)&\multicolumn{4}{c|}{m=10}&m=11\\ \hline
(1,2,3,4)&12&12&6&0&-\\
(1,3,3,3)&12&6&6&6&-\\
(2,2,2,4)&6&12&12&0&-\\
(2,2,3,3)&8&14&8&8&(6,10,4,4)\\
(3,0,1,6)&12&12&6&0&-\\
(3,0,3,4)&0&24&6&0&-\\
(3,1,2,4)&0&18&12&0&-\\
(3,1,3,3)&0&18&6&6&-\\
(3,2,2,3)&0&12&12&6&-
\end{array}
\]

At $m=10$, all triples carry four lines and $e+2o\le2$, hence
$d\ge16$. The table leaves three rows, all with $k=1$. Thus
$e=0$; if $o=1$ then $d=17$, and a singly occupied pair adds
ten to the first energy bound, giving a contradiction. Hence
$d=18$ and $o=0$. The activities summed within the four colours
are respectively
\[
(48,0,42,66),\quad (46,14,28,68),\quad (42,18,48,48).
\]
Indeed all nine available pair supports are full in the last two
rows; in the first, all pairs have colours $2,3$. Distributing
each activity total as evenly as possible among the integer
activities of that colour gives $\sum_i(t_i-14)^2\ge38,42,40$,
respectively. Therefore $E\ge48+36/8+38>84$.

At $m=11$, $d+2t\ge95$ leaves only $k=1$, so $e=0$ and $d\le10$.
If $t\le43$ then $o\ge5$ and the first energy bound exceeds 84.
Thus $t=44$ and $o\ge4$. At $o\ge6$ that bound again exceeds 84;
at $o=5$ add ten for the odd $d$. Hence $(o,d,t)=(4,10,44)$.
All five pair supports are full, giving activities $(28,36,46,46)$
and $\sum_i(t_i-14)^2\ge44$. Now $E\ge48+44>84$.

\subsubsection{Eleven orthogonal roots without a five-point support}
\label{d18:app:pure-eleven}
Let $e,o,d,t$ count the lines on supports of sizes zero, one, two
and three, and let $m$ count distinct triples. Write
$Q=o+2d+3t$. Equation~\eqref{d18:eq:pureenergy} gives
\[
 84\ge E\ge660+4o-4d-12t+(Q-154)^2/11+(178-Q)^2/7 .
\]
The two squares sum to at least $32$. Each two-point support with
only one line improves the displayed bound by ten. A triple with
one, two or three lines improves it by $18,24,18$, respectively.

\paragraph{At least thirteen triples.}
Since $d\le58-o-t$, the bound gives $t-o\ge47$.
Suppose $m\le12$. At $t=47$ one has $o=e=0,d=11$ and an
incomplete triple improves the lower bound by at least $18$.
At $t=48$, the displayed inequality with $32$ in place of the
two squares gives $d\ge8+o$. For $d=8$ this forces $o=0$ and
$Q=160$, where $660-32-576+18^2/7>84$. For $d=9$ there is
a singleton two-point support and $660-36-576+32+10>84$.
The remaining values are $d=10,t=48,e=o=0$. The twelve triples
are full, and if a two-point support is incomplete there must
be at least two such supports, improving the bound by $20$.
Otherwise all activities $t_i$ are even and sum to $164$.
Since $x^2\ge2x$ for even integers $x$, one has
$\sum_i(t_i-14)^2\ge20$. Now $E\ge44+20+14^2/7>84$.

\paragraph{Two parity conditions.}
Over $\F_2$ let $U$ be spanned by the words $(1,\mathbf1_I)$ of
all original supports. Their pairings realize the binary residual
Gram matrix, so the dot product on $U\subset\F_2^{12}$ has rank at
most seven. The ambient form is nondegenerate, whence
$2\dim U-12\le7$ and $\dim U\le9$. Let $D$ consist of the
coordinate words $x$ such that $(0,x)\in U^\perp$.
Thus $\dim D\ge2$. We may select $x,y\in D$ spanning a plane
on which the dot product is singular. Indeed, if $\dim D=2$,
then $\dim U^\perp=3$ and its bilinear rank is at most one:
the restrictions on $U$ and $U^\perp$ have the same radical,
of dimension at least two. If $\dim D\ge3$, take a nonzero
even $x\in D$ and an independent $y\in D$ orthogonal to $x$.

Every nonzero word in $D$ has weight between five and eight.
For a word of weight $b\le4$, every support using its coordinates
contains a pair of them, and each pair is used by at most four
lines. The total activity $T$ on its coordinates is at most
$8\binom b2$. This is inconsistent with
$84\ge(T-16b)^2/b+(2+16b-T)^2/(18-b)$ for $1\le b\le4$.
If $b\ge9$, a triple must use two of these coordinates; each
point outside them therefore colours a matching among the $b$
points. This gives at most $(11-b)\lfloor b/2\rfloor\le8$ triples.

Partition the coordinates by their values on $x,y$. Let $z$ be
the size of the zero part and $a\le b\le c$ the sizes of the
three other parts, denoted $A,B,C$. The three word weights are
$a+b,a+c,b+c$; the singularity of the plane gives
$ab+ac+bc=0$ modulo two. Thus at most one of $a,b,c$ is odd.
The integer conditions $a+b\ge5$, $b+c\le8$, and
$z+a+b+c=11$ give precisely
\[
 (z,a,b,c)=(2,1,4,4),\ (2,2,3,4),\ (1,2,4,4),\ (0,3,4,4).
\]
A triple either uses three points of the zero part, one zero
point and two from a single nonzero part, or one point from
each nonzero part. The first alternative is absent in these
four rows. For each zero point the pairs in any part form a
matching; transversal triples use each pair in $A\times B$ at
most once. The four respective triple capacities are $12,13,13,12$.
Thus only the second and third rows can occur, with $m=13$ and
equality in the capacity bounds of that row.

\paragraph{Part sizes $(1,2,4,4)$.}
Five triples contain the zero coordinate and eight are transversal.
Set $\delta=52-t$, so $0\le\delta\le5$. With every triple full,
the zero coordinate would have activity at least $20$ and the
two-point part activity at least $40$; hence the actual activities
are at least $20-\delta$ and $40-2\delta$. Also
$d\le6+\delta-o$ gives $Q\le168-\delta$.
Keeping the two indicated blocks in~\eqref{d18:eq:pureenergy} yields
\[
 E\ge12+8\delta+3(6-\delta)^2+(10+\delta)^2/7
   =84+\tfrac{22}{7}(\delta-4)^2 .
\]
At $\delta=4$ an incomplete triple strictly improves the inequality.

\paragraph{Part sizes $(2,2,3,4)$.}
The zero points colour one pair of $A$, two pairs of $B$, and four
pairs of $C$. There are also six transversal triples. Call the
colour of the $A$ pair zero. Each colour occurs once on $B$ and
twice on $C$, where the two perfect matchings have distinct edges.
The two $B$ pairs meet. The seven coloured triples thus give two
orthogonal families of residuals, of sizes four and three. Order
the first by $A,B,C,C$ and the second by $B,C,C$. The bipartite
graph of nonzero cross pairings is connected, with seven vertices
and seven edges. Reversing individual residuals makes every sign
on a spanning tree positive. The one remaining sign is $s$ in
the resulting matrix
\[
 H=\begin{pmatrix}1&1&1\\0&1&s\\1&0&0\\1&0&0\end{pmatrix},
 \qquad s=\pm1 .
\]
Positive semidefiniteness requires $4I-H^{\mathsf T}H\succeq0$.
Taking the Schur complement of its first diagonal entry gives
a matrix with diagonal one and off-diagonal $-(2+s)$, forcing
$s=-1$. A realization of these seven vectors, up to their signs, is
\[
 f_4+f_5,\quad f_3+f_6,\quad f_1+f_2,\ f_1-f_2
 \qquad\text{and}\qquad
 f_1+f_5,\quad f_3+f_4,\ f_3-f_4 ,
\]
in orthonormal coordinates. Negating the last column vector gives
the displayed cross pairings for $s=-1$; the two within-family Gram
matrices are $2I$. They span six dimensions.

For any residual $v$, its pairings with these roots are integral.
Consequently the first six coordinates of $v$ are either all
integers or all half-integers. Indeed the two pairs $f_1\pm f_2$
and $f_3\pm f_4$ impose this separately on each pair; the other
three vectors link these coordinates and $f_5,f_6$.
A transversal residual $u$ meets exactly one of the two $C$-pairs
of colour zero. Hence exactly one of $(u,f_1\pm f_2)$ vanishes and
the other is $\pm1$. All six coordinates of $u$ are half-integral.
Since $u^2=2$, each has magnitude $1/2$, and its seventh coordinate
has magnitude $1/\sqrt2$. For every residual $v$ the contribution
of the first six coordinates to $(v,u)$ lies in $\tfrac12\Z$;
thus $v_7=k/\sqrt2$ for some integer $k$. If the first six
coordinates are integral, integrality of $v^2$ makes $k$ even,
and integrality of $(v,u)$ then makes their sum even.
If they are all half-integral, $k$ must be odd, and
$v^2\equiv6/4+1/2=0\pmod2$. In either case $v^2$ is even.
There are therefore no empty or two-point supports.
Now $o=58-t$, $t\le52$, and
$E\ge660+4o-12t+(178-o-3t)^2/7\ge60+16^2/7>84$.

\subsection{Exclusion of the other root systems}\label{d18:sec:rootcases}
It remains to treat the types with at least one component of rank greater than one. The enumerations below sometimes retain pure types through their initial stages. Proposition~\ref{d18:p:pure} accounts for them whenever they are removed. We list the low ranks first, then treat the residual dimensions $18-r\ge3$, $2$, $1$ and $0$.
\subsubsection{Root ranks at most five}
For a unit root direction $e=s/2$, orient the $t$ active lines as $v_i=e+u_i$. Then $u_i^2=4$ and for $i\ne j$, $(u_i,u_j)$ equals $0$ or $-2$. Let $b$ count pairs with value $-2$. With $H=F-16I$,
\[
 84\ge(t-12)^2+112-8b,\qquad
 \|He\|^2=(t-14)^2+60-4b .
\]
The first inequality is the first diagonal block and twice the cross block of $H$ in the decomposition $\R e\oplus e^\perp$. Each of the $b\ge4$ pairs gives the root $u_i+u_j$ orthogonal to $e$. Its inner product with $u_i,u_j$ is $2$ and with any other $u_k$ is nonpositive. It follows that $b\le 2R_\perp(e)$, where $R_\perp(e)$ counts orthogonal root lines. If a component has rank $d_j$ and Coxeter number $h_j$, its unit root lines satisfy $\sum ee^{\mathsf T}=(h_j/2)I$. A fixed line has $2h_j-3$ nonorthogonal root lines including itself. Thus, writing $R$ for the number of all root lines,
\[
 R-2h_j+3\ge2,\qquad
 84\ge\sum_j d_j\max\{0,60-8(R-2h_j+3)\}.
\]
All types up to rank four fail these inequalities. At rank five only $A_5,D_5$ remain. For $A_5$ the weights of subset sizes at least two share one residual direction by Lemma~\ref{d18:r:address}, since their orbit pairings span at least four. The singleton weights have residual norm squared $10/3$. Oriented by singleton, two of these residuals cannot be equal or opposite for distinct lines, as is seen from their projected pairings $5/3,-1/3$. Their mutual inner product is $-2/3$ by the two orbit extrema. There are at most six. Equation~\eqref{d18:r:bound} gives at most $B(6)+6+B(13)=39$ lines. In $D_5$ the vector weights share one residual direction and have equal residual length. On orienting to the same residual, different coordinate indices would give inner product $3$. At most two such lines occur. The spin weights share one residual direction by their orbit pairing width of at least four, for either chirality. The resulting bound is $2+B(6)+B(13)=35$.

\subsubsection{Intermediate root ranks}
For a projected weight class $a$ put $\delta(a)=w(a)-w(2a)/4$. This is additive across the components. For a representative $p$ in the original scale, twice $\delta(a)$ is the difference between $p^2$ and the least inner product $(p,gp)$ on its Weyl orbit, by the orbit-extremum formula in Lemma~\ref{d18:r:address}. If $v=p+z$ with $z\ne0$, then $\delta(a)\le3$. Indeed, $\delta(a)$ is an integer: all Weyl images pair oddly with $v$, and two such pairings differ by $2\delta(a)$ at the extrema. If $\delta(a)\ge4$, a path of simple reflections from $p$ to a minimizing representative passes through an image $gv$ with $(v,gv)=-3$: each step changes the integral pairing by $0$ or $\pm2$. Then $v+gv$ is a square-four root with nonzero residual $2z$.

For each ADE type of rank $r$ and each positive-height anchor, take a common denominator $S$ of the component values $w$. Write $A=S\|z_a\|^2$, $B=S\|z_b\|^2$ and $y=2S(z_a,z_b)$. For an integral combination $ma+nb$ with its specified residual, the minimum squared length, scaled by $S$, is
\[
 S w(ma+nb)+m^2 A+n^2 B+mny .
\]
It is congruent to $(m+n\bmod 2)S$ modulo $4S$ and is zero or at least $5S$. Taking $m=n=1$ makes $y$ integral and fixes its residue modulo $4S$; Cauchy's inequality bounds it by $y^2\le4AB$. Resolving the residual into the anchor direction and its perpendicular complement gives squared energies $y^2$ and $4AB-y^2$ in scale $4AS$. Together with the root-component projections these give the rows tested by~\eqref{d18:r:cut}. Isomorphic factors can be permuted; within an irreducible factor two equal-energy classes are exchanged by diagram automorphisms when necessary (including triality in $D_4$). These changes of coordinates transport the whole configuration and its extension lattice. They need not preserve a chosen extension lattice.

The following count uses all unordered ADE sums, including pure and irreducible types. The first row is the full input; the second counts types admitting some positive-height anchor after the cyclic, pair and trace tests; the last is the exact input to the subsequent support calculation. The additional checker recomputes all $8001$ integer cut inequalities and compares the lists by type.
\[
\begin{array}{c|rrrrrrrrrr|r}
 r&6&7&8&9&10&11&12&13&14&15&\text{total}\\\hline
 \text{all types}&16&24&39&57&88&128&193&276&403&570&1794\\
 \text{positive-height survivors}&15&22&34&48&70&99&141&189&237&242&1097\\
 \text{support inputs}&10&19&32&46&68&98&140&188&236&241&1078 .
\end{array}
\]
The difference of the last two lists consists exactly of the ten pure types $A_1^r$, the irreducible types $A_6,\ldots,A_{10},D_6,D_7$, and the two rank-six types $A_1^4A_2$ and $A_1^2A_2^2$. The pure cases are covered by Proposition~\ref{d18:p:pure}. The root-energy lower bounds for the two rank-six sums are $104$ and $88$, respectively, both greater than $84$.

For $A_d$ write $n=d+1$ and orient each nonzero projected subset to size at most $n/2$. Subsets of size at least two have mutual orbit pairing width at least four, so their residuals share at most one axis. The zero projection gives at most $B(18-d)$ lines. For a singleton projection the residual has squared norm $3+2/n$. Distinct such lines cannot have identical or opposite residuals: with equal residual the total inner product is $3$ or $5$, and with opposite residual it is $-1-4/n$ or $-3-4/n$. The orbit endpoints for two singleton projections then force residual inner product $-1+2/n$. Summing these residuals and taking square norms gives at most $\lfloor4+8/(d-1)\rfloor$ lines of this kind. The resulting bounds $B(d+1)+B(18-d)+\lfloor4+8/(d-1)\rfloor$ for $d=6,\ldots,10$ are $36,34,34,34,33$. In $D_6,D_7$ the vector weights give at most two lines as in $D_5$. Spin weights for both chiralities have mutual pairing width at least four (their Hamming distances range over the admissible parity), so all spin weights share at most one residual axis. This gives $2+B(d+1)+B(18-d)=33,31$, respectively.

The successive totals after support and trace bounds, root-edge energies, sums over components, low-weight residual simplices and two independent residual anchors are respectively $74,37,10,9,7$. In the seven remaining types all factors are $A_1$ or $A_2$. Let $A$ be the matrix of $A_2$ words over $\F_3$, $B$ the matrix of $A_1$ supports over $\F_2$ and $Z$ the real matrix of residual rows. With $J$ the all-ones matrix one has
\[
 3ZZ^{\mathsf T}\equiv-AA^{\mathsf T}\pmod3,\qquad
 3ZZ^{\mathsf T}\equiv J+BB^{\mathsf T}\pmod2 .
\]
These follow directly from the centered singleton coordinates for $A_2$. Both reduced ranks are at most $c$. Short dual words propose further roots. For example, a binary word of weight four orthogonal to all $A_1$ supports gives a sum of four orthogonal unit coordinates pairing evenly with $X$ and $Q$, which would have to belong to $Q$. Likewise a ternary dual word of weight three gives the sum of three shortest $A_2$ weights of squared length $4/3$ in the original scale. Its pairing with every original vector is an even integer by ternary orthogonality. These candidates have a nontrivial class modulo $Q$.

The code tests retain three types, $A_1 A_2^6$, $A_1^2 A_2^6$ and $A_2^7$, with five code cases between them. Every support containing at least two components has a unique residual up to orientation, because for distinct residuals its orbit pairing width would exceed two. The recorded support enumeration and Gram conditions leave $141$ orbits of sets of such labels to test, of which $15$ survive, yielding $421$ assignments of the remaining multiplicities. For each, the residual inner products allowed by Lemma~\ref{d18:r:address} are enumerated and positive semidefiniteness with rank at most $c$ is tested. An independent rational $LDL^{\mathsf T}$ elimination, including the zero-pivot consistency condition, verifies the same terminal counts. The two Gram tests share input support lists. A separate search reconstructs the labelled multi-component supports directly in the five fixed code domains. Of $448\,842$ compatible subsets inspected, respectively $11\,424,30,2,0,128$ satisfy its necessary capacity test. The two enumerations give identical sorted lists, compared entry by entry in all five domains. Enumerating the code-preserving signed permutations also gives the same assignment of every labelled support to its least orbit representative. The five domains give $134,1,1,0,5$ orbits, respectively.

\subsubsection{Root rank sixteen}
Use the finite residual domain of Section~\ref{d18:r:section} with $c=2$. For a fixed positive-height anchor, $A=S z_a^2$, $B=S z_b^2$ and $y=2S(z_a,z_b)$ give the scaled area $4AB-y^2$. An independent partner exists because the residuals span the plane; zero-area partners are included in every positive square-class branch. Rebuilding the discriminant groups from Cartan matrices gives $12\,055$ positive-height anchors and $142\,334$ initial (anchor, area square class) records, before the multiplicity cuts. Across all $815$ types the maximum denominator tested is $210$. The residue representatives are checked against the entire interval $[-10S,10S]$. The successive type counts, after including the separate pure-orthogonal test, are
\[
 815,\ 404,\ 81,\ 50,\ 44,\ 43,\ 20,\ 18,\ 13,\ 10,\ 9.
\]
The nine types give $21$ anchor records. From this frontier, the
weight table and the preceding congruences regenerate $2307$ labelled
word occurrences. The following support counts bound each domain.

\paragraph{Component supports at the rank-sixteen frontier.}
We start from the twenty-one residual domains obtained above. Only their lists of
occurring component shapes are used in the following bounds. The number of domains
and the number of labelled words refer to the input before this argument.
Repeated labels in different domains are counted separately.
Write $U=A_1$, $B=A_2$, $C=A_3$, $D=D_4$ and $E=A_5$.
An absent subscript denotes an end weight in type $A$, or any nonzero weight of $D$;
a subscript denotes a weight up to duality. Factors repeated in a word
denote different components.
\[
\begin{array}{c|rr|l|r}
\text{root system}&\text{domains}&\text{words}&\text{component shapes}&\text{bound}\\\hline
 U^2 B^2 E^2&3&64&E E_2,\ B^2 E,\ U B E,\ U^2 E_3&57\\
 U^2 C^2 D^2&2&91&C D,\ U D^2,\ U C^2,\ U C_2 D,\ U C_2^2&54\\
 U C^5&2&134&C^2,\ C^3,\ U C_2^2&56\\
 B^4 C E&2&80&C_2 E_3,\ B^2 E,\ B^2 C&47\\
 U^3 C^3 D&1&40&U C^2,\ U C_2 D,\ U C_2^2,\ U^3 D,\ U^3 C_2&51\\
 U^7 C^3&3&467&U C^2,\ U C_2^2,\ U^2 C,\ U^3 C_2,\ U^5&54\\
 U^5 C^2 E&3&157&\text{see below}&55\\
U^6 C^2 D&3&430&\text{see below}&56\\
U^{10} C^2&2&844&U C^2,\ U C_2^2,\ U^2 C,\ U^3 C_2,\ U^5&56
\end{array}
\]
For two different realized labels the Weyl pairing interval has width at most two
in the original scale by Lemma~\ref{d18:r:address}. Every common active factor contributes
width at least two: its pairing varies since a nonzero orbit spans the irreducible
factor, and all its values differ by multiples of two by root-weight integrality.
Hence any two different labels share at most one active factor, and any pair of
factors occurs together in at most one label.

Here are the required multiplicity bounds. Within one label we may choose the same
orientation of the residual and projected class, and two different original vectors
then have squared distance eight or twelve. Deleting a factor in which all distinct
weights have squared distance four therefore gives an injective map to the other
coordinates. In particular, the multiplicity of $C^2,C^3,B^2E,UBE,EE_2,
U C^2,U^2 C,CD$ is at most $4,16,9,6,15,8,4,8$, respectively.
The same argument gives $9$ for $B^2 C$.

In $U C_2^2$, $U C_2D$, and $U D^2$ the residual vanishes.
Orient the $U$ coordinate positively. The other coordinates are two signed
orthogonal-axis choices, of respective sizes $(3,3)$, $(3,4)$, or $(4,4)$,
and lengths $\sqrt2$. For $D_4$ this includes both spin orbits by triality
(or by their signed half-coordinate model).
Construct a bipartite graph on the two sets of axes. Occupied pairs are its edges.
Two occupied edges at an endpoint require opposite signs there, so the maximum
degree is two. At the same edge two choices must differ in exactly one sign;
thus there are at most two choices, and doubling an edge requires a degree-one
endpoint at the varying sign. A path on $n$ vertices therefore contributes at
most $(n-1)+\min(n-1,2)\le4n/3$ choices, and a cycle at most $n$. This bounds
the three multiplicities by $8,9,10$.

In $U^3 C_2$ and $U^3 D$ the residual also vanishes. Orient the axis coordinate
positively; at each of its $k=3,4$ axes the other coordinates are three signs.
At the same axis the Hamming distance of two sign strings is two or three;
at different axes it is one or two. Each antipodal pair of strings belongs to
at most one axis, and an axis containing both strings of a pair has no other string.
There are four pairs. If every occupied axis contains a whole pair the total
is at most $2k$; otherwise at most $k-1$ pairs occur doubled. Together with
the four initial pairs this gives multiplicities at most $6,8$ for $k=3,4$.
The ordinary relative bound $B(d)=\lfloor24d/(25-d)\rfloor$ gives $9,11,6$
for $U^2 E_3,C_2 E_3,U^5$, whose residuals vanish and whose ranks are $7,8,5$.

For $U^2 B^2 E^2$, the four multiplicities are multiplied by at most $1,1,4,1$
different labels, respectively, using the pairs $EE,BB,UB,UU$. The total is
$15+9+4\cdot6+9=57$. For $U^2 C^2D^2$ every shape uses a pair of the four
non-$U$ factors. Their six possible pairs have bounds $8$ ($CC$), $9$ (the
four $CD$ pairs), and $10$ ($DD$), giving at most $54$.
For $UC^5$ let $a,b,c$ count the three listed shapes. The ten $CC$ pairs give
$a+3b+c\le10$. Two triples of five $C$ factors must intersect; linearity forces
the intersection to have size one. A third triple meeting both in at most one
point is impossible. Thus $b\le2$. Also $c\le2$, since the $c$ supports
already share $U$ and their $C$ pairs must be disjoint. Hence
$4a+16b+8c\le40+4b+4c\le56$.
For $B^4CE$ the shape $C_2E_3$ occurs at most once. At each of $C,E$ the
other labels realize pairwise disjoint $B$ pairs, hence number at most two.
The bound is $11+4\cdot9=47$.
For $U^3C^3D$, if no label uses all three $U$ factors, every label uses a
pair of $C^3D$; these yield at most $3\cdot8+3\cdot9=51$. A label using
the three $U$ factors is unique and shares a $U$ with every other label, so its
non-$U$ factor cannot occur again. The other labels can then use at most
three pairs at weight at most nine each, plus weight at most eight for the
first label.

Consider the two rows of the form $U^rC^s$ with $r=7,10$. A label using two $C$ factors
has weight at most eight; let $x$ be their number, so $x\le\binom{s}{2}$.
Among the other labels that use $C$, their weights are bounded by twice the
number of occupied $UC$ pairs. Each of the $x$ labels consumes two such pairs.
Excluding $U^5$ labels the total is therefore at most $2rs+4x$.
Two five-subsets in a support system have intersection at most one; on ten
points three such subsets are impossible, since two cover at least nine points
and a third can take at most one point in each plus one outside. Thus $r=10,s=2$
gives $40+4+2\cdot6=56$. On seven points there is at most one $U^5$.
In its absence the bound is $42+12=54$. In its presence the $U^2C$ and
$U^3C_2$ supports attached to a fixed $C$ have disjoint $U$ subsets and each
uses at least one or two of the two $U$ factors outside the five-set, respectively.
Their combined weight at this $C$ is at most eight. The total is then
$3\cdot8+3\cdot8+6=54$.

For the remaining two types use also the residual height $h=z^2$. If $h>1$,
at a fixed oriented residual position the number of lines is at most
$\lfloor4/(h-1)\rfloor$: their projected vectors have squared norm $5-h$ and
pairings at most $1-h$, so positivity of their summed vector gives
$0\le n(5-h)+n(n-1)(1-h)=n[4-n(h-1)]$.
The residual-coordinate lists give the following bounds on the number of
oriented positions at each height, counted up to overall negation.
\[
\begin{array}{c|c|c|c}
\text{root system}&h&\text{shapes}&\text{positions (at most)}\\\hline
U^5C^2E&11/6&CE&2\\
       &5/2&UC&2\\
       &5/6&UCE&2\\
       &0&C_2E_3,\ UC_2^2,\ U^2E_3,\ U^3C_2,\ U^5&1\\\hline
U^6C^2D&3/2&CD,\ U^2C&3\\
       &1&UC^2&2\\
       &0&UC_2D,\ UC_2^2,\ U^3C_2,\ U^3D,\ U^5&1
\end{array}
\]
Distinct labels cannot both use the same $C$ factor at the middle weight,
or the same $E$ factor at weight three: in these self-dual orbits the pairing
width is four or six, respectively. A zero-residual label must meet the
component support of every other label, since otherwise the inner product is zero.

For $U^5 C^2 E$, the $CE$ and $UC$ lines account for at most
$2\cdot4+2\cdot2=12$ lines by the height bound. For $UCE$ the two possible
$CE$ pairs give at most two labels, with multiplicity at most eight by
deleting $E$. The positive total is at most $28$. Among zero-residual labels
other than $U^5$, charge six per $C_2$ factor and nine per $E_3$ factor.
These bound each multiplicity listed above, and each factor is charged at
most once. The total is at most $28+2\cdot6+9+B(5)=55$.

For $U^6 C^2 D$, the $h=3/2$ lines contribute at most $3\cdot8=24$.
The $UC^2$ shape uses the unique $CC$ pair and hence contributes at most
eight, giving at most $32$ positive lines. Denote the five zero-residual
shapes, in table order, by $Q,P,S,T,J$; their respective bounds are
$9,8,6,8,6$. If $J$ occurs, every shape using three $U$ factors meets it
in at least two and hence is excluded. The two available middle $C$ factors
then bound the zero total by $2\cdot9+6=24$. Suppose $J$ is absent.
With no $Q$, at most two $T$ labels fit: their $U$ triples are disjoint since
they already meet in $D$. If two occur their triples partition all six $U$
factors, excluding $S$. This bounds the zero total by $16+8=24$.
With at most one $T$ it is at most $8+2\cdot6=20$.
With one $Q$ there is at most one $T$ (their $U$ supports are disjoint);
only one $C_2$ factor remains for $S$, giving at most $9+8+6=23$.
With two $Q$ labels, their $C$ factors are different, as are their $U$ factors.
At most one $T$ fits, giving at most $26$ zero lines. Every $CD$ label is
excluded since it would meet one of the $Q$ labels in two factors. Every
$U^2 C$ label must meet the $Q$ label on the other $C$ factor at its $U$
factor. Thus all $U^2 C$ labels at a fixed $C$ would share this same $U$
factor as well, and there is at most one. The positive contribution is then
at most $2\cdot4+8=16$. These cases give at most $56$ lines.

\subsubsection{Root rank seventeen}
When $c=1$ the nonzero residual norms squared have one common rational square class $q$, since all residual pairings are rational. The data of a column are its discriminant class and signed height. In all linear combinations the squared residual is included in the minimum from Lemma~\ref{d18:r:address}. For a candidate $s$ of square four with $(s,Q)\subseteq2\Z$ and $s\notin Q$, Lemma~\ref{d18:r:saturate} requires a selected $v$ with
\[
                       (s,v)\notin 2\Z .
\]
At this stage $s$ need not belong to $M^*$, so a nonintegral pairing is possible. For example, in $A_1^3 A_7^2$ take dominant weights with labels $(0,0,0;1,1)$ for $s$ and $(0,0,1;1,7)$ for $v$, and give both residuals the same signed height $\sqrt{2}/2$. These formal columns have $s^2=4$ and $v^2=5$, but $(s,v)=7/4+1/4+1/2=5/2$.
The exact test uses $(s+v)^2=9+2(s,v)$, so even pairing is equivalent to $(s+v)^2\in1+4\Z$. Translation of the projected sum by an element of $Q$ changes its square by a multiple of four because both $s$ and $v$ pair evenly with $Q$. Thus one can test the minimum length in its class and clear all denominators before taking the congruence. Starting with all $1137$ types, the initial tests leave $1059$ pairs (type,$q$), belonging to $320$ types. Pair tests leave thirteen pairs, one for each of thirteen types. Eleven are eliminated by the extra-root and trace tests.

Here are the two terminal square classes directly from the weight table.
A letter $d$ denotes any nonzero $D_4$ class; the $D_5$ symbols
$v,s$ denote its vector and either spin class. The $A$ labels are
listed up to sign, and the two $A_5$ factors of the first type may
be exchanged. The last coordinate $k$ means residual height
$k/\sqrt6$ for the first type and $k/\sqrt{12}$ for the second.
The height convention is exhaustive: $6h$ or $12h$ is an integer
square for the stipulated square class, where $h$ is the residual
squared norm. Taking $w+k^2/6=5$ or $w+k^2/12=5$ gives precisely
the following rows.
\[
\begin{array}{c|ccccc|c|cccc}
&A_3&A_5&A_5&D_4&k&&A_5&A_7&D_5&k\\\hline
 a&0&0&3&d&0 &A&0&0&v&6\\
 b&0&1&2&0&2 &B&0&1&s&3\\
 c&1&0&0&d&3 &C&0&2&v&0\\
 d'&1&1&1&0&1 &D&1&0&v&4\\
 e&2&0&3&0&0 &E&1&2&0&2\\
 &&&&&&F&2&0&v&2\\
 &&&&&&G&3&0&v&0
\end{array}
\]
First discard $C,F$, since twice either row has minimum square $4$.
For the remaining rows the only pair test needed for the preliminary
shape exclusions is that the shortest representative of $x+y$ has
square zero or a multiple of four at least eight. The following
four substitutions give the full details of this test. Retain a
label $j$ in a displayed allowed shape only if
\[
 w(i+j)+(k_i+k_j)^2/m\ \in\ \{0,8,12,16,20\},\qquad m=6,12 .
\]
Here $20$ is sufficient since $(x+y)^2\le20$. Add $A$ labels
modulo the corresponding order, take the sum in the Klein group
for $D_4$, and put $v=2$, $s=\pm1$ modulo four for $D_5$.
Adding the displayed labels with either sign gives the possible shapes
\[
\begin{array}{c|c|l}
\text{fixed row}&\text{previously removed}&\text{possible shapes}\\\hline
 b=(0,1,2,0;2)&-&b,c,a^*,e^*\\
 e=(2,0,3,0;0)&b&c,e\\
 A&C,F&A,E\\
 D&C,F,A&D,E
\end{array}
\]
Here $a^*,e^*$ use their triple in the first $A_5$, and the surviving $b$ in the first row and $e$ in the second use the block order of the corresponding anchor. In the first row every line contributes at least three to
$P_{A_5}+2P_{\rm normal}$, using the first $A_5$. The trace of
this operator is $7$, the trace of its square is $9$, and~\eqref{d18:r:cut}
would give $1102^2\le1508\cdot113$. Thus $b$ is absent. In
the second row the first $A_5$ would be unused, contrary to full
rank of the original system. Hence the first square class only
requires $a,c,d'$. In the third row all $A$ lines have height
squared three, so $(h-1)n\le4$ allows at most two. The $E$ lines
have capacity $24$ as below. In the fourth row the $A_5$ singleton
sign of $D$ is fixed by~\eqref{d18:r:separation}, since the difference
energy is at most $8/3$. Reusing a $D_5$ vector weight on two such
lines would give inner product at least $2+4/3-1/3=3$. There
are ten such weights, so the fourth row has capacity at most
$10+24$. The second square class thus only needs $B,E,G$.

For $A_3A_5^2D_4$ at $q=6$ the preceding cuts allow three shapes,
up to permutation of like factors. For the positive-height shape
using three $A$-singletons, the difference energy for any two labels
is at most $2+8/3+8/3<8$; for the positive-height shape using an
$A_3$ singleton and a nonzero $D_4$ label it is at most $2+2<8$.
At positive height orient along the same normal.
Equation~\eqref{d18:r:separation} allows just one label of each shape.
The balanced zero-height shape occupies one $A_5$ and one $D_4$.
For a fixed $A_5$ its labels differ by energy at most two.
Thus there are at most two balanced labels. For $A_5A_7D_5$ at $q=3$,
the corresponding two positive-height shapes use an $A_5$ singleton
and $A_7$ edge, or an $A_7$ singleton and $D_5$ spin weight. Their
difference energies are at most $8/3+4<8$ and $3+2<8$. The one
balanced label uses the $A_5$ triple and $D_5$ vector class.
This proves directly the label multiplicities used in the bounds
$56$ and $53$, respectively.

For the $A_3A_5^2D_4$ label with three singletons, distinct lines
have inner product $2e-1$ with $e$ equal indices among their three
indices. Its capacity is at most $4\cdot6=24$. The second
positive-height label has eight $D_4$ weights by triality. Two lines
with identical $D_4$ weight would have inner product at least
$2+3/2-1/2=3$, so it contributes at most eight.
For the $A_5A_7D_5$ label with an $A_5$ singleton and an $A_7$ edge,
the inner product is $2(\delta_{ij}+|e\cap f|)-1$. At each of six
singleton indices the edges must form a matching on eight points,
giving at most $24$ lines. For the other positive-height label the
$D_5$ spin weight belongs to one chirality, with $16$ choices.
Reusing a weight would give inner product at least
$5/2+3/4-1/4=3$.

In either zero-height label orient the $D$ weight to the positive
direction of its axis. This gives $t=4$ or $5$ buckets. Within one
bucket the $A_5$ triples intersect in at most one point. There are
at most four triples per bucket: in the absence of a complementary pair their square-three weights pair to $-1$, so their sum has square $n(4-n)\ge0$. A bucket containing complementary triples has size two,
since a further triple cannot meet both in at most one point.
Across buckets the same triple up to complement would produce
inner product $\pm3$. There are ten complementary pairs in total.
If $p$ buckets use both members of a pair, the bound is
$\min(10+p,4t-2p)$, at most $12$ or $13$ according to $t$.

The two totals are $24+8+2\cdot12=56$ and $24+16+13=53$.

\subsubsection{Full root rank}
All residuals vanish. Apply Lemma~\ref{d18:r:address} to integral combinations $ma+nb$ of realized labels and use the block trace bounds~\eqref{d18:r:cut}. The independent verifier starts with all $1599$ types and reconstructs the integers for every discarded anchor. Fourteen types remain. Eleven consist of $A_1$ and even-rank $D$ factors. In standard coordinates the orthogonal change $(x,y)\mapsto((x+y)/\sqrt2,(x-y)/\sqrt2)$ on paired $D$ coordinates sends each allowed weight to an integral vector. This gives an integral realization in $\Z^{18}$. For clarity, the integer-shell test reduces each squared-norm-five vector modulo two, then appends three ones and three zeros. The resulting words of weight four or eight generate a doubly even binary code. The extended words are mutually orthogonal modulo two, since both the original pairing and the contribution of the three appended ones are odd. The all-ones word can be adjoined: each resulting weight $24-j$ is divisible by four when $j$ is. The extension theorem of Huffman and Pless~\cite[Theorem~9.5.5(ii)]{HP} puts this code inside a Type~II code of length $24$. The test uses the classification of the nine Type~II codes of length $24$ by Pless and Sloane~\cite{PlessSloane}.
Let $O,Z$ be the disjoint triples on which ones and zeros are appended.
The possible parity supports in the other $18$ positions are
\[
 \mathcal B(O,Z)=\{w\setminus O:|w|=4\text{ or }8,\ O\subset w,\ w\cap Z=\varnothing\},
\]
where $w$ runs over supports of codewords. All members have size one or five.
Any two of these supports meet in an odd number of coordinates, by
self-orthogonality of the code. Thus a singleton candidate forces every
support to contain its coordinate. This contradicts the frame bound
$\lambda_{\max}(F)<26$, since the frame trace on that coordinate is at
least $58$. Thereafter all vectors have five entries $\pm1$.
On a fixed five-set choose the even-weight sign word in each antipodal
pair and translate one word to zero. The other words have weight two
and their pairs meet. They form a star or a triangle, so each five-set
contributes at most five lines.

Four of the nine codes can now be treated directly.
In the Golay code the $759$ octads form a $5$-$(24,8,1)$ design:
two meet in at most four positions and $759\binom85=\binom{24}5$.
The numbers through three, four and five specified points are $21,5,1$.
By inclusion--exclusion the number through $O$ avoiding $Z$ is
$21-3\cdot5+3\cdot1-\#\{w:O\cup Z\subset w,\ |w|=8\}\le9$.
This gives at most $45$ lines. In the $d_{24}$ code the words of weights
four and eight are unions of coordinate pairs: the other coset
selects one position from each of twelve pairs. Some pair meets $O$ once.
Its other coordinate either belongs to $Z$, leaving no candidates,
or belongs to every candidate.

In either direct sum $e_8^3$ or $d_{16}\oplus e_8$, an octad using
two summands consists of a tetrad in each, since every nonzero word
in a summand has weight at least four. If $O$ lies in one summand,
such an octad yields a singleton candidate from its tetrad in that
summand. Thus all candidates would have to lie in that summand, giving
rank less than $18$. If $O$ meets more than one summand, a surviving octad
must use exactly the two summands met by $O$. In $e_8^3$ this leaves unused coordinates, contradicting
rank $18$. In $d_{16}\oplus e_8$ let $j=1,2$ be the number of marked
ones in the larger summand. Its remaining $d$ coordinates have
$13-j\le d\le16-j$. The trace of $F$ there is $58(4-j)$, because all
selected five-sets meet them in $4-j$ points. For $j=1$, inequality
\eqref{d18:r:cut} already fails at $d=12$, where
$(18\cdot174-290\cdot12)^2=121104>1508\cdot12\cdot6$;
for $j=2$ the endpoint $d=11$ gives
$(18\cdot116-290\cdot11)^2=1214404>1508\cdot11\cdot7$.
In both ranges the discrepancy increases with $d$ and $d(18-d)$ decreases.

The five remaining codes require the following small table. Each column
after the first count removes just the indicated representatives, in order.
The elementary capacity test uses $|\bigcup\mathcal B|=18$ and
$5|\mathcal B|\ge58$. If two positions meet every support, the frame
trace on their span is at least $58$, violating~\eqref{d18:r:cut}.
For the trace column, order positions by their number of incident
supports and use each prefix $S$. The sum of the smallest $58$ of the
numbers $|B\cap S|$, taking each $B\in\mathcal B$ five times,
bounds the actual trace from below; the largest $58$ bound it from above.
\[
\begin{array}{l|r|rrrr|r}
 \text{tetrad system}&\text{cuts}&\text{singleton}&\text{capacity}&
 \text{two positions}&\text{trace}&\text{left}\\\hline
 d_{12}^2&66&8&54&4&0&0\\
 d_{10}e_7^2&195&23&162&8&2&0\\
 d_8^3&129&10&108&10&1&0\\
 d_6^4&138&9&117&7&4&1\\
 d_4^6&34&3&29&1&0&1\\\hline
 &562&53&470&30&7&2
\end{array}
\]
Each entry counts representatives of ordered pairs $(O,Z)$ under the
coordinate permutations used for that code.
For $d_{2k}$ first permute the $k$ coordinate pairs and switch any even
number of pairs; for $e_7$ use invertible linear maps on the seven
nonzero vectors of $\mathbb F_2^3$. List the two disjoint coloured
subsets in each component up to these operations and take products
having three points of each colour in total. Finally interchange
equal components by the supplied permutations preserving the code.
For each row the orbit sizes sum to $\binom{24}3\binom{21}3$.
The two remaining representatives have orbit sizes $32$ and $11\,520$,
respectively.

For the representative in $d_4^6$ (orbit size $11\,520$), join two
supports when they meet in three coordinates. The graph has three
components of each of the following two forms. In the two-support component the supports share
a triple, and their remaining pairs are disjoint. The four sign
patterns of the triple up to inversion cannot be shared by distinct
supports, since this would yield inner product of absolute value
three. At fixed oriented triple the last two signs of distinct
vectors must be opposite. This gives at most eight lines.

The four-support component consists of a common coordinate and
four disjoint coordinate pairs $P_1,P_2,Q_1,Q_2$, with supports
$\{u\}\cup P_i\cup Q_j$. Orient each vector to be positive at
$u$. In each coordinate pair use its two diagonal axes. After
division by $\sqrt2$ the residuals have the form $s e_i+t f_j$
on two disjoint sets of four axes. Their mutual inner products are
$0$ or $-1$. Regard each as a signed edge $ij$. At a common endpoint,
edges to distinct neighbors must have opposite signs there. Thus the
underlying simple graph is a disjoint union of paths and cycles.
Each edge occurs at most twice: two sign pairs on the same edge
must differ at exactly one endpoint. That endpoint cannot have
another neighbor. A cycle on $v$ vertices contributes at most $v$
vectors, an isolated edge at most two, and a path on $v\ge3$ vertices
at most $(v-1)+2=v+1$. At most two vertex-disjoint paths on at least
three vertices fit on eight vertices. The component therefore has
capacity at most ten. Altogether the bound is $3(8+10)=54$.

For the representative in $d_6^4$ (orbit size $32$), the same graph has three components. In each, the supports are the union of a star of $K_4$ (three edge coordinates) and one of three disjoint pairs of coordinates. Orient a vector to be positive on the first coordinate of its pair. Its pair and its sign on the second coordinate determine one of six groups. Within a group the pair contributes two to every inner product. If $k$ distinct stars appear, their signs on the $k-1$ shared coordinates are fixed by opposition; on the $4-k$ remaining coordinates the sign words on each star have distance at least two. Thus a group has size at most four. Across all groups there are only sixteen signed stars up to negation. Two different groups cannot reuse such a star, since their pair contribution is zero. A group using both signs of one star has size exactly two. If $h$ groups do so, the component has at most $\min(16+h,24-2h)\le18$ points. The three components give at most $54$ lines.

\paragraph{The remaining three full-rank types.}
For $A_9^2$ the only shapes are a single five-subset, or a singleton
and a two-subset in different blocks. Use signed subset sizes $(a,b)$
in the two blocks. The inner product of labels $(a,b),(a',b')$ is
congruent to $-(aa'+bb')/5$ modulo $2\Z$, by the centered-indicator
formula. Two mixed labels in opposite block orders would give numerator
$0$ or $\pm4$, impossible. Orient each mixed vector positively on the
singleton block. Opposite signs on the two-subset block give numerator
$1-4=-3$, again impossible. After a common reflection the mixed
inner products are $2(\delta_{ij}+|e\cap f|)-1$; hence each edge is
used at most once, giving $\binom{10}2=45$ mixed lines. Pure lines
from different blocks would be orthogonal. Thirteen pure lines in one
block would give a $13\times10$ balanced sign matrix $U$ with
$UU^{\mathsf T}$ having diagonal $10$ and other entries $\pm2$.
Put $C=U^{\mathsf T}U$. In each row of $C$ the nine off-diagonal
entries are odd and sum to $-13$. Since $(t+1)(t+3)\ge0$ for odd
$t$, their squares sum to at least $52-27=25$. Now
$13\cdot100+13\cdot12\cdot4=1924=\tr(C^2)\ge
10(13^2+25)=1940$, a contradiction. The upper bound is $45+12=57$.

In $A_1A_3A_7^2$, multiply the component energies by four. The
possible entries in the $A_3$ and $A_7$ factors are $0,6,8$ and
$0,7,12,15,16$, respectively. Summing to $20$ gives, up to the
exchange of the two $A_7$ factors, the rows
\[
 (0;6,7,7),\qquad(4;0,16,0),\qquad(0;8,12,0).
\]
Twice a word of the last shape has nonzero class of square four.
In the first shape, each of the three singleton labels has two possible
signs. Changing a subset of signs gives difference energy equal to the
corresponding subsum of $2,3,3$. By~\eqref{d18:r:separation}, two realized
labels of this shape agree up to simultaneous negation. Orient to that
common label. Two vectors now have inner product $2e-1$, where $e$
counts their equal singleton indices. Any two indices determine the
third, giving at most $4\cdot8=32$ lines. All other lines have unit
projection on the $A_1$ axis, giving at most $25$ more.

In $A_1^3A_5^3$ the scaled energies $3w$ are $0,3$ on each $A_1$ and
$0,5,8,9$ on each $A_5$. To sum to $15$ one takes either three
$A_5$ singletons or two $A_1$ axes and one balanced $A_5$ triple.
For two singleton words whose signs differ in $j$ places, the difference
energy is $8j/3$. Equation~\eqref{d18:r:separation} gives $j=0,3$.
As above these lines contribute at most $6^2=36$. For the other
words, reusing the same two axes with a different $A_5$ gives difference
energy six; reusing the same $A_5$ with different pairs gives two.
Hence at most three such labels occur.

For one such label orient the first $A_1$ coordinate positively.
The second sign gives two groups of triples on six points, with
intersection at most one within each group. A group of $m$ triples
has binary distances at least four, whence
$4\binom m2\le6m^2/4$ and $m\le4$.
At equality every distance is four and each of the six coordinates
occurs in exactly two triples. The triples are therefore the four
stars of a $K_4$ whose six edges are the coordinates.
The pairs of opposite edges form a perfect matching on these six
coordinates. The four stars, taken up to complement, are all four
transversal axes of that matching: the eight transversals are the
four stars and their complements. Any two perfect matchings on six
points have a common transversal, by alternating the two colours on
each cycle of their union (including doubled edges). If both groups
had size four, some triples would be equal or complementary across
the groups. Their $A_1$ contributions cancel, leaving inner
product $3$ or $-3$. Each label thus has at most seven lines,
giving $36+3\cdot7=57$.

\begin{proof}[Conclusion of the proof of Theorem~\ref{d18:r:unit}]
A unit vector in $L^*$ produces a lattice $M$ with nonempty root system by Lemma~\ref{d18:r:saturate}. Proposition~\ref{d18:p:pure} and the rank cases just considered exclude every such root system. Hence $L^*$ contains no unit vector.
\end{proof}

\subsection{Lattice points at a common height}\label{d18:s:height}
We now use Theorem~\ref{d18:r:unit} to pass from the original $58$ line representatives to a section of a lattice of rank $24$. In a larger Euclidean space take $w\perp V$ with $w^2=6$ and set
\[
          x_i^\pm=\frac w2\pm\frac{v_i}{\sqrt2}.
\]
The $116$ points have common squared length $6/4+5/2=4$. Two from different pairs have inner product $3/2\pm1/2\in\{1,2\}$. Each pair sums to $w$, so the points lie in $w/2+V$ and are interchanged by $x\mapsto w-x$.

Call a lattice vector of squared length two a \emph{root}. An even unimodular lattice has an integral Gram matrix with even diagonal and determinant one in a lattice basis.

\begin{proposition}\label{d18:p:lift}
There is a positive definite even unimodular lattice $P$ of rank $24$ containing $w$ and all $x_i^\pm$. The space $V$ contains no root of $P$. There are two roots in $P$ perpendicular to both $w$ and $V$ with inner product $-1$.
\end{proposition}
\begin{proof}
Apply Lemma~\ref{lem:paired-sections} with $q=5$ and $T=A_2$. This realizes the paired points and the two simple roots of $T$, with pairing $-1$, in an even unimodular lattice $P$ of rank 24. The same pairing of shell points occurs in Degtyarev's work~\cite[\S2.5]{Deg}.

By Theorem~\ref{d18:r:unit}, no unit vector $u\in V$ pairs integrally with all $v_i$. If $r\in P\cap V$ were a root, then $u=r/\sqrt2$ would be such a vector, since $(u,v_i)=(r,x_i^+)\in\Z$.
\end{proof}

\begin{lemma}[Common height]\label{d18:c:order}
There is a vector $t\perp(\R w+V)$ such that $(t,r)\ne0$ for every root $r\perp w$ of $P$. Every $x_i^\pm$ has height $(t,x_i^\pm)=0$.
\end{lemma}
\begin{proof}
A root perpendicular to $w$ and lying in $\R w+V$ would belong to $V$, contrary to Proposition~\ref{d18:p:lift}. Thus each equation $(t,r)=0$ defines a proper hyperplane in $(\R w+V)^\perp$. Finitely many proper hyperplanes cannot cover that space. Every $x_i^\pm$ lies in $w/2+V$, so its height is zero.
\end{proof}

For example, in the lattice of type $A_{24}$ write $[z]=z-(\sum_i z_i/25)\mathbf1$ in the sum-zero subspace of $\R^{25}$, and abbreviate $[\mathbf1_S]$ to $[S]$. The root lattice has discriminant group $\Z/25\Z$, with coset of $[z]$ given by $\sum_i z_i$ modulo $25$ for $z\in\Z^{25}$. The unimodular extension uses the subgroup of order five. For a square-four vector $x=[z]$, add a constant to the entries of $z$ to make $m=\sum_i z_i$ one of $0,5,10,15,20$. Then $\sum_i z_i^2=4+m^2/25$. Since $\sum_i z_i^2\ge m$, the cases $m=10,15$ cannot occur. In the case $m=0$ the vector has two entries $+1$ and two $-1$. In the cases $m=5,20$, equality $\sum_i z_i^2=m$ makes $z$ a zero-one vector, so $x$ is $[S]$ or $-[S]$ for $|S|=5$.

Now take $w=[T]$ with $|T|=10$. For $m=0$, $(w,x)\le2$. For $x=\pm[S]$, the pairing is $\pm(|S\cap T|-2)$. Thus the candidates with $(w,x)=3$ are exactly the $[S]$ with $S\subset T$ and $|S|=5$. See Degtyarev's \href{https://archive.mpim-bonn.mpg.de/id/eprint/851/1/preprint_2019_51.pdf}{MPIM preprint}, \S5.4, for this case. The ten coordinates $t_i$ for $i\in T$ are distinct because $e_i-e_j$ is a root perpendicular to $w$. Now $(t,[S])=\sum_{i\in S}t_i$. For this sum to vanish, any four positions determine the fifth if it exists. Each eligible five-set contains five four-sets, giving at most $\binom{10}{4}/5=42$ candidates of height zero.

In general, a list along which height increases contains at most one selected point. One can certify a link $x\to y$ by writing $y-x$ as a nonzero nonnegative combination of roots with positive height. A partition of $n$ candidates into lists with $e$ such links therefore bounds the number of selected points by $n-e$. Section~\ref{d18:rb:section} explains how to choose a common orientation for these roots.

\subsection{Height bounds in the Niemeier lattices}\label{d18:rb:section}

Let $T$ be the span of the two roots in Proposition~\ref{d18:p:lift}. Fix $w$ and take candidate points of square four pairing to three with $w$. A selected set is closed under $x\mapsto w-x$. Its nonpartners have inner product one or two. All its points are perpendicular to $T$. The Leech lattice has no roots, and the roots of $A_1^{24}$ on distinct axes are orthogonal. Niemeier classified the $24$ positive definite even unimodular lattices of rank $24$ up to isometry~\cite{Niemeier}. The other $22$ types are covered below.

\subsubsection{Chain bounds in eighteen types}
Here we verify lists of candidates directly. Let $s_r z=z-(z,r)r$ for $r^2=2$. Use a system of roots perpendicular to $w$ whose reflections preserve the candidates. If the mark $T$ is imposed, take roots perpendicular to $T$ as well. Choose a reference vector $h$ with $(h,r)\ne0$ for all these roots; call $r$ positive if $(h,r)>0$. For the height vector $t$ in Lemma~\ref{d18:c:order}, choose $t'$ in its reflection orbit maximizing $(h,t')$. This orbit is finite, since the group permutes the roots and fixes the orthogonal complement of their span. For every positive $r$,
\[
 0\le(h,t')-(h,s_rt')=(h,r)(t',r).
\]
As $(t',r)\ne0$, this puts $t'$ on the positive side of every such root. Applying the same transformation to the selected points preserves their heights. A link $x\longrightarrow y$ is now valid whenever $y-x$ is a nonzero nonnegative combination of positive roots. The height increases, so links with at most one incoming and one outgoing at each point form disjoint paths. On $n$ candidates, $e$ links leave $n-e$ paths, including isolated vertices. To construct a partition, take two copies of the candidates and join the left copy of $x$ to the right copy of $y$ for each allowed link $x\to y$. A matching of size $e$ gives $e$ links with distinct starts and distinct ends. Strictly increasing height rules out cycles. Any maximum matching therefore yields at least as many links as the exhibited certificate; verifying the exhibited links alone suffices for the bound.

\paragraph{Enumerating the candidates.}
Choose simple roots $\alpha_i$ of the indicated root lattice $Q$ with block diagonal
Cartan matrix $C=((\alpha_i,\alpha_j))$. For the column of pairings
$u_i=(z,\alpha_i)$, set $b=C^{-1}u$. Then
\[
 z=\sum_i b_i\alpha_i,\qquad z^2=u^{\mathsf T}C^{-1}u,\qquad
 z\in P\quad\Longleftrightarrow\quad b\bmod\Z^{24}\in H,
\]
where $H$ is the glue subgroup in $Q^*/Q$. Here is a finite way to obtain a model
for each type. On every component enumerate the nonnegative integral columns $u$
of square at most seven, increasing one entry at a time. The inverse Cartan entries
are nonnegative, so no extension of a partial column above the norm bound can
return below it. The enumeration checks that every component coset occurs: its list of residues has
$\det C_i$ entries for the component matrix $C_i$. Since reflections preserve
the residue and bring any weight to a dominant one, this gives the minimum of
every coset. Add these minima over components. Starting with the zero subgroup,
adjoin residues until $|H|^2=\det C$, retaining only subgroups whose nonzero
cosets have even minimum at least four. The resulting extension is even, has
determinant one, and has exactly the roots of $Q$. Thus any one such $H$
models the required lattice by Niemeier's classification~\cite{Niemeier}.

Combine dominant component weights whose class word lies in $H$ and whose
squares sum to six. These are all the dominant $w$. For square four, apply all simple reflections to each dominant point and
keep all points satisfying $(x,w)=3$. The searches retain the component
labels, including the order of isomorphic factors. For an ordinary link
$x\to y$, the certificate checks $b_i(y)=b_i(x)$ on every fixed index and
$b_i(y)\ge b_i(x)$ on the others, with at least one strict inequality.
The allowed indices are precisely the simple roots orthogonal to $w$ and,
when a mark is imposed, to both marking roots. Thus $y-x$ is a positive
combination of the allowed simple roots. Computing the endpoints in the full
square-four list verifies the bound even if some candidates are isolated.

\begin{table}[htbp]
\centering\small
\caption{Choices of the norm-six vector and bounds on signed points in eighteen lattice types.}\label{d18:table:ordinary-height}
\begin{tabular}{@{}lrr@{\qquad}lrr@{}}\toprule
Root type & Choices of $w$ & Bound & Root type & Choices of $w$ & Bound\\\midrule
$D_{24}$&3&42 & $A_{17}E_7$&20&72\\
$D_{16}E_8$&8&59 & $A_{12}^2$&40&80\\
$E_8^3$&10&49 & $A_{11}D_7E_6$&51&83\\
$A_{24}$&9&60 & $E_6^4$&52&80\\
$D_{12}^2$&13&66 & $A_9^2D_6$&83&90\\
$D_{10}E_7^2$&24&73 & $A_8^3$&108&96\\
$A_{15}D_9$&25&75 & $A_6^4$&228&109\\
$D_8^3$&38&78 & $A_7^2D_5^2$&154&101\\
$D_6^4$&94&91 & $D_4^6$&404&112\\\bottomrule
\end{tabular}
\end{table}
The eighteen types in Table~\ref{d18:table:ordinary-height} give $1364$ choices of $w$. Each bound is the largest over the computed cases for its type, in signed points. In $D_4^6$, $218$ of the $404$ choices require the mark. All $2646$ marked refinements give at most $112$ points. To see that the refinements suffice, take the roots perpendicular to $w$. Within their span a generic vector orthogonal to $T$ has exactly the roots of $T$ perpendicular to it: the roots in that plane form $A_2$. Carrying the vector into a closed chamber identifies $T$ with two adjacent simple roots in the system perpendicular to $w$. All such pairs are included.

\subsubsection{Enumeration in the four exceptional types}
Each record specifies $T$ and $w$. Two nonorthogonal roots lie in one component. In an $A$ component their span is equivalent to the differences on the first three coordinates. The same holds in $D_4$ by signed permutation, using the fourth coordinate to adjust the parity of sign changes. A permutation of equal components is used only when it preserves the full lattice. The possible marked components in $A_5^4D_4$ are $A_5$ and $D_4$; they give $131$ and $35$ records. The other three types give $28$, $17$ and $5$ records, as in Table~\ref{d18:tab:marks}.

\paragraph{The marked enumerations.}
For the four exceptional lattices, the following coordinate description gives
the finite domains behind the orbit counts. In the order
$A_5^4D_4,A_4^6,A_3^8,A_2^{12}$, use scales $n=6,5,4,3$. In an
$A_{n-1}$ block the scaled coordinates $a_i=n z_i$ sum to zero and
are congruent to one residue modulo $n$. We label this residue by $-k$ for
$n=6,5$ and by $k$ for $n=4,3$. In a $D_4$ block at scale six,
a label $k=0,1,2,3$ gives $a_i\equiv 0,0,3,3\pmod6$, respectively,
and $\sum_i a_i\equiv6k\pmod {12}$. These labels add by bitwise
addition in $\mathbb F_2^2$. In each lattice combine blocks only for
words of its glue subgroup. Their orders are $72,125,256,729$,
respectively~\cite[Table~16.1]{CS}; each order squares to the determinant
of its root lattice, and each nonzero coset has even minimum at least four.
The component orbits under the checked glue-preserving coordinate isometries have sizes $(4,1),(6),(8),(12)$, respectively. Within one component the numbers of possible marking planes in $A_5,D_4,A_4,A_3,A_2$ are $20,16,10,4,1$ and each set is a single Weyl orbit. These assertions can be checked by listing pairs of roots with pairing $-1$, identifying their spans, and applying the permutations to the glue words. The following descriptions specify the glue subgroups.
\begin{itemize}
\item At scale six take $b\in\mathbb F_2^4$ of even weight, indexed from
zero, and $i,j\in\mathbb F_3$. The first four labels form
$3b+2(i,i+j,i-j,j)\bmod6$.
The fifth label is $b_1\cdot1+b_2\cdot2+b_3\cdot3$ in $\mathbb F_2^2$,
where $1,2$ denote a basis and $3=1+2$.
\item At scale five take all words $((a+bi+ci^2)_{i\in\mathbb F_5},c)$.
\item At scales four and three take the words $(v,vG)$.
At scale four, $v\in(\Z/4\Z)^4$ and
\[
 G=\begin{pmatrix}2&1&1&1\\-1&2&1&-1\\-1&-1&2&1\\-1&1&-1&2\end{pmatrix}.
\]
At scale three, $v\in\mathbb F_3^6$. Index $G$ by
$\infty,1,2,3,4,0$. Its diagonal is zero, the other entries in the
$\infty$ row and column are one, and the other entries are
$G_{ij}=\chi(i-j)$ for the quadratic character $\chi$ of $\mathbb F_5$.
\end{itemize}
A mark in an $A$ or $D$ factor requires $a_1=a_2=a_3$. In $A_2^{12}$ the
whole marked factor therefore vanishes; we shorten the ternary code at
that position. For square six enumerate all component tuples with
$|a_i|\le\lfloor n\sqrt6\rfloor$, the required residues and block conditions,
and $\sum_i a_i^2=6n^2$. Within each factor normalize by Weyl operations
preserving the mark. For the list of remaining shapes, take images of
each proposed representative under coordinate isometries preserving the
mark and the glue code. Check preservation on every codeword, check that
the image sets are disjoint, and remove them from the list of shapes.
Its exhaustion verifies coverage even when the chosen isometries generate
only a subgroup of the full stabilizer. The sizes of the initial lists
and the resulting numbers of representatives are
\[
\begin{array}{c|rrrrr}
  (P,\text{marked factor})&
  (A_5^4D_4,A_5)&(A_5^4D_4,D_4)&(A_4^6,A_4)&(A_3^8,A_3)&(A_2^{12},A_2)\\
  \text{shapes}&401&347&629&1080&1749\\
  \text{representatives}&131&35&28&17&5 .
\end{array}
\]

For each representative form the square-four list by the same rule with
$\sum_i a_i^2=4n^2$, and filter by the scaled pairing $3n^2$ with the
scaled $w$. A link must have difference a positive multiple of a root
perpendicular to $w$ and the mark, using $e_i-e_j$ or, in $D_4$,
$e_i+e_j$, with $i<j$ in a single block. A common reference vector
with positive, strictly decreasing coordinates in each block makes all
these roots positive. After checking link endpoints, this gives
$74\,308$ links on $93\,264$ candidate occurrences. Exactly $201$ of
the $216$ marked records have fewer than $116$ lists. The remaining
fifteen are treated by the binary rank test. Each count is over
marked records, allowing the same line system to have several markings.
\begin{table}[ht]
\caption{Marked sections in the four remaining lattice types.}\label{d18:tab:marks}
\centering\begin{tabular}{lrrr}
\toprule
Lattice type & Markings & Excluded by height & Requiring the rank test\\
\midrule
$A_5^4D_4$&166&165&1\\ $A_4^6$&28&26&2\\
$A_3^8$&17&9&8\\ $A_2^{12}$&5&1&4\\
\midrule
Total&216&201&15\\
\bottomrule
\end{tabular}
\end{table}

\subsection{Binary rank bounds for the remaining marked sections}\label{d18:hc:section}

For a fixed marking $(T,w)$, take one vector $y=2x-w$ from each pair $\{x,w-x\}$. The partners give opposite $y$, and $y^2=10$. If $58$ pairs are selected, their directions span a space of dimension at most $18$, and different selected directions satisfy $|(y,y')|=2$.

Each choice of the marked root span and $w$ has a record in its host model. We retain the original record numbers, starting at zero in each model. Isometric sections in different models may therefore undergo the same test more than once. The preceding height bounds leave eight records in $A_3^8$, two in $A_4^6$, four in $A_2^{12}$, and one in $A_5^4D_4$.

\subsubsection{The bound in a binary section}
Fix the positive roots as in the height argument. Two directions are incompatible if $|(y,y')|\ne2$ or if $y+y'$ or $y-y'$ is a nonzero vector in the cone generated by these roots or its negative. Indeed the height functional vanishes on all selected directions and has strict sign on either cone. For the test it suffices to use the difference roots of equally valued coordinates of $w$ in the same $A$ factor, fixing the three coordinates of the mark individually whenever it lies in that factor. The marked $D_4$ factor is fixed pointwise. Order the coordinates in each free block. A vector lies in the cone precisely when it vanishes on the fixed coordinates, the sum on each free block is zero, and all initial partial sums in each ordered free block are nonnegative. Indeed, if $z_1+\cdots+z_m=0$ in a block, then $(t,z)=\sum_{j<m}(t_j-t_{j+1})(z_1+\cdots+z_j)$. Order the block so $t_j>t_{j+1}$. For the incompatibility test at least one coefficient must be positive.

Keep all candidate axes. As cells take the connected components of the graph joining pairs with forbidden absolute inner product. Within each cell put an edge between directions passing both compatibility tests. Every selected subset within such a component is a clique. Disregarding conditions between components still gives an upper bound. The largest component in the fifteen records has $36$ axes.

Clear the coordinate denominators of the centered directions and form an integer matrix with these vectors as rows. Keep $d$ independent columns spanning its column space. If a nonempty collection of these columns sums to an even column, replace one member by half that sum. This change is invertible over $\mathbb Q$ and halves the absolute value of any fixed nonzero $d$-row minor. Thus the process terminates with $d$ columns independent modulo two. Let $z_y\in\mathbb F_2^d$ be the resulting row. The binary rank of selected rows is at most their rational rank, so their common annihilator contains a space $H$ of dimension $k=d-18$.

For $U\subseteq H$, let $D_U$ be the axes with $z_y\perp U$, put $n_U=|D_U|$, and let $b(U)$ be the sum, over the cells, of the largest compatible subset in the cell restricted to $D_U$. Every selected axis lies in $D_U$, hence $b(U)\ge58$. Each cell bound is obtained by choosing a vertex and taking the larger of the two answers obtained by omitting it or taking it and retaining only its neighbors.

\subsubsection{Reduction to pairs of equations}
Write $n$ for the total number of axes tested and $n_a=n_{\langle a\rangle}$. An axis outside $D_H$ is deleted by exactly $2^{k-1}$ nonzero normals of $H$. Averaging first in $H$ and then in $H/\langle a\rangle$ inside $D_{\langle a\rangle}$ gives independent $a,b\in H$ with
\[
 (2^k-1)(n-n_a)\le 2^{k-1}(n-58),\qquad
 (2^{k-1}-1)(n_a-n_{\langle a,b\rangle})
      \le 2^{k-2}(n_a-58).
\]
Here $k=2$ or $3$. For each nonzero $a$ compute $n_a$ and $b(\langle a\rangle)$. Test all $a,b$ satisfying both inequalities for which $a,b,a+b$ each have single-equation bound at least $58$. The identity $2n_{\langle a,b\rangle}=n_a+n_b+n_{a+b}-n$ follows by checking the four binary pairings with $a,b$. Finally compute $b(\langle a,b\rangle)$ on the actual intersection. Table~\ref{d18:table:binary-planes} lists the number of pairs reaching this last test and the number with bound at least $58$. For each fixed $a$ we retain one of $b,a+b$. Different first normals may lead to the same plane.

\begin{table}[htbp]
\centering\small
\caption{Terminal binary plane tests with a distinguished first normal.}\label{d18:table:binary-planes}
\begin{tabular}{llrrr}
\toprule
Host & Indices (in order) & $d$ & Pairs tested (in order) & Passing\\
\midrule
$A_3^8$ & $0,1$ &20& $0,0$ &0\\
$A_3^8$ & $5,6,10,14,15,16$ &21& $0,12,0,0,0,36$ &0\\
$A_4^6$ & $26,27$ &21& $0,0$ &0\\
$A_5^4D_4$ & $32$ (mark in $D_4$) &21& $3$ &0\\
$A_2^{12}$ & $1,2,3,4$ &21& $3,36,0,472614$ &0\\
\bottomrule
\end{tabular}
\end{table}
The direct enumeration makes $472\,704$ tests with distinguished first normal. No tested plane has bound at least $58$. This excludes all fifteen marked records.

\begin{proof}[Proof of Theorem~\ref{d18:t:main}]
A system of $58$ lines would yield $116$ points as in Proposition~\ref{d18:p:lift}. Section~\ref{d18:rb:section} excludes two types by their roots and eighteen by the height bounds. In the four remaining types the marking enumeration reduces $216$ cases to fifteen, all excluded in Section~\ref{d18:hc:section}. Thus $N(18)\le57$. Proposition~\ref{d18:l:construct} attains the bound. Lemma~\ref{d18:a:angle} gives the assertion on the span and angle of every extremal system.
\end{proof}

\subsection{Verification of the finite computations}\label{d18:sec:verification}
The root stage checks the bounds and rank cases in Sections~\ref{d18:r:section}--\ref{d18:sec:rootcases}. The height stage reconstructs the lattice and marking domains in Section~\ref{d18:rb:section} and verifies the links in each list. The binary stage reconstructs the matrices of Section~\ref{d18:hc:section} and counts the compatible rows. The exact input data, integer witnesses and verification programs accompany the article. Use Python 3.11 or newer with NumPy, Node.js, and a C++17 compiler. After installing the dependencies listed in the reproduction guide, run the following commands from the repository root
\begin{verbatim}
cd dimension-18/reproduction
unset PYTHONOPTIMIZE
shasum -a 256 -c MANIFEST.sha256 &&
    python3 -m zipfile -e EL18_58_Search_Checkpoint.zip . &&
    python3 -B audit/referee/replay.py EL18_58_search --second
\end{verbatim}
The additional check rebuilds the terminal graphs and binary kernels independently and counts the same sections by a Walsh transform. Both implementations use the lattice models of Section~\ref{d18:rb:section}.

\section{Dimension nineteen}\label{part:19}
\subsection{Statements and reduction}\label{d19:sec:statement}
Let $N(d)$ be the largest number of equiangular lines through the origin in $\R^d$. The construction of Asche gives $N(19)\ge72$; see \cite[p.~124]{Taylor} and \cite[Example~2]{LY}. The upper bound $N(19)\le74$ was obtained in \cite{GSY1}. Zhou reports the improved upper bound $N(19)\le73$ in the preprint cited in the introduction. We prove $N(19)=72$ by excluding a system of 73 lines, using the following two theorems. Throughout, a line system is represented by exactly one vector on each line unless explicitly stated otherwise.
\begin{theorem}[A: dual unit-vector exclusion]\label{d19:thm:A}
Suppose $v_1,\ldots,v_{73}\in\R^{19}$ satisfy $\sq{v_i}=5$ and $\ip{v_i}{v_j}\in\{-1,1\}$ for $i\ne j$. Then the dual of $L=\sum_i\Z v_i$ contains no vector of squared norm one.
\end{theorem}
\begin{theorem}[B: marked-section exclusion]\label{d19:thm:B}
Let $P$ be a positive definite even unimodular lattice of rank 24, and let $w,r_0\in P$ satisfy $w^2=6$, $r_0^2=2$, and $w\cdot r_0=0$. Set
\[
\mathcal F(w,r_0)=\{x\in P:x^2=4,\ x\cdot w=3,\ x\cdot r_0=0\}.
\]
If $\mathcal X\subseteq\mathcal F(w,r_0)$ satisfies $w-\mathcal X=\mathcal X$, and $V=\operatorname{span}_{\R}(\mathcal X-w/2)$ has dimension at most 19 and contains no norm-two vector of $P$, then $|\mathcal X|\le144$.
\end{theorem}
Theorem~\ref{d19:thm:A} is reduced below to explicit rational certificates and four exceptional root types. Theorem~\ref{d19:thm:B} uses Golay-code arguments and finite shell enumerations, including an independent computation for all twenty types with at most six irreducible root components. The conclusion from these two statements is as follows.
\begin{corollary}[Main conclusion]\label{d19:cor:main}
$N(19)=72$.
\end{corollary}
\subsubsection{Angle and dimension}\label{d19:sec:angle}
The standard integrality theorem of Neumann \cite[Theorem~3.4]{LS} implies that 73 equiangular lines in $\R^{19}$ have common cosine $1/q$ for an odd integer $q\ge3$. The trace inequality \cite[Lemma~6.1]{vanLintSeidel} gives $73\le19(1+72/q^2)$, whence $q^2\le76/3$. At cosine $1/3$ there are at most 36 lines in dimension 19 by \cite[Theorem~4.5]{LS}, so $q=5$. Scale the vectors to squared norm 5. Writing $F=\sum_i v_iv_i^{\mathsf T}$, we have
\begin{equation}\label{d19:eq:frame}
\tr F=365,\qquad \tr F^2=7081,\qquad \|F-19I\|_F^2=70.
\end{equation}
The rank is at least $365^2/7081=1825/97>18$; thus the vectors span $\R^{19}$.
For later use let
\begin{equation}\label{d19:eq:bound}
B(d)=\min\left\{\left\lfloor\frac{24d}{25-d}\right\rfloor,\ b_d\right\},\quad 0\le d\le19,
\end{equation}
where $b_{14}=28,b_{16}=40,b_{17}=48,b_{18}=59$, and $b_d=+\infty$ otherwise. The first term follows by the trace inequality at cosine $1/5$; the improvements follow from \cite{GSY1,GSY2,GS}. In particular, these arguments do not assume that $N(18)=57$.
\begin{lemma}[Projection support]\label{d19:lem:projection}
If $\pi:V\to U$ is linear and $\rk\pi=s$, then at least $73-B(19-s)$ of the lines have nonzero projection. In particular the lower bounds for $s=1,2$ are 14 and 25.
\end{lemma}
\begin{proof}The lines in the kernel lie in dimension at most $19-s$. Apply \eqref{d19:eq:bound}.\end{proof}
\subsection{Saturation and the root deficit inequality}
\subsubsection{Construction of a saturated lattice}
Assume the contrary in Theorem~\ref{d19:thm:A}: $u\in L^*$ and $u^2=1$. Since the Gram entries are integral and the vectors span, $L$ is a full lattice (it is contained in the dual of the lattice generated by a basis of the $v_i$). Starting with $L+\Z(2u)$, repeatedly adjoin a vector $s$ of square 4 pairing evenly with the current lattice, whenever it is not already contained therein. All intermediate lattices are integral and lie between $L$ and $L^*$, so this process terminates. Let $M$ denote an endpoint, $R=\{s\in M:s^2=4\}$, $Q=\sum_{s\in R}\Z s$, and $W=\R Q$. Then $M=L+Q$, and
\begin{equation}\label{d19:eq:parity}
(x,y)\equiv\tau(x)\tau(y)\pmod 2,\qquad x^2\equiv\tau(x)\pmod4,
\end{equation}
where $\tau(x)=(x,v_1)\bmod2\in\{0,1\}$. Indeed all original generators have square $1\pmod4$ and pair oddly, and each adjoined generator has square 4 and pairs evenly with the preceding lattice. The claims follow by expanding integral combinations. Also,
\begin{equation}\label{d19:eq:sat}
s^2=4,\quad (s,M)\subseteq2\Z\quad\Longrightarrow\quad s\in M.
\end{equation}
There is no unit vector in $M$: such a vector $a$ would have $(a,v_i)$ odd for all $i$ and hence $a^{\mathsf T}Fa\ge73$, contradicting $\lambda_{\max}(F)\le19+\sqrt{70}$. Consequently $(s,v_i)\in\{0,\pm2\}$ for $s\in R$; otherwise $v_i\mp s$ would have square one. By \eqref{d19:eq:parity}, every $s\in R$ pairs evenly with $M$, and reflection $x\mapsto x-(x,s)s/2$ preserves $M$. Thus $R/\sqrt2$ is a nonempty simply laced root system with root lattice $Q/\sqrt2$ \cite[Ch.~VI]{Bourbaki}.
\subsubsection{Local active-vector graphs}
Put $H=F-19I$, $r=\dim W$, and $c=19-r$. For a unit root direction $e=s/2$, orient the active vectors as $v_i=e+z_i$, where $z_i\perp e$ and $z_i^2=4$. Define $G_e$ on these $t_e$ vectors by adjacency if $(z_i,z_j)=-2$, and let $b_e$ be its edge count. Their Gram matrix is $4I-2A(G_e)\succeq0$, so $\rho(A(G_e))\le2$ and $b_e\le t_e$. Each edge gives a norm-four vector $z_i+z_j\in R\cap e^\perp$.
\begin{lemma}\label{d19:lem:rankedges}
For each $e$, $b_e\le 2(r-1)$. If equality holds, $r-1$ is even.
\end{lemma}
\begin{proof}
Let $D_e$ be the span of the edge sums. A nontrivial connected component on $k$ vertices has either $k-1$ edges (a tree) or $k$ edges (a cycle, by equality in the average-degree bound). In the tree case the unsigned incidence matrix has rank $k-1$, and the residual Gram matrix has nullity at most one by Perron--Frobenius \cite[Proposition~1.3.9 and Theorem~2.2.1]{BH}. The edge sums therefore span at least $k-2$ dimensions; for $k=2,3$ direct use of the positive definite Gram matrix gives ranks 1,2. In the cycle case the incidence rank is $k$ or $k-1$ according as $k$ is odd or even, giving a lower bound of $k-1$ or $k-2$ on the edge-sum rank. These bounds give at least half the edge count, with equality possible only for a four-cycle, whose edge-sum rank is two. Different components are orthogonal. Thus $b_e\le2\dim D_e\le2(r-1)$, and equality forces $r-1$ even.
\end{proof}
The identities
\begin{align}
\|He\|^2&=(t_e-19)^2+4t_e-4b_e\label{d19:eq:rootenergy}\\
&=(t_e-17)^2+72-4b_e\notag
\end{align}
follow from $Fe=t_e e+\sum z_i$. For an irreducible component $W_j$ of rank $d_j$ and Coxeter number $h_j$, with one $e$ per antipodal pair of roots, Weyl invariance gives $\sum_e ee^{\mathsf T}=(h_j/2)P_j$, where $P_j$ is the orthogonal projection onto $W_j$. (The trace fixes the scalar, since there are $d_jh_j/2$ root directions.) Lemma~\ref{d19:lem:rankedges} then gives $70\ge r(80-8r)$, excluding $1\le r\le9$. If $r=10$, equality $b_e=18$ is impossible; thus $b_e\le17$. Writing $T=\tr(P_W F)$ and using Jensen's inequality in the ten-dimensional root space and the nine-dimensional complement yields
\[
70\ge40+\frac{(T-170)^2}{10}+\frac{(194-T)^2}{9}\ge40+\frac{576}{19}>70.
\]
Hence $11\le r\le19$.
\begin{lemma}[Root-system deficit]\label{d19:lem:deficit}
Let $T_j=\tr(P_jF)$ and $D_j=(2/h_j)\sum_{e\in W_j}(t_e-b_e)$. Then $D_j\ge0$ and
\begin{equation}\label{d19:eq:deficit}
\sum_j\frac{(T_j-19d_j)^2}{d_j}+4\sum_jD_j+
\frac{(365-\sum_jT_j-19c)^2}{c}\le70.
\end{equation}
For $c=0$ the last term is omitted and $\sum_jT_j=365$.
\end{lemma}
\begin{proof}Use \eqref{d19:eq:rootenergy} in its first form, the tight-frame identity and Jensen's inequality in each component; on the complement use $\|P_{W^\perp}HP_{W^\perp}\|_F^2\ge(\tr(P_{W^\perp}H))^2/c$.\end{proof}
\subsection{Weight labels and necessary inequalities}\label{d19:sec:weights}
For $v=p+z$ ($p\in W$, $z\in W^\perp$), the projection $p/\sqrt2$ pairs with all roots of $Q/\sqrt2$ in $\{0,\pm1\}$. It is therefore zero or a minuscule weight in each irreducible component. This follows by Weyl conjugation to a dominant weight: the highest-root pairing is at most one, so only fundamental weights at coefficient-one nodes can occur. Each minuscule orbit is the set of shortest representatives of its discriminant class \cite[Corollary~2.2.5]{Lurie}. Write $a$ for the class in $(Q/\sqrt2)^*/(Q/\sqrt2)$, and $w(a)$ for twice its minimum norm. A \emph{label} is $(a,z)$ up to simultaneous negation. Its multiplicity is the number of original lines with that label. Its support is the set of nonzero root components of $a$. The empty support labels have multiplicity one. We use $h=5-w(a)=z^2$ for the residual square (not a Coxeter number).
\begin{lemma}[Orbit extrema]\label{d19:lem:extrema}
For distinct labels $(a,z)$ and $(b,z')$ up to sign, every pairing between their Weyl images belongs to $\{-1,1\}$. The maximum and minimum of the projected pairing are
\[
\frac{w(a)+w(b)-w(a-b)}2,\qquad
\frac{w(a+b)-w(a)-w(b)}2.
\]
The interval of possible pairings has width at least two in any component where both classes are nonzero. For $A_d$ it has width $2\min(k,l,d+1-k,d+1-l)$, with the usual subset indices $k,l$.
\end{lemma}
\begin{proof}All Weyl images lie in $M$ and have the original class modulo $Q$ and residual. Pairings are odd. A pairing $\pm5$ makes the vectors equal up to sign, and a pairing $\pm3$ makes a difference or sum a norm-four vector in $Q$; either identifies labels. For the maximum choose dominant representatives; their difference pairs with each positive root in $\{0,\pm1\}$, hence is shortest in its class. Apply the same argument to the negative class for the minimum. A direct computation in the minuscule weight orbits gives the asserted widths; in $A_d$ use the centered subset vectors $\sqrt2(\mathbf1_S-|S|\mathbf1/(d+1))$. For $D_d$ use vector and half-sign spin weights; in $E_6,E_7$ the maximum-minus-minimum follows from their discriminant norm data in Table~\ref{d19:tab:slots}.\end{proof}
Consequently, distinct labels share at most one component, and a support of size at least two has at most one label. Any class with orbit pairing width at least four against itself can appear in at most one label in its component (up to sign). The required component data are as follows. We call each row a slot, and may group the two signs of a class since only bounds are used.
\begin{table}[ht]
\centering\small
\begin{tabular}{lllll}\toprule
Type & Class & $w$ & $\delta=w-w(2a)/4$ & Orbit size\\\midrule
$A_d$ & $k\le(d+1)/2$ & $2k(d+1-k)/(d+1)$ & $k$ & $\binom{d+1}{k}$\\
$D_d$ & vector & $2$ & $2$ & $2d$\\
$D_d$ & spin (either) & $d/2$ & $\lfloor d/2\rfloor$ & $2^{d-1}$\\
$E_6$ & nonzero & $8/3$ & $2$ & $27$\\
$E_7$ & nonzero & $3$ & $3$ & $56$\\\bottomrule
\end{tabular}
\caption{Nonzero minuscule classes in square-four normalization. $D_4$ has three separate, triality-equivalent slots.}\label{d19:tab:slots}
\end{table}
If $h>0$, necessarily $\delta(a)\le3$: the least pairing of $v$ with its Weyl images is $5-2\delta(a)$. Along simple reflections the pairing with $v$ changes by $0$ or $\pm2$. If $\delta(a)\ge4$, some pairing equals $-3$, and the sum is a norm-four vector with nonzero residual $2z$, impossible.
\subsubsection{Multiplicity bounds}
A label of multiplicity $m$, active dimension $d$, and residual square $h$ satisfies
\begin{equation}\label{d19:eq:cap}
m\le B(d+\mathbf1_{h>0}),\qquad (h-1)m\le4.
\end{equation}
Indeed its vectors span at most the specified dimension, and orienting their residuals identically gives $m^2h\le\|\sum v\|^2\le m^2+4m$. For a fixed choice on the other components, an $A_d$ slot with $k$-subsets can contain at most $\lfloor\binom{d+1}{k-1}/k\rfloor$ representatives: two different subsets sharing $k-1$ entries differ by a norm-four root, forcing pairing 3. Similarly, a vector slot of $D_d$ has fiber bound 2 (different indices give a difference of square four), and a spin slot has fiber bound at most $\lfloor2^d/(d+1)\rfloor$, by disjoint radius-one binary Hamming balls \cite{Hamming}. The $D_4$ spin slots have fiber bound 2 by triality. For the two exceptional types the full orbit size is a safe fiber bound. These bounds multiply by the orbit sizes of all other slots. For five $A_1$ components and zero residual, the maximum is 5: orient signs to even parity; the translated size-two supports intersect pairwise and thus form a star or a triangle.
The singleton-only labels on one $A_d$ factor ($d\ge2$) have total multiplicity at most $\lfloor4+8/(d-1)\rfloor$. Orient by the singleton class. Each such label has multiplicity one, and the residual Gram off-diagonals are $-1+2/(d+1)$ by Lemma~\ref{d19:lem:extrema}. Each residual vector has squared norm $3+2/(d+1)$. The norm of the sum gives the claim. For an $A_1$ singleton support the same argument gives mutually orthogonal residuals, so the bound is $c$. The nonsingleton $A_d$ slots, each $D_d$ slot (with spin signs grouped), and the exceptional slots can each appear in at most one label per component. This is the width-four test of Lemma~\ref{d19:lem:extrema}; for odd $D_d$ the spin signs are a single class up to sign, and for even $D_d\ge6$ either choice of chiralities has width at least four.
\subsubsection{Deficit of a single label}
For two lines let $d_{ij}=v_i-(v_i,v_j)v_j$, of squared norm 8. A pair contributes an edge in $G_e$ exactly if $|(d_{ij},e)|=2$. If this occurs, $d_{ij}$ is a sum of two orthogonal square-four roots up to sign, hence belongs to $Q$ and the two lines have the same label. Each $G_e$ therefore splits over labels. Its deficit on each label is nonnegative, by the induced Gram matrix argument. For a label with active component energies $w_i$, write $T_i=mw_i$, and let $C_i$ be an upper bound on the weighted edge count per pair in component $i$. The total deficit of this label is bounded below by
\begin{equation}\label{d19:eq:labeldef}
\max\left(0,\ m\sum_iw_i-C\binom m2,\ \sum_i\max\{0,mw_i-C_i\binom m2\}\right),
\end{equation}
where $C$ bounds the sum of weighted contributions of a pair. Indeed the weighted active counts are $mw_i$ by the tight-frame identity. The norm-eight vector in $Q$ can have either two nonzero components of square four, or one of square eight. The first case has total weight $2/h_i+2/h_j$ with Coxeter numbers $h_i,h_j$. In the second case its component contribution is at most 2 by the tight-frame identity, and in $A_d$ it is exactly $8/(d+1)$ when it contributes: the unscaled vector has two entries $+1$ and two entries $-1$, giving four root directions. A component with $\delta=1$ cannot support such a within-component contribution: any difference of its shortest representatives has square at most 4. For the zero-residual sign ambiguity, replace the second representative by its negative as necessary; the same shortest-orbit argument applies if the class is self-inverse, and otherwise a nontrivial sum cannot belong to that component's root lattice. Thus the maxima over the listed alternatives give $C,C_i$.
For a label using only $A$ singleton slots, orient consistently by a nonzero residual or a non-self-inverse class. The representatives are words on alphabets of sizes $q_i=d_i+1$. Distances between different words are two or three, since their full squared difference is four times the Hamming distance. A contributing pair has distance two. For a fixed pair of coordinates $i,j$, the other coordinates partition the words into fibers of size at most $L_{ij}=\min(q_i,q_j)$. The number of pairs at distance two on these coordinates is at most
\[
P(m,L)=\lfloor m/L\rfloor\binom L2+\binom{m\bmod L}2.
\]
The edge weight for such a pair is $2/q_i+2/q_j$. Hence another lower bound on the deficit is
\begin{equation}\label{d19:eq:hammingdef}
m\sum_i w_i-\sum_{i<j}(2/q_i+2/q_j)P(m,L_{ij}).
\end{equation}
The componentwise bounds obtained by keeping only the terms incident with a given coordinate can also be summed after truncation at zero. The only case not using this refinement is the zero-residual all-$A_1$ label. All these inequalities are encoded in \texttt{deficit(p,m)} in the accompanying verifier.
\subsection{Exact certificates for Theorem~\ref{d19:thm:A}}\label{d19:sec:certs}
\subsubsection{The finite relaxation}
Exclude a component of type $E_8$ immediately, since its only minuscule projection is zero, contradicting the spanning condition. For $D_d$, $d\ge11$, the spin weight has square greater than 5; the vector weights would give $T_j\le146$, whereas $(19d-146)^2/d>70$. Thus the component types in the search are $A_1,\ldots,A_{19}$, $D_4,\ldots,D_{10}$, $E_6,E_7$.
Let $n_t$ count the components of type $t$. A pattern $p$ is a multiset of nonzero slots of Table~\ref{d19:tab:slots} on distinct components (retaining their types, but forgetting their indices and the signs of the classes). Enumerate the patterns of total energy at most 5, discarding those with positive residual and $\delta>3$. For each allowed multiplicity $m$ from the bounds above, let $x_{p,m}\ge0$ count the labels of that pattern and multiplicity. Include a single variable for the number of empty-support labels. The following constraints are necessary; they are sufficient for the certificate calculation, not sufficient for realization as lines.
\begin{enumerate}
\item $\sum_{p,m} m x_{p,m}=73$. The count of empty-support labels is at most $B(c)$, and the count of zero-residual lines is at most $B(r)$.
\item Write $s_t(p)$ for the number of slots of type $t$ in $p$. Two distinct labels cannot use the same pair of factors. Therefore
\[
\sum_{p,m}x_{p,m}\binom{s_t(p)}2\le\binom{n_t}2,
\quad \sum_{p,m}x_{p,m}s_t(p)s_u(p)\le n_t n_u\quad(t\ne u).
\]
The one-label-per-slot and singleton multiplicity bounds from Section~\ref{d19:sec:weights} give the corresponding inequalities summed over all factors of each type.
\item At each component of type $t$, all the other components occurring in nonsingleton supports containing it are distinct. If $k\ge2$, put $a=n_u-\mathbf1_{t=u}$. Then
\[
\sum_{p,m}x_{p,m}s_t(p)\mathbf1_{s_u(p)-\mathbf1_{t=u}\ge k}
\le n_t\lfloor a/k\rfloor.
\]
If $n_t=5$, at most two distinct labels use at least three factors of this type, since three triples with pairwise intersections at most one have union of size at least six.
\item For each choice of $k_t\le n_t$ factors of every type, put $d=\sum_t k_t\rk(t)$. Averaging the bound on the lines avoiding those factors gives
\begin{equation}\label{d19:eq:averaged}
\sum_{p,m}m x_{p,m}\prod_t\frac{\binom{n_t-s_t(p)}{k_t}}{\binom{n_t}{k_t}}\le B(19-d).
\end{equation}
Restricting the sum to zero-residual labels gives also the bound $B(r-d)$.
\item A zero-residual label using an $A_d$ component of subset size $k$ can share that component with another label only when the two subset sizes (after changing signs) are $1$ and $(d+1)/2$. To see this, the maximum and minimum in that component must be $1,-1$; for subsets of sizes $k,l\le(d+1)/2$ they are $2\min(k,l)-2kl/(d+1)$ and $-2kl/(d+1)$. Their width and midpoint give the assertion. If the slot is not half-size, the energy on that factor is consequently at most $mw$ plus the maximum contribution of a single half-size slot in a compatible pattern, if such a slot can occur. Here compatibility requires disjoint remaining components, and the half-size slot has a single label by the width test. Discard any variable for which this upper bound is less than
\begin{equation}\label{d19:eq:blocklower}
\frac{365d}{19}-\sqrt{(70-16/19)d(19-d)/19}
\end{equation}
for the factor's rank $d$. This follows by centering $F$ at $(365/19)I$ and applying Cauchy--Schwarz to its trace against a rank-$d$ projection.
\item When $c=1$, all positive residual squares lie in a single rational square class: the pairings of residuals are rational by subtracting the rational projected pairings, and their squares multiply to the square of a rational pairing. We may check each possible square class separately; there must be a positive residual since the original vectors span.
\end{enumerate}
All the preceding constraints use only integer or rational arithmetic (the comparison with \eqref{d19:eq:blocklower} is made by squaring positive rational quantities). For the energy bounds, merge the equal-type factors into blocks, using Jensen's inequality. Let $d_i$ be the block ranks, including the complementary block when $c>0$. Write $T_i=\sum m w_i(p)x_{p,m}$ for the energy in that block, and $D$ for the sum of the lower bounds \eqref{d19:eq:labeldef}--\eqref{d19:eq:hammingdef}, with the coefficient $x_{p,m}$ on each label. By Lemma~\ref{d19:lem:deficit},
\[
\sum_i (T_i-19d_i)^2/d_i+4D\le70.
\]
In particular $4D\le70$, and for any rational vector $(g_i)$,
\begin{equation}\label{d19:eq:tangent}
2\sum_i g_i T_i+4D\le70+\sum_i d_i(38g_i+g_i^2).
\end{equation}
The verifier includes the two-sided trace bounds from \eqref{d19:eq:blocklower} on every union of the blocks (rounding the radius upwards to the next hundredth), and appends any additional inequalities \eqref{d19:eq:tangent} listed in the certificate. If $c=0$, the residual energy must be zero; otherwise it is included as a block of dimension $c$.
\subsubsection{Additional cuts for the exceptional patterns}
Three small strengthening arguments are used in the certificate generator, and are recorded here to make the relaxation auditable.
\begin{enumerate}
\item Suppose there are seven $A_2$ components and $c<5$. Let $p$ count labels supported on three of these components and nowhere else. Their residuals have square one, and pair with absolute value $1/3$ if the supports meet and 1 otherwise. If two supports are disjoint, their residuals are collinear and every other support must meet both (a support cannot miss one and meet the other). Any such third triple uses one point from each of the disjoint triples and the remaining seventh point; there are at most three of them by the pair-packing bound. Thus $p\le5$ in this case. Otherwise the residual directions form cosine-$1/3$ equiangular lines in dimension $c$, giving $p\le\lfloor 8c/(9-c)\rfloor$. Take the maximum of the two bounds. In the special type $A_1A_2^7$ of rank 15, let $q$ count $A_1A_2^2$ labels. Then $q\le3$ because the pairs in the seven-set form a matching; $p\le6$ by the preceding argument. If $p=6$, the six pure triples consume 18 of the 21 pairs. The unused three pairs have all degrees even and form a triangle, forcing $q\le1$. Hence $2p+q\le13$.
\item For $A_1^3D_6D_{10}$, the rank is 19 and there is no residual. A label on $A_1$, a $D_6$ vector slot and a $D_{10}$ vector slot has representatives $e+\sqrt2\epsilon f_i+\sqrt2\eta g_j$ with the $A_1$ coordinate positive. Fixing $i,\epsilon$ forces $j$ to be the same for all vectors in the fiber (otherwise the pairing is 3), so there are at most two choices of $\eta$. Thus its multiplicity is at most $4\cdot6=24$.
\item For $A_5^3D_4$, also of full rank, a label with the middle subset weight of $A_5$ and any nonzero $D_4$ slot has multiplicity at most 12. By triality orient its representatives as $\sqrt2(\mathbf1_T-\tfrac12\mathbf1)+\sqrt2 e_i$, $|T|=3$, $1\le i\le4$. In a fixed group $i$, the triples intersect in at most one point. A complementary pair precludes a third triple; otherwise the $A_5$ projections have mutual pairing $-1$ and norm squared 3, so at most four occur. Distinct groups cannot use the same antipodal triple direction, and there are ten such directions. If $b$ groups have complementary pairs, the count is bounded by both $10+b$ and $2b+4(4-b)$, hence by 12. In particular the type itself can be excluded directly: the only other possible pattern of energy five is three $A_5$ singleton slots, of multiplicity at most $6^2=36$. The bound is $36+3\cdot12=72$.
\end{enumerate}
\subsubsection{Verification and coverage}
\paragraph{Geometric source of the certificate rows.}
Apply the saturation construction \eqref{d19:eq:parity}--\eqref{d19:eq:sat} to a hypothetical 73-line system with a dual unit vector. Each actual label contributes one to the variable for its pattern and multiplicity, and each empty-support label contributes one to the empty-support variable. The minuscule weights, the residual-square restriction, and the multiplicity and fiber bounds of Section~\ref{d19:sec:weights} give the pattern and multiplicity filters in \texttt{cap} and \texttt{patterns}. These variables account for all 73 lines.

In \texttt{setup}, Lemma~\ref{d19:lem:extrema} and the singleton residual bounds give the slot, pair and higher support-packing bounds. The subspace bound \eqref{d19:eq:bound} gives the empty-support, zero-residual and averaged projection rows \eqref{d19:eq:averaged}. The frame identity \eqref{d19:eq:frame} gives the block-trace bounds through \eqref{d19:eq:blocklower}, and combining these with the orbit extrema gives the zero-residual filter \texttt{allowed}. The one-dimensional residual argument gives the square-class row when a branch is used. The additional rows for the specified types use the ternary support packing and the two full-rank fiber bounds proved above.

The same-label edge counts give the coefficients in \texttt{deficit} through \eqref{d19:eq:labeldef} and \eqref{d19:eq:hammingdef}. Lemma~\ref{d19:lem:deficit} then gives the total deficit bound and each stored tangent row \eqref{d19:eq:tangent}. Thus the label counts of a counterexample would satisfy the rows reconstructed for its root type and residual branch, as applicable.

Let $Ax\le b$ be the resulting system of rational inequalities with nonnegative variables. A stored certificate consists of rational tangent parameters $g$ and a nonnegative rational vector $y$ of row multipliers \cite[Corollary~7.1f]{Schrijver}. Since $\sum x_{p,m}\le73$, set $\ell=\min(0,\min_j(y^{\mathsf T}A)_j)$. Then
\[
y^{\mathsf T}Ax\ge\ell\sum x_{p,m}\ge73\ell.
\]
Consequently $y^{\mathsf T}b<73\ell$ certifies infeasibility. This is the sole arithmetic acceptance test of \texttt{verify\_dual}. Floating-point optimization may be used to discover a certificate, but is not needed to verify it.
The supplied files \nolinkurl{root19_all_mixed_certificates.jsonl} and \nolinkurl{root19_homogeneous_DE_recheck.jsonl} contain 6571 and 140 certificates, respectively. The former includes seven one-dimensional residual branches for $A_1^2A_3^4D_4$; the expected square classes are generated from the permitted patterns (including the empty support). For coverage, enumerate all unordered multisets of the specified component types with rank between 11 and 19, and compare the resulting keys with the disjoint union of the two certificate type lists and the four types in the next subsection. Equivalently, before excluding $E_8$ and large $D_d$ components, the rank counts are the coefficients of
\[
\prod_{j\ge1}(1-t^j)^{-1}\prod_{j\ge4}(1-t^j)^{-1}
\prod_{j=6}^8(1-t^j)^{-1}.
\]
Table~\ref{d19:table:root-coverage} gives the counts.
\begin{table}[htbp]
\centering\small
\caption{Root-type counts for the dimension-nineteen rational certificates.}\label{d19:table:root-coverage}
\begin{tabular}{rrrr}\toprule
Rank & All types & Allowed components & Mixed certificate types\\\midrule
11&128&124&118\\12&193&185&172\\13&276&263&254\\14&403&380&366\\15&570&533&519\\16&815&754&735\\17&1137&1042&1025\\18&1599&1449&1422\\19&2207&1979&1954\\\midrule
Total&7328&6709&6565\\\bottomrule
\end{tabular}
\end{table}
The set comparison is checked by \texttt{check\_coverage.py} in the accompanying files: $6709=6565+140+4$. This verification concerns all root types, not just those for which the certificate search happens to have succeeded.
\subsubsection{The four remaining types}\label{d19:sec:special}
The four types not covered by the rational certificates are \(A_1^{11}\), \(A_1^{12}\), \(A_1^{13}\) and \(A_2^7\).
\paragraph{Orthogonal roots.}
Suppose the root system is $A_1^r$, $r=11,12,13$. Write $v_i=p_i+z_i$ with $p_i\in W$ having coordinates $0,\pm1$ in the orthonormal root directions. Put $I_i=\supp p_i$. Reversing all its signs shows $|I_i|\ne4$, since if $|I_i|=4$, $v_i+gv_i=2z_i$ would be a norm-four root outside $W$. By Lemma~\ref{d19:lem:extrema}, distinct supports meet in at most one coordinate; their residuals are orthogonal if they meet and have absolute pairing one if they are disjoint. A quintet (support size five) has residual zero, so must meet every other support exactly once. For a support of size at least two all its lines have the same label. The size-two, size-three and size-five multiplicities are at most 2, 4 and 5: orient to the same residual when nonzero, and use binary sign words with distances two or three. For size five use the even-parity orientation and the intersecting two-set argument above.
Let $t_j$ be the number of active lines in coordinate $j$, $Q_1=\sum_jt_j$, $B_1=\sum_j b_{e_j}$, and $Z=\sum_i z_ip_i^{\mathsf T}$. Splitting the frame energy into blocks gives
\begin{equation}\label{d19:eq:A1energy}
70\ge72r+\sum_j(t_j-17)^2+\|Z\|_F^2-4B_1+
\frac{(4+19r-Q_1)^2}{19-r}.
\end{equation}
For a pair of lines contributing an edge at some coordinate, their signed difference of square 8 is a sum of two orthogonal roots. It therefore contributes at exactly two coordinates on the same support. Contributions to $Z$ from distinct supports are Frobenius orthogonal, by either disjointness of projected supports or orthogonality of residuals. Counting sign pairs gives the contribution of each occupied support to $\|Z\|_F^2-4B_1$ in Table~\ref{d19:table:orthogonal-support}, where $m$ is the multiplicity.
\begin{table}[htbp]
\centering\small
\caption{Contributions to the energy for a single orthogonal support.}\label{d19:table:orthogonal-support}
\begin{tabular}{c|rrrrr}
$|I|\backslash m$&1&2&3&4&5\\\hline
1&4&8&12&16&20\\2&6&$-8$&--&--&--\\3&6&0&$-18$&$-48$&--\\5&0&$-8$&$-24$&$-48$&$-80$
\end{tabular}
\end{table}
For the singleton support the residuals are mutually orthogonal, so the value is $4m$ for any allowed $m$. For sizes two and three orient the residuals equally: the $Z$-term is $(5-|I|)\|\sum p_i\|^2$, and the $B_1$-term is twice the number of sign pairs at distance two. For size three a complementary pair precludes a third word; in its absence the words have common parity. For quintets the residual is zero and every pair contributes at two axes. These observations give the table including its distance-three exception.
\paragraph{No quintet.} Let $o,d,t$ count lines on singletons, pairs and triples. From \eqref{d19:eq:A1energy} and Jensen's inequality,
\begin{equation}\label{d19:eq:A1no5}
70\ge72r+4o-4d-12t+(4+2r)^2/19.
\end{equation}
Since $o+d+t\le73$, this is impossible for $r=13$; for $r=12$ it gives $t\ge68$. If $t=68$ it forces $d=5,o=0$; then the contribution lower bound in \eqref{d19:eq:A1energy} is $28+(214-204)^2/12+(232-214)^2/7>70$. Hence $t\ge69$ and at least 18 distinct triples occur. Their residuals have squared norm two and pair to 0 or $\pm1$ according as their supports meet or are disjoint. In any selection of 19 triples, the 57 incidences on 12 points force at least $9\binom52+3\binom42=108$ intersecting pairs. For the residual frame $T$ on these 19 vectors, $\tr T=38$, $\tr T^2\le202$, impossible in dimension seven. Thus exactly 18 triples occur. The 54 incidences force at least $6\binom52+6\binom42=96$ intersecting pairs, so the frame on these residuals satisfies $\tr T=36$, $\tr T^2\le186$. The rank inequality and even integrality of $\tr T^2$ give equality $\tr T^2=186$ and rank seven. Then $P=T-5I$ satisfies $\tr P=\tr P^2=1$. It preserves the full even integral lattice generated by the residuals, because their inner products are integers. The product of its nonzero eigenvalues is a nonzero integer, whereas for $k\ge2$ such eigenvalues the product has absolute value at most $(1/k)^{k/2}<1$. Thus $P$ is a rank-one orthogonal projection. Every residual of square two splits into two orthogonal vectors of that even lattice, so lies either in the image or kernel of $P$. But $\tr(PT)=6$ forces three distinct residuals on the same line, impossible.
For $r=11$, \eqref{d19:eq:A1no5} forces $t\ge59$, hence at least 15 distinct triples. The residual Gram matrix modulo 2 is realized by the words $(1,\mathbf1_I)$ in $\F_2^{12}$ for all supports. Indeed $(v_i,v_j)$ is odd and projected pairings reduce to the intersection parity. The rank of the resulting bilinear form on their span $U$ is at most eight. Nondegeneracy of the ambient dot product implies $2\dim U-12\le8$, and $U^\perp$ has a nonzero vector of the form $(0,\mathbf1_A)$. Thus all supports meet $A$ in an even number of points. If $a=|A|\le4$, total activity on $A$ is at most $8\binom a2$ (each pair in $A$ can support at most four lines); its energy deviation $(19a-T_A)^2/a$ alone exceeds 70, where $T_A$ is this activity. Hence $a\ge5$. Triples contained in the complement of size $b=11-a$ number at most $P(b)=\lfloor b\lfloor(b-1)/2\rfloor/3\rfloor$ (take $P(0)=0$, $P(5)=2$). All other triples consist of a pair from $A$ and one point from its complement. They give matchings indexed by the latter point, so there are at most $\min\{\binom a2,b\lfloor a/2\rfloor\}$. For $a=5,\ldots,11$ the combined bounds are $14,17,13,13,8,5,0$. Only $a=6$ can support at least 15 triples. At most two triples lie in its five-point complement, leaving at least 13 edges of $K_6$ colored by five matchings. The matching-size vectors, after sorting, are exactly
\[
(3,3,3,3,3),\ (3,3,3,3,2),\ (3,3,3,3,1),\ (3,3,3,2,2).
\]
The verifier \texttt{root\_rank19\_pure\_A1\_certificate.py} enumerates such disjoint matchings, normalizing one size-three matching by a vertex permutation. The residual Gram has diagonal two and nonzero off-diagonals $\pm1$ exactly between disjoint edges of distinct colors. Signs on a spanning forest may be made positive by switching residual signs. A triangle of nonzero pairings must have positive product whenever another residual has a nonzero pairing with exactly one vertex of that triangle: a negative product would give a null relation on the three triangle vectors with all coefficients nonzero, contradicting the extra pairing. Enumerate all remaining signs subject to these necessary conditions. The counts of matching records are $2,24,24,288$, and of signed matrices $128,768,384,13248$, respectively. Every matrix has rank at least nine over $\F_3$, contradicting real rank at most eight. This enumeration is finite and uses only integer elimination.
\paragraph{At least one quintet.} There are at most three quintets, since four have union size at least $20-\binom42=14$. For clarity we give the constraints for the integer checks used in the same verifier. They are relaxations obtained from \eqref{d19:eq:A1energy}, replacing each coordinate block's sum of squares by the square of its total deviation divided by its size. In each case the total number of lines is 73. A quintet of multiplicity $a$ contributes $-4a(a-1)$ to the energy. The contributions from other occupied supports are bounded below by $4,-4,-12$ per line for sizes 1,2,3 respectively. For $q$ quintets carrying $k$ lines in total, let $C(k,q)$ be the minimum of $-4\sum a_i(a_i-1)$ over $1\le a_i\le5$, $\sum a_i=k$.

\emph{One quintet:} with multiplicity $a$, put $s=r-5$, $c=19-r$. All other supports meet it exactly once. Let $d,t$ count pair and triple lines. There are at most $\min\{\binom s2,5\lfloor s/2\rfloor\}$ triples, and the incidences on cross pairs give $\lceil d/2\rceil+2\lceil t/4\rceil\le5s$. With $U=d+2t$, the energy is at least
\[
72r+292-4a^2-8U+\frac{(4a-12)^2}5+
\frac{(U-17s)^2}s+\frac{(19r-69-4a-U)^2}c.
\]
\emph{Two quintets:} they meet at one point, with eight private points and $s=r-9$ exterior points. Other supports either contain the common point, with line counts $p_1,p_2,p_3$ by size, or meet both private sets, with counts $a_2,a_3$. The activity totals on the common, private and exterior blocks are
\[
T_0=k+p_1+p_2+p_3,\quad T_1=4k+2(a_2+a_3),\quad T_2=p_2+2p_3+a_3.
\]
The necessary packing inequalities are $\lceil p_2/2\rceil+2\lceil p_3/4\rceil\le s$, $p_2\le2s$, $p_3\le4\lfloor s/2\rfloor$, and $\lceil a_2/2\rceil+\lceil a_3/4\rceil\le16$, $a_2\le32$, $a_3\le\min(64,16s)$. The latter bounds use the 16 pairs between the two private sets and the four-by-$s$ pairs between either private set and the exterior. In fact the sharper bound $a_3\le16s$ follows by each triple consuming one such pair. The energy lower bound is
\begin{align*}
72r+C(k,2)+4p_1-4(p_2+a_2)-12(p_3+a_3)
&+(T_0-17)^2+\frac{(T_1-136)^2}{8}\\
&+\frac{(T_2-17s)^2}{s}+\frac{(4+19r-\sum T_i)^2}{c}.
\end{align*}
\emph{Three quintets with distinct intersections:} their union has 12 points, of which three are pairwise intersections and nine are private. Put $s=r-12$. Other supports either contain one intersection and one private point of the third quintet (possibly with an exterior point), or consist of three private points, one from each quintet. Write $d,t,u$ for the corresponding line counts. Then $d\le18$, $t\le4\min(9,3s)$, $u\le36$, and $\lceil d/2\rceil+\lceil t/4\rceil\le9$, by counting the nine available intersection--private pairs, the $3s$ intersection--exterior pairs, and the nine pairs between two of the private triples. The activity totals on the three blocks are $(2k+d+t,3k+d+t+3u,t)$. The energy bound is $72r+C(k,3)-4d-12(t+u)$ plus the three block-deviation terms and the complementary trace term, with block sizes $3,9,s$. Omit the exterior term if $s=0$ (then $t=0$).
\emph{Three concurrent quintets:} their union has 13 points. Thus $r=13$ and any other support is either the common singleton or a triple of private points, one per quintet. If $p,t$ are the line counts, $t\le64$ by the 16 pairs between two private sets. The energy bound is $72\cdot13+C(k,3)+4p-12t$ plus the block trace terms for activities $k+p,4k+3t$ on blocks of size 1 and 12, and for the complement.

The accompanying verifier computes the minima of these explicit lower bounds over the necessary integer domains, as recorded in Table~\ref{d19:table:quintet-bound}. Here ``empty'' means that the integer constraints already have no solution.
\begin{table}[htbp]
\centering\small
\caption{Minimum energy bounds on the integer domains with quintet supports.}\label{d19:table:quintet-bound}
\begin{tabular}{rrrrr}\toprule
$r$ & $q=1$ & $q=2$ & $q=3$ distinct & $q=3$ concurrent\\\midrule
11&$454/5$&$280$&--&--\\
12&$2578/35$&$2014/21$&empty&--\\
13&$2002/15$&$287/4$&$1648/9$&$239/3$\\\bottomrule
\end{tabular}
\end{table}
All finite minima exceed 70, completing the exclusion of these types.
\paragraph{The type \texorpdfstring{$A_2^7$}{A2 (seven copies)}.}
Here the residual space is five-dimensional. Each nonzero component has projection square $4/3$ and labels have supports of size at most three. Size-two and size-three labels have multiplicity at most three and nine by the fiber bound. Singleton and empty supports have multiplicity one per label. Write $m_\ell,s_\ell,h_\ell=5-4s_\ell/3$ for the multiplicity, support size and residual square, and $n_j$ for the number of lines active in component $j$. Put $J=\sum_j n_j$. The energy inequality implies
\[
(4J/3-266)^2/14+(270-4J/3)^2/5\le70,
\]
so $J\ge190$. Distinct supports of size at least two meet in at most one component; any two distinct labels have residual absolute pairing $1/3$ if their supports meet and 1 otherwise. Consequently the residual frame $S=\sum_i z_iz_i^{\mathsf T}$ satisfies the exact identities
\begin{align}
\tr S&=365-4J/3,\label{d19:eq:res-trace}\\
\tr S^2&=73^2-\frac89\sum_jn_j^2+\sum_\ell m_\ell^2(h_\ell^2-1+8s_\ell/9).\label{d19:eq:res-gram}
\end{align}
For \eqref{d19:eq:res-gram}, start with one for every ordered pair of original lines; subtract $8/9$ for every active intersection and correct all same-label pairs to $h_\ell^2$. The number $p$ of triple supports is at most seven. Pair packing and the multiplicity bounds give $J\le73+3(21-3p)+18p=136+9p$, hence $p\ge6$. For $p=6$, equality $J=190$ forces 54 triple lines, nine pair lines, and ten singleton lines; the three unused pairs must be used. Thus \eqref{d19:eq:res-gram} yields
\[
\tr S=335/3,\quad
\tr S^2\le73^2-8\cdot190^2/63+1296+168+400/3=147559/63.
\]
This gives $(\tr S)^2-5\tr S^2\ge47780/63>0$, contrary to the rank inequality.

\begin{samepage}
If $p=7$, the triples partition the pairs and no two-support label occurs. Write $63-\delta$ for the number of triple lines, $o$ for the singleton count, and $e=10+\delta-o$ for the empty count. Since $J=3(63-\delta)+o\ge190$ and $o\le10+\delta$, we have $0\le\delta\le4$, $o\le14$. The sum of the squared triple multiplicities is at most $567-18\delta+\delta^2$ (write them as $9-\delta_i$ and use $\sum\delta_i^2\le\delta^2$). Substituting this and $\sum n_j^2\ge J^2/7$ into \eqref{d19:eq:res-gram} gives
\[
9(5\tr S^2-(\tr S)^2)\le
-528\delta^2/7+912\delta o/7-2736\delta-152o^2/7+72o-396<0.
\]
For the last inequality, bound the mixed term using $o\le14$ and note that the quadratic in $o$ has negative maximum. Together with the certificate verification above, this proves Theorem~\ref{d19:thm:A}.
\par
\end{samepage}
\subsection{Proof of the reduction}\label{d19:sec:embed}
\begin{proof}[Proof of Corollary~\ref{d19:cor:main}]
Use the angle reduction of Section~\ref{d19:sec:angle}. Take perpendicular vectors $w,r_0$ of squares 6 and 2 in a two-dimensional space perpendicular to the 19-dimensional span of the $v_i$, and put
\[
x_i^\pm=w/2\pm v_i/\sqrt2.
\]
Their squares are four, they pair to three with $w$, and they are perpendicular to $r_0$. The two partners sum to $w$; other pairings are 1 or 2. Lemma~\ref{lem:paired-sections} with $q=5$ and $T=\Z r_0$ embeds the paired points and the marking in an even unimodular lattice $P$ of rank 24. The span $V$ of the centered points cannot contain a norm-two vector $s\in P$, since then $s/\sqrt2$ is a unit vector in $L^*$, by $(s/\sqrt2,v_i)=(s,x_i^+)\in\Z$, contradicting Theorem~\ref{d19:thm:A}. The 146 points $x_i^\pm$ now contradict Theorem~\ref{d19:thm:B}. The lower bound is the construction cited in Section~\ref{d19:sec:statement}.
\end{proof}
\begin{samepage}
\subsection{Marked Niemeier sections}\label{d19:sec:niemeier}
For the proof of Theorem~\ref{d19:thm:B}, suppose there are at least 73 pairs in the set $\mathcal X$. Choose exactly 73 pairs and redefine $V$ as their centered span. For distinct nonpartners $x,y\in\mathcal X$, integrality and Cauchy--Schwarz imply $(x,y)\in\{0,1,2,3\}$: the range of the pairing is $[-1,4]$ and the extreme values correspond to partners and equality. If $(x,y)=3$, then $x-y$ is a norm-two vector in $V$. If $(x,y)=0$, the vector $x+y-w$ is a norm-two vector in $V$. Both are forbidden, so the remaining pairings are 1 or 2 \cite[\S2.5]{Deg}. Thus $v_i=\sqrt2(x_i-w/2)$ for one representative per pair satisfy the hypotheses of Theorem~\ref{d19:thm:A}, and in particular they span dimension 19. For this section only the pairings, rank, and projection bounds of Section~\ref{d19:sec:angle} are used; there is no need to invoke Theorem~\ref{d19:thm:A} again.
\par
\end{samepage}
\begin{lemma}[Common height]\label{d19:lem:height}
There exists $t\in(\R w\oplus V)^\perp$ such that $(t,s)\ne0$ for every norm-two vector $s\in P\cap w^\perp$. All points of $\mathcal X$ have height $(t,x)=0$.
\end{lemma}
\begin{proof}A root perpendicular to $w$ cannot belong to $\R w\oplus V$, since then it belongs to $V$. Each forbidden pairing defines a proper hyperplane in the indicated complement. Their finite union cannot fill it. The last assertion follows from $x\in w/2+V$. This is the same common-height argument as in Lemma~\ref{d18:c:order}.\end{proof}
We use the classification of Niemeier lattices \cite{Niemeier,CS}. The Leech lattice is impossible because of the specified root $r_0$. The remaining 23 types are treated in the following three subsections. For any root subsystem perpendicular to both $w$ and $r_0$, we can apply its Weyl group to make $t$ positive on a chosen positive system. Thus differences in the nonzero positive-root cone cannot occur between selected points. We call a choice of a norm-six vector and a perpendicular root a marking.
\subsubsection{The type \texorpdfstring{$A_1^{24}$}{A1 (24 copies)}}\label{d19:sec:golay2}
Use the model $P=\frac1{\sqrt2}\{z\in\Z^{24}:z\bmod2\in G_{24}\}$, where $G_{24}$ is the extended binary Golay code. We use the weight set $\{0,8,12,16,24\}$, the octad Steiner system $S(5,8,24)$, and the fact that distinct octads meet in 0, 2 or 4 points \cite[Chapters~11 and~16]{CS}. Write $b=\sqrt2w$, $z=\sqrt2x$, and $v=z-b/2$ for the line representatives. Then
\[
b^2=12,\quad z^2=8,\quad (z,b)=6,\quad v^2=5.
\]
The marked root is in a coordinate $m$ with $b_m=z_m=0$. Arbitrary coordinate sign changes preserve the code model; normalize nonzero $b_i$ to be positive, and make the height coefficients positive at the zero coordinates of $b$. Such coefficients are nonzero by Lemma~\ref{d19:lem:height}. The possible nonzero entries of $b$, by its norm and the code weights, are three 2's, eight 1's and one 2, or twelve 1's. Every norm-eight $z$ is either two signed 2's or eight signed 1's on an octad.
\paragraph{Three 2's.} Let $T$ be their support. The inner product condition forces $z$ to be supported on an octad containing $T$, with all entries on $T$ equal to one. The number of such octads avoiding $m$ is $\lambda_3-\lambda_4=21-5=16$, using the Steiner system parameters. For each, the centered axis has five signed unit coordinates in the complement of $T$ in the octad. The common-height equation partitions five positive coefficients into equal-sum parts, up to complement. There are at most four such partitions: a singleton partition excludes every other one, and otherwise choose the side of size two. Two disjoint such two-sets would already use the entire weight, leaving a positive fifth coefficient, so the two-sets intersect pairwise. An intersecting family of two-sets on five points has size at most four. The bound is $16\cdot4=64$ axes.
\paragraph{Eight 1's and one 2.} Let $O$ be the octad of 1's, $a$ the position of 2. There are at most eight axes with $v_a\ne0$: orient their representatives with $z_a=0$ (the partner then has $z_a=2$), so $z$ must be the octad $O$ with just one negative entry to pair to six with $b$. A vector with two signed 2's and $z_a=0$ cannot pair to six. Lemma~\ref{d19:lem:projection} therefore eliminates all these axes. The other axes have $v_a=0$, and $z$ is supported on $T\cup\{a\}\cup B$, with $|T|=4$, $T\subset O$, $|B|=3$, $B\subset O^c\setminus\{a,m\}$. The inner product condition forces entries on $T$ and $a$ to be one. In particular
\begin{equation}\label{d19:eq:binregular}
v=s_T+\sum_{i\in B}\eps_i e_i,\qquad s_T=\mathbf1_T-\tfrac12\mathbf1_O.
\end{equation}
Here each $B$ corresponds to exactly one complementary pair $\{T,O\setminus T\}$, if it occurs. Indeed adding the octad $O$ in the code replaces $T$ by its complement; choosing either side, any two four-subsets of $O$ can be made to intersect in at least two points. Two different partitions with the same $B$ would then violate the octad intersection bound. Similarly two different triples $B,B'$ sharing at least two points would, by choosing the sides in $O$, give octads sharing at least five points.
Before removing $m$, the triples $B$ on $O^c\setminus\{a\}$ form an $S(2,3,15)$. To count triples containing a given $b'\ne a$ outside $O$, use the unique octad through $\{a,b'\}\cup U$ for each $U\in\binom O3$. Such an octad must meet $O$ in four points, so is counted four times; the dual octad for its $B$ contributes another four. The resulting number is $\binom83/8=7$. Together with the intersection property, this gives the Steiner system. After deleting $m$, there are 28 triples on 14 points, each point of degree six.
For a fixed $B$ orient the axis by a fixed choice of $T$. Equation \eqref{d19:eq:binregular} and positive height coefficients in $B$ imply that the allowed sign choices form an antichain of $2^B$, of size at most three \cite{Sperner}. To obtain 73 axes, at least $73-2\cdot28=17$ triples must attain three. Every exterior coordinate must occur, since otherwise the six triples through it are missing and at most $3(28-6)=66$ axes occur. Each exterior coordinate occurs on at least 14 axes by Lemma~\ref{d19:lem:projection}; hence at least two of its six triples must attain three. In a three-element antichain of size three all subsets are singletons or all are doubletons: the presence of a singleton and a disjoint doubleton gives a maximal antichain of size two, while the empty and full sets cannot occur. Thus the difference vectors on a saturated triple span its two-dimensional sum-zero space. The hypergraph of saturated triples is connected: its minimum degree is two and it is linear (intersections at most one), so every component has at least five vertices. A disconnected hypergraph on 14 vertices splits into sizes $5+9$, $6+8$ or $7+7$; the pair bound $\lfloor\binom h2/3\rfloor$ on $h$ vertices gives at most 15, 14 or 14 triples, contradicting the 17. Consequently $V$ contains the whole 13-dimensional exterior sum-zero space.
There are at least $\lceil73/3\rceil=25$ used triples, which give distinct partitions of $O$: the same partition for different triples would give octads sharing $T\cup\{a\}$, five points. The corresponding $s_T$ span the seven-dimensional sum-zero space on $O$. To prove this, a nonzero normal $u$ of sum zero would yield at least 50 four-subsets $T$ of zero $u$-sum. Choose $i,j$ with $u_i\ne u_j$. Of the 40 four-subsets containing exactly one of $i,j$, at most 20 can have zero sum, by swapping $i,j$. There are just 30 containing both or neither. Equality throughout would imply all four-subsets of the other six coordinates have zero sum, hence all those six entries are zero. Then $u_i=-u_j\ne0$, contradicting the existence of 20 zero-sum sets containing exactly one. The projection of $V$ onto $O$ thus has rank seven. Combined with the exterior subspace, $\dim V\ge20$, impossible.
\paragraph{Twelve 1's.} Let $D$ be their support, a dodecad. Every axis has form $v=s_T+\eps_i e_i+\eps_j e_j$ where $|T|=6$, $T\subset D$ and $E=\{i,j\}\subset D^c\setminus\{m\}$. As before, $D$ in the code replaces $T$ by its complement. For a fixed $E$ there is at most one such partition: distinct octads would meet in at least five points after choosing the sides ($\max(k,6-k)+2\ge5$). Orient by a side, and apply the height antichain argument on two exterior coordinates to get at most two axes per $E$. If $a$ exterior coordinates occur, Lemma~\ref{d19:lem:projection} yields $14a\le146$; also $73\le2\binom a2$. Thus $a=10$. Make a graph on these ten points, with edges the pairs giving two axes. At any vertex at least 14 axes occur and there are at most nine distinct pairs, so the graph has minimum degree at least five and is connected. Each such pair has sign choices $(+,-)$ and $(-,+)$ by the antichain property, giving the difference $e_i-e_j$ in $V$. Hence the nine-dimensional exterior sum-zero space lies in $V$.
There are $q\ge37$ used pairs. Let $M_D$ be the Gram matrix of their projected $s_T$. It has diagonal 3; its off-diagonals vanish when the pairs intersect and have absolute value one when disjoint. Indeed the relevant octads and their partners intersect in 0, 2 or 4 points, which forces $|T\cap T'|=3$ in the first case and 2 or 4 in the second. If $d_i$ are the degrees of the used-pair graph, then
\[
\tr M_D=3q,\quad \tr M_D^2=q^2+10q-\sum_i d_i^2\le3q^2/5+10q.
\]
It follows that $\rk M_D\ge9/(3/5+10/q)>10$, hence at least 11. Again $\dim V\ge9+11=20$, a contradiction. This finishes $A_1^{24}$.
\subsubsection{The types \texorpdfstring{$A_3^8$ and $A_2^{12}$}{A3 (eight copies) and A2 (twelve copies)}}\label{d19:sec:niemeier-compute}
Here we use the explicit code models \cite[Table~16.1]{CS} and the companion program \texttt{niemeier\_19\_certificate.py}. Write $A_{m-1}=\{x\in\Z^m:\sum x_i=0\}$, for $m=4$ or 3. The dual has coordinates with common fractional part in multiples of $1/m$. Let $C$ be a glue code over $\Z/m$ on $n=24/(m-1)$ factors. A lattice vector is described by $n$ such local dual vectors whose residue classes lie in $C$. For $m=4$, take the kernel of the matrix
\[
H_4=\begin{pmatrix}3&3&2&3&1&0&0&0\\3&0&3&2&3&1&0&0\\3&0&0&3&2&3&1&0\\3&0&0&0&3&2&3&1\end{pmatrix}
\quad\bmod4.
\]
For $m=3$ take the row code of $(I_6\mid P_3)$ over $\F_3$, where
\[
P_3=\begin{pmatrix}0&1&1&1&1&1\\1&0&1&1&2&2\\1&1&0&2&1&2\\1&1&2&0&2&1\\1&2&1&2&0&1\\1&2&2&1&1&0\end{pmatrix}.
\]
To verify these models directly, the code enumerator checks that $|C|=m^{n/2}$, each glue class has even minimum norm $\sum_i c_i(m-c_i)/m$, and every nonzero class has minimum at least four. This establishes an even unimodular extension with exactly the prescribed roots, by polarization and determinant. It is the required Niemeier lattice by classification. The ternary code is the extended ternary Golay code, in particular of minimum nonzero weight six, with weights 6,9,12. These latter properties also follow immediately by enumeration of its 729 words.
\paragraph{Complete marking and shell enumeration.}
Multiply all coordinates by $m$ to use integers. A local dual vector is an integer $m$-tuple summing to zero with all entries congruent modulo $m$. Enumerate these tuples by looping through the first $m-1$ coordinates in the norm-bounded box and determining the last by the sum. To enumerate $w$, loop over glue words, enumerate local tuples of norm at most $6m^2$, and impose the total norm $6m^2$. Sort each local tuple in decreasing order by Weyl symmetry. A root perpendicular to $w$ arises from a repeated coordinate value in one factor; permutation of the positions with that value makes any chosen pair equivalent. This gives all markings. The code verifies a small list of signed permutations of the factors preserving the glue code, then reduces the markings by their generated subgroup. Since only verified automorphisms are used, no assumption of the full automorphism group is needed.
For each orbit representative, enumerate all $x$ by the same local procedure with total norm $4m^2$, dot product with $w$ equal to $3m^2$, and the marked coordinates equal. The recursion prunes only by minimum remaining norm and Cauchy--Schwarz. These are necessary conditions and cannot discard a valid completion. The enumeration checks closure under $x\mapsto w-x$ in the integer model. Form centered axes $y=2x-w$, taking one from each antipodal pair. In unscaled coordinates they have squared norm 10 and desired mutual absolute pairing 2.
\paragraph{Compatibility and local cliques.}
In each component, form the ordered free blocks of coordinates with equal $w$-value, leaving the two marked positions fixed. Their difference roots preserve the marking. For the height vector of Lemma~\ref{d19:lem:height}, order its entries decreasingly in each such block, after applying permutations within the blocks. A vector lies in the resulting positive-root cone precisely when it is zero on the fixed positions, sums to zero in each block, and every initial partial sum in every ordered free block is nonnegative. Indeed the simple roots in a block are consecutive differences of coordinate unit vectors. Thus no two distinct selected axes can have $y+y'$ or $y-y'$ in the nonzero cone or its negative.
On all candidate axes first form the graph joining those with forbidden absolute pairing (not 2), and partition it into its connected components $C_j$. Within $C_j$ form the compatibility graph $\Gamma_j$ by imposing both the correct absolute pairing and the two cone tests. Any selected subset of $C_j$ is a clique, so $\sum_j\omega(\Gamma_j)$ bounds the number of axes. We need not impose compatibility between different components for this upper bound. The maximum-clique routine is exact: recursively enumerate vertices with their still-possible neighbors, and prune only using a proper greedy coloring of the remaining graph, whose number of colors bounds its clique number. For coloring bounds in clique search, see \cite{TomitaSeki}. If a linear functional has nonzero values on a subset whose corresponding sum of clique bounds is at most 13, Lemma~\ref{d19:lem:projection} permits deleting that entire subset in a hypothetical 73-axis solution. The program tests coordinate and difference functionals within factors and a few more rational directions, deleting only when this condition holds. Each such test remains valid after previous deletions. Table~\ref{d19:table:small-markings} summarizes the marking coverage.
\begin{table}[htbp]
\centering\small
\caption{Computed and analytic marking records for the two small glue models.}\label{d19:table:small-markings}
\begin{tabular}{lrr}\toprule
 & $A_3^8$ & $A_2^{12}$\\\midrule
Weyl-normalized marking records&13192&43980\\
Verified subgroup orbits&29&88\\
Computed bound at most 72 axes&29&81\\
Analytic cases below&0&7\\\bottomrule
\end{tabular}
\end{table}
All arithmetic in the enumeration is integer arithmetic; NumPy arrays hold the small scaled coordinates and inner products exactly. The two largest-component cases flagged previously in $A_3^8$ have 472 and 600 signed candidates; the clique-sum bounds for the particular markings of interest are 68 and 67 axes, respectively. The complete marking loop, not just those two cases, is what ensures coverage. The seven analytic cases all have the same squared-norm shape of $w$: the marked $A_2$ component has norm $8/3$, five others have norm $2/3$, and six have norm zero.

\phantomsection\label{d19:sec:K6}\paragraph{The analytic ternary Golay cases.}
Let $S$ be the six nonzero components of $w$, $j\in S$ the marked one, and $Z=S^c$. In the marked component, $w_j=-2\lambda$ with $\lambda^2=2/3$ and $x_j=k\lambda$ by perpendicularity to the marked root. The squared norms of $x$ and $w-x$ force $k=-2,-1,0$. By symmetry, to exclude $k=-2,0$ it suffices to exclude zero. If $k=0$ and the glue word of $x$ is zero, then $x$ is the sum of roots in two components (the $A_2$ root lattice has no square-four vector). Its pairing with $w$ is at most two. If the glue word is nonzero, then minimum weight six and $x^2=4$ force exactly six components of minimal norm $2/3$. The pairing $(x,w)=3$ would require all five nonmarked nonzero components of $w$ to be used: four identical weights and one opposite-class weight of pairing $1/3$. The glue words would differ at just three coordinates (the marked one, the opposite-class one, and an outside coordinate), violating distance six. Thus $x_j=w_j/2$.
Write $c,d$ for the glue words of $w,x$. Both $d$ and $c-d$ have weight six. The marked coordinate is $d_j=-c_j$. Among $S\setminus\{j\}$ let $a,b,z$ count coordinates with $d_i=c_i,-c_i,0$ respectively. Comparing the weights of $d$ and $c-d$ gives $a=z$, $b=5-2z$. If $z=1$, the word $c+d$ has weight three, impossible. Thus $z=0$ or 2. If $z=0$, $d=-c$ and the centered vectors have zero projections on $Z$. The five internal choices each give two opposite centered directions: the two minimal weights in the opposite class that pair to $1/3$ with the chosen weight of $w$. The positive-height partition argument from Section~\ref{d19:sec:golay2} gives at most four axes here.
If $z=2$, $d$ has exactly two nonzero coordinates in $Z$. For each $E\in\binom Z2$, the shortened code $C_E$ supported on $S\cup E$ has dimension exactly two. Its dimension is at least two by the four zero-coordinate constraints; every nonzero word has weight six (length at most eight), so a three-dimensional subcode would have total weight $26\cdot6=156$, exceeding the coordinate bound $8\cdot18=144$. There are four projective directions in $C_E$: the direction of $c$, the kernel of coordinate $j$, and two more, whose vectors with coordinate $-c_j$ are $d,c-d$. Each of the latter has outside support $E$. Indeed a nonzero word with that fixed coordinate has either $z=0$ or 2 by the preceding comparison when both $d$ and $c-d$ have weight six; the latter is guaranteed here because neither is zero and the length is at most eight. The case $z=0$ is exactly $-c$, already accounted for. Hence each $E$ corresponds to exactly one pair of dual glue orbits, and there are 15 pair cells indexed by the edges of $K_6$.
For a fixed orientation of a glue orbit in a pair cell, the possible $x$ form a product of three local sets: two weights at the unique internal opposite-class coordinate, and three minimal weights at each of the two exterior coordinates. All other coordinates are fixed by the pairing and norm constraints. In each local set the differences are roots perpendicular to $w$; the height order is strict and a chain. Thus common height selects an antichain in $[2]\times[3]\times[3]$. Its width is five, by the symmetric chain decomposition of a product of chains \cite{dBvTK} (or direct enumeration of its rank levels). Let $n_E\le5$ count axes in the cell; the special orbit has $n_0\le4$. Since $73=n_0+\sum_E n_E$, some $n_E=5$.
If $E=\{u,v\}$ has $n_E=5$, both exterior projections have rank two. Otherwise, since the three minimal weights in a fixed nonzero class are pairwise noncollinear, at most one of them could be used at the rank-one component; the product reduces to $[2]\times[3]$, of width two. By Lemma~\ref{d19:lem:projection}, a rank-two projection onto factor $u$ requires at least 25 active axes. Only its five incident cells can contribute, each with capacity five, so all five must be saturated. The same applies to $v$ and to every newly saturated edge. Since the line graph of $K_6$ is connected, all fifteen cells become saturated, giving at least 75 axes, impossible. This excludes all seven analytic markings.
\subsubsection{Independent verification of the twenty remaining root-system types}\label{d19:sec:import}
We give a direct computation for the twenty remaining root-system types. Degtyarev's work~\cite{Deg}, including the \href{https://archive.mpim-bonn.mpg.de/id/eprint/851/}{full preprint} (MPIM Preprint Series 2019, no.~51), provides a related approach to marked sections of Niemeier lattices. Our computation derives the needed bounds from the lattice models and shell enumerations described below. The driver \texttt{verify\_niemeier20.py} runs the following checks.

\paragraph{Construction of a lattice model.}
For the specified root system write $Q=Q_1\perp\cdots\perp Q_k$. Let $C_j$ be the Cartan matrix of the irreducible factor $Q_j$, $d_j=\det C_j$, and $K_j=d_j C_j^{-1}$. An element of $Q_j^*$ is represented by the integer column $u$ of pairings with the simple roots. Its coordinates in the simple-root basis, its squared norm, and its discriminant class are, respectively,
\[
 C_j^{-1}u,\qquad d_j^{-1}u^{\mathsf T}K_ju,\qquad K_ju\bmod d_j\Z^{\rk Q_j}.
\]
Here the class representation is injective on $Q_j^*/Q_j$ and addition is componentwise residue addition. All dominant pairing columns with squared norm at most seven are generated recursively by assigning nonnegative coordinates. Every entry of $K_j$ is positive, so a partial column exceeding the norm budget cannot extend to a feasible column. The script checks that these lists contain all $d_j$ discriminant classes. For each class its minimum squared norm is therefore exactly the minimum in the dominant list: every vector can be moved to the dominant chamber without changing its norm or class, and any shorter vector would occur in the same list.

Form the finite product of discriminant groups. Construct a subgroup $H$ of order $\sqrt{\det Q}$ all of whose classes have even minimum squared norm, and whose nonzero classes have minimum at least four. The code starts with the zero subgroup and adjoins generators, checking all cosets in each larger subgroup, until it reaches the required order. Closure under addition and the stated norm conditions are checked. The union of these classes in $Q^*$ is even, since changing a representative by a root-lattice vector changes its square by an even integer. Integrality follows by polarization. The order gives determinant one, and the minimum condition ensures that no new roots are present. By the uniqueness part of the Niemeier classification \cite{Niemeier,CS}, this model is the required lattice. No choice of a special glue embedding is assumed.

\paragraph{Enumeration of markings.}
We enumerate every possible $w$ up to the product of the Weyl groups by combining the dominant local lists, requiring the glue tuple to lie in $H$ and the sum of the local squared norms to be six. Factor transpositions are used to reduce the list only when the components have the same type and the transposition preserves $H$; each such transposition is tested on every element of $H$. It is not necessary to find the full lattice automorphism group. For a dominant $w$, the roots perpendicular to its projection in a component form the parabolic root subsystem on the simple roots of zero pairing. Each irreducible component of this subsystem has a transitive Weyl group action on its roots. Thus every possible marked root is equivalent, under an isometry fixing $w$, to one of the simple roots chosen from the connected components of the zero-pairing subdiagram. The program first tries the unmarked section, and only enumerates these marked representatives if its unmarked bound exceeds 144. The unmarked check is valid for every marking because it merely drops the extra perpendicularity restriction.

\paragraph{Local orbits and the first-pass bound.}
The local dual-lattice shell of squared norm at most four is obtained by applying the simple reflections to every relevant dominant column, stopping when no new column appears. In pairing coordinates, reflection in a root with simple-root coefficient column $a$ is the map $u\mapsto u-(a^{\mathsf T}u)C_ja$. For a marking, retain only local choices with zero pairing with its marked root in that component. A local candidate $x_j$ with squared norm $n$ and pairing $p$ with $w_j$ may be discarded if
\[
 (3-p)^2>(4-n)(6-w_j^2),
\]
since the other components would violate Cauchy--Schwarz. This is only a necessary test; the exact total norm and pairing are imposed at assembly.

For each component, let $R_j^+$ consist of its positive roots perpendicular to $w$ and, when imposed, to the marked root. They are obtained by filtering the enumerated norm-two roots. The indecomposable positive roots (those not expressible as a sum of two positive roots of this subsystem) give its simple-root basis. The local candidates are grouped by their norm, discriminant class and orbit under the reflection group of this subsystem, by reflecting to its dominant chamber. Each group has constant pairing with $w_j$. The admissible common-height vector can be moved into these positive chambers simultaneously by an isometry that fixes the prescribed data. It follows that selected differences cannot belong to the nonzero cone generated by these positive roots.

For a local orbit $\mathcal O_j$, let $a_j=|\mathcal O_j|$. Its compatibility graph joins two distinct elements if their squared difference is in $\{4,6,10\}$ and their difference belongs to neither the positive cone nor its negative. The allowed squared differences follow from the original pairings (including the antipodal partner). Let $b_j$ be any upper bound on this graph's clique number. The program uses exact clique recursion for groups of size at most 120, a proper greedy coloring for groups of size 121 through 500, and the trivial bound $a_j$ for larger groups. All three give upper bounds. For a product cell $\mathcal O=\prod_j\mathcal O_j$ in the global shell, fixing the coordinates other than $j$ leaves a clique in the local graph. Consequently
\begin{equation}\label{d19:eq:localfiber}
 |\mathcal X\cap\mathcal O|\le \min_j b_j\prod_{i\ne j}a_i.
\end{equation}
For each tuple in $H$ we assemble all products of local groups having that tuple of discriminant classes, total squared norm four and total pairing three. The products partition the shell. Summing \eqref{d19:eq:localfiber} proves the first-pass bound. The recursion that assembles these products discards a branch only when the minimum possible squared norm of its remaining factors is too large.

\paragraph{Second pass.}
For each marking whose first-pass upper bound exceeds 144, we enumerate its entire signed shell, now keeping each actual vector. The checker compares this shell's cardinality with the independent product count from the first pass, checks distinctness and closure under $x\mapsto w-x$, and chooses one representative $y=2x-w$ per antipodal pair. These representatives have square ten. Define $A_0$ on the resulting axes by absolute pairing two. Let $C_1,\ldots,C_s$ be the connected components of its complement graph (ignoring loops). For each component retain only the edges of $A_0$ for which neither $y+y'$ nor $y-y'$ belongs to the positive cone or its negative. These are necessary compatibility conditions by the common-height lemma. A selected system is a clique on its restriction to each $C_i$, so the sum of the clique upper bounds is an upper bound on the number of axes. The algorithm uses exact maximum-clique recursion for $|C_i|\le80$, and the minimum of two proper greedy-coloring upper bounds otherwise. It suffices that each coloring is proper, not that it is optimal.

For clarity, all cone comparisons are exact. Let $R$ be the matrix of simple-root coefficient rows for the perpendicular subsystem in one factor, and let $M=RC_jR^{\mathsf T}$. For a vector with pairing row $u$ its projection onto this subsystem has coefficient row $(uR^{\mathsf T})M^{-1}$. Multiplying by $\det M$ gives integer comparison coordinates. In a local orbit all differences lie in the subsystem span. In the global test we additionally compare the complementary coordinates: for a global scale $d=\operatorname{lcm}_j d_j$ and scaled simple-root coordinates $B=duC_j^{-1}$, the integer vector
\[
 (\det M)B-d(uR^{\mathsf T}\operatorname{adj}M)R
\]
is zero exactly for vectors in the subsystem span. Thus testing its sum or difference for zero, together with the sign of the comparison coordinates in every factor, tests precisely the required cone membership. The empty subsystem has zero cone. All Gram products, norm tests and cone tests are computed using integers. The implementation checks $d_j\le25$, $\max K_j<200$, root coefficients of absolute value at most 8, $\det M\le 2^{24}$ and adjugate entries of absolute value at most $2^{23}$. These bounds, the rank limit 24, and the pairing bound $|u_i|\le4$ for the centered vectors ensure that the signed 64-bit integer products used here do not overflow.

The complete first pass gives 540 exceptional marked records. Every one is excluded by the second-pass clique upper bound, without the optional projection-support deletions. The exact final bounds on signed points, obtained by taking the maximum of all passing first-pass bounds and twice the second-pass axis bound, are in Table~\ref{d19:tab:independent-niemeier}.
\begin{table}[ht]
\centering\small
\begin{tabular}{lrr}\toprule
Root system & Second-pass records & Signed-point bound\\\midrule
$D_{24}$ & 2 & 40\\
$D_{16}E_8$ & 3 & 90\\
$E_8^3$ & 0 & 59\\
$A_{24}$ & 4 & 56\\
$D_{12}^2$ & 2 & 88\\
$A_{17}E_7$ & 2 & 78\\
$D_{10}E_7^2$ & 4 & 90\\
$A_{15}D_9$ & 0 & 82\\
$D_8^3$ & 0 & 116\\
$A_{12}^2$ & 0 & 83\\
$A_{11}D_7E_6$ & 0 & 108\\
$E_6^4$ & 0 & 104\\
$A_9^2D_6$ & 0 & 113\\
$D_6^4$ & 0 & 130\\
$A_8^3$ & 0 & 100\\
$A_7^2D_5^2$ & 0 & 136\\
$A_6^4$ & 0 & 132\\
$D_4^6$ & 510 & 144\\
$A_5^4D_4$ & 1 & 144\\
$A_4^6$ & 12 & 142\\
\bottomrule\end{tabular}
\caption{Independent enumeration of the twenty remaining root systems. All bounds are on signed lattice points, not axes.}\label{d19:tab:independent-niemeier}
\end{table}

\paragraph{Independent shell-count check.}
An extra check recomputes the local orbit sizes as the Weyl group order divided by the Weyl group order of the stabilizer, which is the parabolic group on the simple roots of zero pairing. For squared norms at most four these sizes are compared with the actual local enumeration. Combining the orbit-size counts over $H$ gives the full lattice shell counts at squared norms two, four, and six. For each of the twenty models they agree with the theta-series identity for a rank-24 even unimodular lattice \cite[Ch.~4]{CS},
\[
 \Theta_P(q)=E_4(q)^3+(|P_2|-720)\Delta(q),\qquad P_2=\{x\in P:x^2=2\}.
\]
In particular, writing $|P_2|=24h$,
\[
 |P_4|=196560-576h,\qquad |P_6|=16773120+6048h.
\]
This is a cross-check rather than a substitute for the completeness proof of the enumeration. The driver \texttt{verify\_niemeier20.py} reruns both passes and the shell-count check, verifies that every type and every second-pass record has been covered, and prints \texttt{ALL 20 TYPES VERIFIED}. All twenty bounds are at most 144. Combined with the three types handled in the preceding subsections and the exclusion of the Leech lattice by the marked root, this proves Theorem~\ref{d19:thm:B} by the stated finite checks. No numerical bound from a table of geometric sets is an input to this argument.
\subsection{Verification of the finite computations}\label{d19:sec:verification}
The mathematical part of the root-system argument identifies the necessary inequalities used by the finite calculation. The exact rational certificates verify infeasibility of the stated relaxations. The reproduction directory provides a driver for the coverage check, the mixed and homogeneous root-type certificates, the pure $A_1$ calculation, and both Niemeier computations. The driver compares the results with the accompanying verification logs. From the repository checkout, run
\begin{verbatim}
cd dimension-19/reproduction
python3 -m venv .venv
source .venv/bin/activate
python3 -m pip install -r requirements.txt
unset PYTHONOPTIMIZE
shasum -a 256 -c MANIFEST.sha256 && python3 -B audit/referee/replay.py
\end{verbatim}
Use Python 3.10 or newer with assertions enabled. The Niemeier computations use NumPy and NetworkX. The independent twenty-type enumeration also uses SymPy for exact matrix inverses. The small-model replay enumerates 88 marking orbits in $A_2^{12}$ and verifies the computed bounds for 81. It identifies the seven remaining orbits by the marked shape of $w$ and reports them as \texttt{ANALYTIC}. The \hyperref[d19:sec:K6]{analytic ternary Golay argument} excludes these seven using common height and rank-two projection support. The proof combines the verified finite checks with the geometric reductions and analytic cases above.

\section{Dimension twenty}\label{part:20}
\subsection{Statement and reduction to sections of rank 24}
Write $N(d)$ for the maximum cardinality of a real equiangular line system in $\R^d$. The exact value in dimension twenty is a classical small-dimensional problem. Taylor's construction \cite[p.~123]{Taylor} gives $N(20)\ge90$, while Greaves--Syatriadi--Yatsyna \cite[Theorem~5.5]{GSY1} proved $N(20)\le94$. We close this gap by combining a lattice embedding with combinatorial and energy estimates. The accompanying finite checks are part of the proof of the following theorem.
\begin{theorem}\label{d20:main}
$N(20)=90$.
\end{theorem}
The argument is organized as follows. We first reduce to norm-six sections of rank-24 even unimodular lattices. Section~\ref{d20:sec:orbit-energy} establishes a uniform orbit capacity bound, and Section~\ref{d20:sec:rooted} applies it to the nonbinary rooted types. Section~\ref{d20:sec:golay} handles the binary rooted type by design theory. In Section~\ref{d20:sec:leech}, the rootless case is reduced to a generalized Hamming weight \cite{Wei}, which is bounded by Golay incidence counting and a lattice-parity obstruction. Finally Section~\ref{d20:sec:lower} gives a direct count for the lower bound and assembles the proof. No assertion that every 91-line configuration extends to the Conway configuration is required.

For comparison, the stronger previously published bounds $N(19)\le74$ and $N(18)\le59$ are due to Greaves--Syatriadi--Yatsyna and Greaves--Syatriadi, respectively \cite{GSY1,GS}. Neither result is used in the proof. We first establish the weaker fixed-angle estimates that suffice here, and also give a root-system proof of the classical $1/3$ bound \cite[Theorem~4.5]{LS}.

\begin{lemma}\label{d20:thirdangle}
$N_{1/3}(20)\le38$.
\end{lemma}
\begin{proof}
Choose representatives $v_i$ of square 3, with nonpartner pairings $\pm1$, and a unit vector $w$ perpendicular to their span. Put $r=\sqrt2w$ and $x_i=(w+v_i)/\sqrt2$ \cite[\S\S4--6]{Yoshino}. The lattice generated by $r,x_i$ is even integral of rank at most 21, because its generators have square 2 and integral pairings. It is discrete: choose a real basis from the generators, whose integral span has dual containing all generators. Its norm-two vectors form a finite simply laced root system, since they span and reflections in them preserve the lattice. Both $x_i$ and $r-x_i$ are roots pairing to 1 with $r$. Every such root lies in the irreducible component $R$ containing $r$. Its rank $a$ is at most 21.

The operator $\sum_{\beta\in R}\beta\beta^{\mathsf T}$ on the component span is invariant under the Weyl group, and hence scalar by irreducibility; its trace gives the scalar $2|R|/a$. Consequently $\sum_{\beta\in R}(\beta,r)^2=4|R|/a$. The two roots $\pm r$ contribute 8, and all other pairings are $0,\pm1$. Negation interchanges the pairings 1 and $-1$. Thus exactly $2|R|/a-4$ roots pair to 1 with $r$. They are paired without fixed points by $\beta\mapsto r-\beta$, and the original lines determine distinct such pairs. The number of lines is at most $|R|/a-2$. For $A_a,D_a,E_6,E_7,E_8$ these bounds are respectively $a-1,2a-4,10,16,28$, by the usual root counts \cite[Ch.~VI]{Bourbaki}. With $a\le21$ the maximum is at most 38.
\end{proof}

\begin{lemma}\label{d20:sixty}
$N_{1/5}(18)\le60$ and $N_{1/5}(19)\le76$.
\end{lemma}
\begin{proof}
For $m$ representatives of square five and pairings $\pm1$ spanning dimension $d$, the Gram trace inequality \cite[Lemma~6.1]{vanLintSeidel} gives $25m^2\le d(m^2+24m)$. In particular the relative bounds in dimensions 18 and 19 are 61 and 76. Suppose 61 lines occur in dimension 18. The same inequality forces their span to have dimension 18. Write their Gram matrix as $5I+S$. The eigenvalue $-5$ of the integral Seidel matrix $S$ has multiplicity 43. Its other 18 eigenvalues $\lambda_i$ have sum 215 and sum of squares 2585. Setting $t_i=\lambda_i-12$ gives
\[
\sum_{i=1}^{18}t_i=-1,\qquad \sum_{i=1}^{18}t_i^2=17.
\]
The $t_i$ are the roots of a monic integral polynomial. If all are nonzero, the arithmetic--geometric mean inequality on their squares gives $\sum t_i^2\ge18$. Thus 12 is an eigenvalue. The characteristic polynomial of $S$ modulo 2 is $x(x+1)^{60}$, because $S\equiv J-I$. Hence this even integer eigenvalue has multiplicity exactly one. The 17 remaining $t_i$ have nonzero integral product; the same inequality, now with equality, forces each to be $\pm1$. Their sum is $-1$, so the spectrum of $S$ is $[-5]^{43},[11]^9,[12]^1,[13]^8$.

The matrix $p(S)$ for $p(x)=-(x+5)(x-11)(x-13)$ is positive semidefinite of trace 17. But each diagonal entry equals $425-(S^3)_{ii}$, since $S_{ii}=0$ and $(S^2)_{ii}=60$. Each $(S^3)_{ii}$ is an even integer, by pairing terms indexed by $(j,k)$ and $(k,j)$. Therefore all 61 diagonal entries of $p(S)$ are nonnegative odd integers, giving trace at least 61, a contradiction. This proves the fixed-angle estimate by the Seidel-matrix argument of Greaves \cite[Remarks~4.7 and~6.6]{Greaves18}, included here to make the numerical input explicit.
\end{proof}

\begin{lemma}\label{d20:angle}
A set of 91 equiangular lines in $\R^{20}$ has common cosine $1/5$ and spans $\R^{20}$.
\end{lemma}
\begin{proof}
The cosine $\alpha$ lies in $(0,1)$, since 91 distinct mutually orthogonal lines cannot occur in dimension 20. Write the Gram matrix of unit representatives as $I+\alpha S$, with $S$ an integral Seidel matrix. The eigenvalue $-1/\alpha$ of $S$ has multiplicity at least 71. The degree of its irreducible polynomial over $\Q$ is at most $91/71$, so $1/\alpha=q$ is a positive integer. If $q$ were even, the characteristic polynomial of $S$ modulo 2 would be divisible by $t^{71}$, contrary to $S\equiv J-I\pmod2$ and $\det(tI-S)\equiv t(t+1)^{90}$. Thus $q$ is odd \cite[Theorem~3.4]{LS}. The trace inequality gives $91\le20(1+90/q^2)$, so $q^2\le1800/71$. Since $q>1$, $q=3$ or $5$. Lemma~\ref{d20:thirdangle} excludes $q=3$. Rescale the vectors to $v_i^2=5$; their other pairings are $\pm1$. With $F=\sum_i v_iv_i^{\mathsf T}$, we have $\tr F=455$, $\tr F^2=10465$, and $\rk F\ge455^2/10465=455/23>19$.
\end{proof}

Take $v_1,\dots,v_{91}$ as in Lemma~\ref{d20:angle}, and choose $w$ of square 6 orthogonal to their span. Put $x_i^\pm=w/2\pm v_i/\sqrt2$. These have square 4 and pairing 3 with $w$; partners sum to $w$, and nonpartners from different lines pair to 1 or 2. Lemma~\ref{lem:paired-sections} with $q=5$ and $T=0$ embeds this configuration in a rank-twenty-four even unimodular lattice $P$ \cite{Niemeier,CS}.

For any such $P$ and any $w\in P$ of square 6, define
\[
\mathcal F_P(w)=\{x\in P:x^2=4,\ \ip{x}{w}=3\}.
\]
We represent each antipodal pair $\{x,w-x\}$ by a centered vector $y=2x-w$, of square 10, or its negative \cite[\S2.5]{Deg}. Two different axes are compatible precisely when $|\ip{y}{y'}|=2$. The task is to exclude 91 mutually compatible axes whose span has dimension at most 20. Such a hypothetical set spans dimension 20 by Lemma~\ref{d20:angle}'s trace calculation. For a rank-$s$ linear map on its span, let $n(\rho)$ count the nonzero projections. Lemma~\ref{d20:sixty}, applied to the kernel of the map, gives
\begin{equation}\label{d20:support}
 n(\rho)\ge15\quad(s=1),\qquad n(\rho)\ge31\quad(s\ge2).
\end{equation}
All uses of projection deletion below are conditional on this hypothetical set.

\subsection{Parabolic orbits and the block-energy bound}\label{d20:sec:orbit-energy}
For this section normalize the axes by $v_x=\sqrt2(x-w/2)$, so $v_x^2=5$ and compatibility means absolute pairing one. Put
\[
R_0=\{r\in P:r^2=2,\ \ip{r}{w}=0\},\qquad W_0=W(R_0).
\]
The group $W_0$ fixes $w$ and preserves the section. It is useful to identify its orbits intrinsically, although the verifier will not construct a graph.
\begin{lemma}\label{d20:orbitgraph}
The connected components of the graph of incompatible axes in $\mathcal F_P(w)$ are exactly the $W_0$-orbits on axes.
\end{lemma}
\begin{proof}
For distinct nonpartner points $x,x'$ the possible pairings are $0,1,2,3$, by integrality and the strict Cauchy--Schwarz inequality for their centered vectors. If the pairing is three, $r=x-x'$ is a root perpendicular to $w$ and $\ip{x}{r}=1$, so reflection in $r$ takes $x$ to $x'$. If the pairing is zero, $r=x+x'-w$ is such a root with $\ip{x}{r}=1$; the reflection takes $x$ to $w-x'$. Conversely, for a root $r\in R_0$, the integral pairing $\ip{x}{r}$ has absolute value at most two. A reflection that changes the axis changes the point to one with pairing $4-\ip{x}{r}^2\in\{0,3\}$ with $x$. The reflection generators thus connect every orbit, and every incompatibility edge stays inside an orbit.
\end{proof}

Let $U_1,\ldots,U_k$ be the orthogonal spans of the irreducible components of $R_0$, of dimensions $d_j$. A signed orbit of centered vectors consists of vectors with fixed complementary part $z$ and with independent Weyl-orbit choices $p_j\in U_j$. Put $h=z^2$ and $e_j=p_j^2$, so $h+\sum e_j=5$. Omit the inactive components ($e_j=0$). The linear span of the orbit has dimension $r=\mathbf1_{h>0}+\sum_{e_j>0}d_j$: the Weyl orbit in each active irreducible component has average zero and spans that component, and its differences span it as well. The remaining fixed projection contributes just the line $\R z$, if nonzero. In particular $r\le23$.

A signed orbit may coincide with its negative, in which case $h=0$ and its axes are obtained by dividing its cardinality by two. Otherwise the two signed orbits give the same set of axes, and we choose representatives in one of them. Write $o_j$ for the local orbit sizes, $o=\prod_j o_j$, and $\epsilon=1$ for a signed orbit coinciding with its negative and zero otherwise. Thus the number of available axes in the cell is $o/(1+\epsilon)$.

\begin{lemma}[Block-energy and centroid bounds]\label{d20:blockenergy}
The capacity of this cell is at most
\begin{equation}\label{d20:cellcap}
 b=\min\left\{\frac{o}{1+\epsilon},\ \left\lfloor\frac{24}{h^2+\sum_{e_j>0}e_j^2/d_j-1}\right\rfloor,\ B_h\right\},
\end{equation}
where $B_h=\lfloor 4/(h-1)\rfloor$ when $h>1$, and $B_h=+\infty$ otherwise.
\end{lemma}
\begin{proof}
For $m$ selected axes, choose representatives in the signed orbit and let $F=\sum_{i=1}^m v_iv_i^{\mathsf T}$. Then $\tr F^2=m^2+24m$. Its orthogonal diagonal blocks on the complementary line and on $U_j$ have squared Frobenius norms at least $m^2h^2$ and $m^2e_j^2/d_j$, respectively. Hence
\[
m^2+24m\ge m^2\left(h^2+\sum_{e_j>0}e_j^2/d_j\right).
\]
The denominator in \eqref{d20:cellcap} is positive, since the parenthesized expression is at least $25/r>1$ by Cauchy--Schwarz. Also $m^2h\le\|\sum_i v_i\|^2\le m^2+4m$, giving the centroid bound. For a self-negative signed orbit $h=0$, so the orientation needed for this last inequality is irrelevant. The size bound is immediate.
\end{proof}

\begin{lemma}[Fiber bound]\label{d20:fiber}
When $e_j<4$, fixing every coordinate outside $U_j$ in an oriented cell allows at most $\lfloor4/(4-e_j)\rfloor$ distinct compatible vectors. Thus, provided $h>0$ or another root component is active, the capacity is also at most
\[
\left\lfloor\frac{1}{1+\epsilon}\frac{o}{o_j}\left\lfloor\frac4{4-e_j}\right\rfloor\right\rfloor.
\]
\end{lemma}
\begin{proof}
For $m'$ vectors in a fiber, their $U_j$ projections have square $e_j$, and distinct axes give pairings at most $e_j-4$ there. The square of their sum is therefore at most $m'[4-m'(4-e_j)]$, which is nonnegative. In the self-negative case include both orientations of each selected axis. The fixed nonzero complementary projection in the fiber prevents both signs of the same axis occurring there. Summation over the $o/o_j$ fibers and division by $1+\epsilon$ gives the assertion.
\end{proof}

\subsection{The nonbinary rooted lattices}\label{d20:sec:rooted}
We apply \eqref{d20:cellcap} to all cells and sum their capacities. The cells can be obtained using local parabolic-dominant pairings, without enumerating any global section or any full local Weyl orbit. The twenty general models are checked using only the block-energy and centroid bounds and orbit cardinalities.

\subsubsection{Local Cartan data and completeness}
Let $Q_i$ be a root factor, $C_i$ its Cartan matrix, $D_i=\det C_i$, and $K_i=D_i C_i^{-1}$. A vector in $Q_i^*$ has integral simple-root pairing column $u$, squared norm $u^{\mathsf T}K_i u/D_i$, and discriminant class represented by $K_i u\bmod D_i$. The verifier enumerates nonnegative pairing columns of square at most seven using the positive entries of $K_i$. At least one representative of each of the $D_i$ classes occurs, checked in the program; the minima therefore occur in this list. The supplied glue subgroup in \texttt{glue\_models.json} can be regenerated using the subgroup construction in Section~\ref{d19:sec:import}. It is checked for closure, order $\sqrt{\det Q}$, even minimum norms and minimum at least four outside the zero class. Its preimage is even integral and unimodular with precisely the stated roots. Uniqueness in the Niemeier classification identifies it with the required lattice.

For local marking $w_i$ with dominant pairing column $a$, the simple roots of the perpendicular subsystem are exactly the nodes $Z=\{j:a_j=0\}$: a positive root is a nonnegative linear combination of simple roots. For the point $x_i$ let $u$ be its pairing column. Choose it dominant on $Z$, with no sign restriction at the other nodes. In a connected component $I$ of the zero-node subdiagram, write $C_I$ for its Cartan matrix. The centered projection has squared norm
\[
e_I=2u_I^{\mathsf T}C_I^{-1}u_I.
\]
All $u_I\ge0$ with $e_I\le5$ are enumerated using the positive inverse Cartan entries. At nodes outside $Z$, the pairing with a simple root lies in $\{-2,-1,0,1,2\}$ since $x_i^2\le4$. Check the two local norm budgets $x_i^2,(w_i-x_i)^2\le4$ and the local centered square at most five. These tests include every possible global point. Every orbit has a unique parabolic-dominant representative. The local orbit size on each connected component is its Weyl-group order divided by the Weyl-group order of its zero-pairing subdiagram. This is the usual stabilizer formula for a closed Weyl chamber \cite[Ch.~VI]{Bourbaki}. The partner is tested by moving the column $a-u$ into the same parabolic dominant chamber by its simple reflections. In particular we know whether the local orbit is partner-fixed.

No full local Weyl orbit need be generated: all this is done by enumeration of nonnegative local pairing columns and finitely many entries at nonzero marking nodes. As a cross-check of the model and dominant marking lists, the norm-two, four, six shell cardinalities obtained by convolving local orbit-size histograms agree with $n_2$, $196560-24n_2$, $16773120+252n_2$ from the standard theta-series formula \cite[\S3.2, equation~(4)]{NT}.

\subsubsection{The weighted histogram formula}
A scalar marking profile consists of the glue word $c=(c_i)$ of $w$ and its local squared norms $q_i$, with $\sum q_i=6$. We do not choose a global marking vector. For each factor, local marking class and local square, prepare the histogram of possible local parabolic orbits of $x_i$. A signature records the class of $x_i$, the two squares $x_i^2,(w_i-x_i)^2$, the active triples $(d_I,e_I,o_I)$, and whether the local orbit is partner-fixed. The multiplicity of a signature is the number of distinct local parabolic-dominant columns giving it. If the marking class and square have more than one local dominant lift, retain its separate possible histograms.

For any choices $H_i$ of local histograms, sum over all glue words of $x$ and all tuples of signatures satisfying both total norm budgets $4,4$. Write $k_i$ for the signature multiplicities. For a tuple $\sigma$ compute $h=5-\sum e_I$, $o=\prod o_I$ and $\epsilon=1$ precisely when all the local orbits are partner-fixed. Let $b(\sigma)$ be \eqref{d20:cellcap}. The total upper bound is
\begin{equation}\label{d20:scalarformula}
\mathcal B(H_1,\ldots,H_t)=\frac12\sum_{\sigma}(1+\epsilon(\sigma))b(\sigma)\prod_i k_i.
\end{equation}
Indeed, a pair of distinct signed orbits has two equal contributions; a self-negative orbit gives one contribution with weight two before division. The formula thus sums one valid bound for each set of axes. This remains an upper bound if some actual histogram is replaced by the entrywise maximum of all its alternatives: every summand is nonnegative. It is linear in each individual histogram.

There is at most one ambiguous factor in any scalar norm-six profile. A shortest vector of a discriminant class in an irreducible simply laced root lattice has one dominant representative (also directly checked for the local data here). If the marking glue word is nonzero its minimum square is at least four, so at most one local square can exceed its class minimum by a positive even integer. If its glue word is zero, the norm-two choices in a factor have one Weyl orbit, and at most one factor can have square four or six. Thus the maximum of \eqref{d20:scalarformula} can be taken by a single local histogram maximization, with the remaining histograms fixed. This property is also explicitly checked on every scalar profile. The program does not generate global dominant marking vectors or reduce their Weyl orbits. It considers all allowed scalar norm budgets; all tests are arithmetic on local tables and signature multiplicities.

\subsubsection{Results for the twenty general models}
The program \texttt{verify\_scalar.py} checks the upper bounds in Table~\ref{d20:table:scalar-profiles}, using exact rational arithmetic and Lemma~\ref{d20:blockenergy}. They hold without the dimension-20 restriction. The table entries are maxima of computed upper bounds, not assertions that configurations attaining them exist.
\begin{table}[htbp]
\centering\small
\caption{Scalar profiles and upper bounds in the twenty general rooted models.}\label{d20:table:scalar-profiles}
\begin{tabular}{lrr@{\qquad}lrr}\toprule
Type&Scalar profiles&Bound&Type&Scalar profiles&Bound\\\midrule
$D_{24}$&2&31&$D_{16}E_8$&6&43\\
$E_8^3$&10&36&$A_{24}$&5&43\\
$D_{12}^2$&9&47&$A_{17}E_7$&13&51\\
$D_{10}E_7^2$&19&51&$A_{15}D_9$&16&53\\
$D_8^3$&29&56&$A_{12}^2$&24&56\\
$A_{11}D_7E_6$&39&58&$E_6^4$&52&57\\
$A_9^2D_6$&57&63&$D_6^4$&74&63\\
$A_8^3$&72&67&$A_7^2D_5^2$&120&70\\
$A_6^4$&164&72&$D_4^6$&344&84\\
$A_5^4D_4$&290&82&$A_4^6$&540&87\\
\bottomrule
\end{tabular}
\end{table}
There are 1885 scalar profiles; the entrywise-max envelope passes for 1464 of them, while 421 require the single-factor refinement. All refinements pass. No coloring, compatibility graph, projection deletion, rank search or global marking orbit list is involved for these twenty types.

\subsubsection{The two small glue models}
For $A_2^{12}$ use the ternary row code $(I_6\mid P_3)$, where
\[
P_3=\begin{pmatrix}
0&1&1&1&1&1\\1&0&1&1&2&2\\1&1&0&2&1&2\\1&1&2&0&2&1\\1&2&1&2&0&1\\1&2&2&1&1&0
\end{pmatrix}.
\]
For $A_3^8$ use the kernel over $\Z/4\Z$ of
\[
H_4=\begin{pmatrix}3&3&2&3&1&0&0&0\\3&0&3&2&3&1&0&0\\3&0&0&3&2&3&1&0\\3&0&0&0&3&2&3&1\end{pmatrix}.
\]
The code checks verify order $m^{n/2}$, even minimum squares $\sum c_i(m-c_i)/m$ and minimum at least four outside zero, for $(m,n)=(3,12),(4,8)$. This gives the two even unimodular glue extensions \cite[Table~16.1]{CS}. Use local coordinates in $A_{m-1}^*$ scaled by $m$, with integral entries congruent modulo $m$ and sum zero. The local marking can be sorted; its stabilizer permutes positions having equal marking entries. This allows the preceding parabolic data to be generated by local coordinate permutations alone. The program \texttt{verify\_small.py} checks all scalar norm profiles (3840 and 1256), reduces them by verified signed permutations preserving the code (22 and 10 scalar orbits), and tests every alternative of the at most one ambiguous local lift (24 and 13 tests, respectively). It never constructs global candidate vectors.

The additional bounds are Lemma~\ref{d20:fiber} and two small product bounds. In the self-negative orbit with five independent binary sign coordinates of norm one, orient every axis to even parity. Its distinct words have Hamming distance two; translating one word to zero gives an intersecting family of two-subsets of a five-set, of size at most four (a star or a triangle). Hence there are at most five axes. The other special self-negative orbit has one norm-one binary sign block and two squared-norm-two $A_3$ weight orbits of size six. Each latter orbit consists of three orthogonal antipodal directions: in ordinary coordinates it is the orbit of $(1,1,-1,-1)/\sqrt2$. Indeed, its unscaled dual-root projection has squared norm one, and in $A_3^*$ the only vectors of this norm are the permutations of $(1,1,-1,-1)/2$ (use coordinates congruent modulo one and of sum zero). Orient the first sign block positively, so the vectors have form $g+\sqrt2\varepsilon e_i+\sqrt2\eta f_j$, $1\le i,j\le3$. The unit directions here are orthogonal. The occupied supports $(i,j)$ form a bipartite graph with maximum degree two: at a shared endpoint the signs of distinct supports must be opposite. Each support has multiplicity at most two, as the two sign words must differ at exactly one coordinate. A doubled support has a degree-one endpoint, because at its differing coordinate no other support can have the required opposite sign to both. If there are $E$ occupied edges and $n_1,n_2$ vertices of degrees one and two, then $n_1+2n_2=2E$ and $n_1+n_2\le6$. The total multiplicity is at most $\min(2E,E+n_1)\le\min(2E,12-E)\le8$.

The projection tests can also be done on entire signed-orbit pairs. A factor's projected norm is constant on such a cell. For a rank-one projection bound, group its possible nonzero projections by the one-dimensional subspaces they span, using the local coordinate orbit alone. Fixing a projected vector in a cell removes that factor's active blocks and adds their energies to $h$ (its now fixed projection); the other local choices are unchanged. Apply the same capacity formulas to this slice. In a prescribed one-dimensional subspace, choose one of the two directions. In a self-negative cell orient each axis so that its nonzero projection points in that direction. In any other cell sum the slice capacities for both directions in the subspace. For each subspace, sum the slice bounds across the cells. The largest of these sums bounds the size of every rank-one nonzero support. The scalar pairing with the local marking is also constant up to sign, and its nonzero support is a union of whole cells.

By \eqref{d20:support}, for each scalar pairing with a local marking, delete its entire nonzero support if the sum of the capacities on that support is at most fourteen. For a factor projection, delete its entire nonzero support if its total capacity is at most fourteen, or if its rank-one support bounds are at most fourteen and its total capacity is at most thirty. If all rank-one support bounds are at most fourteen and the zero support has capacity at most ninety, the factor is forced to have rank at least two. For $k$ forced factors and at most $t$ such active factors per retained cell, counting incidences gives $31k\le91t$. The three-root marking of $A_2^{12}$ has $k=9,t=3$, impossible. All these tests are applied to local signature products and local projections, without pairwise compatibility matrices. Table~\ref{d20:table:small-models} gives the computed upper bounds grouped by nonzero component squared norms (norms in the table are scaled by $m^2$).
\begin{table}[htbp]
\centering\small
\caption{Bounds for the two small glue models grouped by scaled norm pattern.}\label{d20:table:small-models}
\begin{tabular}{cl r}\toprule
$m$&Scaled norm pattern&Bound\\\midrule
3&$(54)$&33\\
3&$(18,18,18)$&90\\
3&$(18,6^6)$&90\\
3&$(24,6^5)$&65\\
3&$(6^9)$&84\\\midrule
4&$(96)$&28\\
4&$(32,64)$&88\\
4&$(32,32,32)$&85\\
4&$(44,16,12^3)$&86\\
4&$(48,12^4)$&84\\
4&$(48,16^3)$&88\\
4&$(32,16^4)$&78\\
4&$(32,16,12^4)$&82\\
4&$(16^3,12^4)$&80\\
4&$(12^8)$&84\\
\bottomrule
\end{tabular}
\end{table}
Here powers in a tuple denote repeated entries. For example, $6^6$ in this table denotes six entries equal to six, not exponentiation. A bound obtained after projection deletion or incidence contradiction is conditional on the hypothetical 91-line system, not necessarily an unrestricted section bound. This completes all twenty-two nonbinary rooted models.

\subsection{The binary Golay model}\label{d20:sec:golay}
For $A_1^{24}$ use $P=\{z/\sqrt2:z\in\Z^{24},\ z\bmod2\in G_{24}\}$, where $G_{24}$ is the parity extension of the cyclic binary code generated by $X^{11}+X^9+X^7+X^6+X^5+X+1$ \cite[Definition~3.5.20]{NW}. The checker verifies its order 4096, weights $0,8,12,16,24$ and that the 759 octads form $S(5,8,24)$ \cite[Ch.~11]{CS}. Distinct octads meet in 0, 2 or 4 points, by linearity and the weight set. The lattice is even unimodular by doubly evenness and the code order, and has exactly the orthogonal coordinate roots. Write $b=\sqrt2 w$ and $z=\sqrt2 x$. We require $b^2=12,z^2=8,\ip{z}{b}=6$, with centered representatives $v=z-b/2$. Independent coordinate sign changes preserve this model. Thus take $b$ nonnegative: its possibilities are three entries 2; eight entries 1 and one entry 2; or twelve entries 1. A norm-eight $z$ is two signed 2s or eight signed 1s on an octad. We exclude all three cases without assuming that the span is root-free.

\subsubsection{Three entries 2}
Let $T$ be their support. Only octads containing $T$ with signs positive there can occur. Removing $T$ gives a $2$-$(21,5,1)$ design with 21 blocks, each point in five blocks. On each block the centered vectors are five sign coordinates. Orient to even parity, so distinct representatives have Hamming distance two. Translating one to the zero binary word, the others are pairwise intersecting two-subsets of a five-set, of which there are at most four: any intersecting family of edges is a star or lies in a triangle. Thus each block has capacity five. Since the span has dimension at most 20, choose a nonzero normal $h$ in the 21-dimensional exterior coordinate space. Write $J=\supp h$. Any block meeting $J$ has capacity at most four, since the Gram matrix of five compatible vectors has diagonal 5 and off-diagonal $\pm1$ and is strictly diagonally dominant. Any block meeting $J$ in exactly one point has capacity zero. If $|J|\ge5$, five points of $J$ meet at least $25-\binom52=15$ blocks, giving $M\le105-15=90$. If $1\le |J|=k\le4$, each point of $J$ lies in at least $5-(k-1)$ blocks containing no other point of $J$. These blocks are distinct for different points, so at least $k(6-k)\ge5$ blocks have capacity zero; $M\le80$.

\subsubsection{Eight entries 1 and one entry 2}
Let $O$ be the octad of ones, $a$ the position of the two, and $Z$ the remaining 15 coordinates. At most eight candidate axes have nonzero centered coordinate $a$. Indeed, both $z$ and $b-z$ have square eight, so $|z_a|\le2$ and $|2-z_a|\le2$, whence $z_a\in\{0,1,2\}$. If $z_a\ne1$, the antipodal operation $z\mapsto b-z$ permits choosing $z_a=0$. With $z_a=0$, the two-2 type pairs at most four, and the only octad giving pairing six is $O$ with exactly one negative sign. By \eqref{d20:support}, all eight exceptional axes must be discarded. The remaining axes have the form
\[
 v=s_T+q_B,\qquad s_T=\mathbf1_T-\tfrac12\mathbf1_O,\quad |T|=4,\qquad q_B=\sum_{i\in B}\varepsilon_i e_i,
\]
where $T\cup\{a\}\cup B$ is an octad, $B\subseteq Z$ has size three, and the signs on $T$ and $a$ are positive. Adding $O$ to the octad in the binary code replaces $T$ by $O\setminus T$. For fixed $B$ at most one complementary partition of $O$ occurs: otherwise two distinct octads with the same $B$ can be oriented by complements to share at least two points of $O$, hence at least six points overall. Distinct triples cannot share two points, since after complementing one $T$ the two octads would share at least $2+1+2=5$ points. The exterior triples form $S(2,3,15)$: for a fixed $b'\in Z$, the unique octad on $\{a,b'\}\cup U$, $U\in\binom O3$, meets $O$ in four points and is counted four times; the complementary octad with the same $B$ accounts for four more. There are $56/8=7$ triples through $b'$. Thus there are 35 triples in all, each pair used at most once and hence exactly once.

Orient each used triple by one of its two complementary $T$ choices. Its sign words have mutual Hamming distance two or three. A complementary pair precludes a third word; otherwise all words have the same parity. Hence each triple has capacity four, and three selected signs form three vertices of a regular tetrahedron. In particular, the span of their differences has a one-dimensional orthogonal complement in $\R^3$, generated by a vector with every coordinate nonzero; four selected signs affinely span $\R^3$.

Suppose 91 axes are selected, and let $q\ge23$ be the number of used triples. Their internal vectors $s_T$ span the whole seven-dimensional sum-zero space on $O$. To see this, let $G$ be their Gram matrix and let $d_i$ be the degrees of the used triple system. Its diagonal is 2. Distinct triples meeting in one point have $|T\cap T'|=2$ (use the octad intersection bound on both complementary choices), so their internal pairing is zero. For disjoint triples, the octad intersections have sizes 2 or 4, giving $|T\cap T'|=1$ or 3, so the absolute internal pairing is one. Therefore
\[
 \tr G=2q,\qquad \tr G^2=q^2+6q-\sum_i d_i^2\le\tfrac25q^2+6q.
\]
Since $\rk G\ge4/(2/5+6/q)>6$, the claimed span is seven-dimensional.

If at least two exterior points are inactive, at most $35-7-7+1=22$ triples remain, giving at most 88 axes. If all 15 points are active, choose independent normals $h,h'$ in the ambient $7+15=22$ dimensional space of internal and exterior coordinates. The exterior columns $H_i=(h_i,h'_i)$ span $\R^2$: otherwise a nonzero linear combination of the normals would be purely internal, contradicting the established internal span. Divide the nonzero columns into classes up to proportionality, of sizes $s_1,\ldots,s_k$, and let $z$ be the number of zero columns. Here $k\ge2$. A triple with a nonzero column has capacity at most three; if it contains three selected signs, the tetrahedron difference property implies that all three columns are nonzero and proportional. Fix a nonzero column. Its point lies on at least 15 selected axes by \eqref{d20:support}; of its seven incident triples, at least one must contain three axes. Thus every nonzero class has size at least three, and $z\le9$.

For any set of $s$ points in the $S(2,3,15)$ design, the number of triples entirely contained in that set is at most
\[
b(s)=\min\left(\left\lfloor\tfrac{\binom s2}{3}\right\rfloor,
\left\lfloor\tfrac{s\lfloor(s-1)/2\rfloor}{3}\right\rfloor,
35-7s+\tbinom s2\right).
\]
The first two terms count pairs and internal vertex degrees respectively. For the third, if $b_j$ counts triples meeting the set in $j$ points, then $\sum b_j=35$, $\sum jb_j=7s$, and $b_2+3b_3=\binom s2$. Consequently $b_0=35-7s+\binom s2-b_3\ge0$. A triple has capacity at most two except that a triple wholly in a nonzero class has capacity three, and an all-zero triple has capacity four. It follows that
\[
M\le70+\sum_{j=1}^{k}b(s_j)+2b(z).
\]
The following small integer-partition check gives the maximum of the excess term, with $k\ge2$, $s_j\ge3$ and $\sum s_j=15-z$:
\[
\begin{array}{c|rrrrrrrrrr}
z&0&1&2&3&4&5&6&7&8&9\\\hline
\text{maximum}&18&14&11&12&10&14&13&18&16&18
\end{array}
\]
These entries follow by enumerating nondecreasing partitions of integers at most 15 into parts at least three, or by the recurrence $F_k(n)=\max_{3\le s\le n-3(k-1)}(b(s)+F_{k-1}(n-s))$, with $F_1(n)=b(n)$ for $n\ge3$. The accompanying Golay checker verifies all ten entries directly. Thus $M\le88$.

If exactly one exterior point is inactive, 28 triples survive. In the ambient space of dimension $7+14=21$ choose one nonzero normal. Some exterior coefficient is nonzero by the internal span. Its point lies in six surviving triples and on at least 15 selected axes; at least three of those triples therefore contain three axes. The difference property forces at least seven nonzero exterior coefficients. Write $z$ for the number of zero coefficients, so $z\le7$. If $z\le6$, there are at most $b(z)\le4$ all-zero triples; the bound is $M\le28\cdot3+4=88$. If $z=7$, there are seven nonzero coefficients as well. The numbers of all-zero and all-nonzero triples are each at most seven. All mixed triples have capacity two, so $M\le28\cdot2+7\cdot2+7=77$. This completes the second marking case.

\subsubsection{Twelve entries 1}
Let $D$ be the dodecad. A candidate $z$ must be an octad meeting $D$ in exactly six points with positive signs there: intersections of an octad with $D$ cannot exceed six by the code weight set. Every axis is thus $s_T+\varepsilon_i e_i+\varepsilon_j e_j$ with $|T|=6$, $T\subset D$, $s_T=\mathbf1_T-\mathbf1_D/2$ and $i,j\in D^c$. A fixed exterior pair has at most one complementary partition $\{T,D\setminus T\}$, since otherwise complementing a part makes two octads intersect in at least $3+2=5$ points. Orient by one part. The squared internal norm is three and the sign words of length two must differ in both positions for compatibility, so at most two axes occur per pair. Suppose there are 91 selected axes. If the number $a$ of active exterior points is at most ten, the bound $2\binom a2\le90$ already suffices. Hence $a=11$ or 12. Let $q\ge46$ be the number of used exterior pairs and let $d_i$ be the degrees of the graph of these pairs. Their associated internal vectors $s_T$ have square three. For two distinct pairs meeting in one exterior point, the octad intersection constraints, applied also after replacing $T$ by $D\setminus T$, force the internal intersection to have size three. For disjoint exterior pairs, the two complementary choices force internal intersection size two or four. Thus the internal Gram matrix $H$ has off-diagonal entries zero for pairs sharing a vertex and absolute value one for disjoint pairs. Consequently
\[
\tr H=3q,\qquad \tr H^2=q^2+10q-\sum_i d_i^2
\le (1-4/a)q^2+10q.
\]
The rank inequality implies
\[
\rk H\ge \frac{9}{1-4/a+10/q}>10
\quad (a\le12,\ q\ge46).
\]
The internal span is therefore the entire eleven-dimensional sum-zero space on $D$.

The ambient space of the selected vectors now has dimension $11+a$, and their span has dimension 20. Choose a basis of the normal space, of dimension $a-9$, and write $H_i'$ for the column of normal coefficients on exterior point $i$. These columns span $\R^{a-9}$: a nonzero linear combination of the normals with zero exterior part would contradict the internal span. Form the graph $G_2$ of doubled exterior pairs (those with two selected axes). On such a pair $\{i,j\}$ the two sign words are opposite; the difference of the corresponding vectors is $2(\varepsilon_i e_i+\varepsilon_j e_j)$. Hence $H_i'=\pm H_j'$. Each connected component of $G_2$ contributes at most one dimension to the span of these columns. On the other hand, every active point occurs on at least 15 axes by \eqref{d20:support}. At most $a-1$ used pairs contain it, and a doubled pair contributes one additional axis. Thus $G_2$ has minimum degree at least $16-a$. Each component has size at least $17-a$, so there are at most $\lfloor a/(17-a)\rfloor$ components. This is at most one for $a=11$ and at most two for $a=12$, strictly below the required column rank $a-9$. This contradiction completes the binary Golay proof.

\subsection{The Leech lattice: a Golay-code obstruction}\label{d20:sec:leech}
We use the standard coordinate model of the Leech lattice, with integral coordinates divided by $\sqrt8$, and the classical fact that its norm-six vectors form one orbit under its isometry group \cite[p.~493]{Allcock}. In particular we may take
\[
w=W/\sqrt8,\qquad W=(5,1,\ldots,1).
\]
In this model the minimal vectors with even integral coordinates have either two entries of absolute value four, or eight entries of absolute value two on an octad of $G_{24}$. All coordinates of a lattice vector have the same parity. For a point $x$ in $\mathcal F_P(w)$, exactly one of $x,w-x$ has even coordinates. The condition that its integral coordinate vector have scalar product 24 with $W$ forces the former type to have positive fours at coordinate zero and one other coordinate, and the latter type to have positive twos on an octad through zero. Both types with these signs belong to the Leech lattice in its standard model. Thus, writing $G$ for the punctured Golay code at coordinate zero and $\mathcal H$ for its heptads, the section consists of the axes with representatives
\begin{align*}
p_j&=\tfrac14(3,-1,\ldots,-1)+2e_j && (1\le j\le23),\\
q_H&=-\tfrac14\mathbf1_{24}+\mathbf1_H && (H\in\mathcal H).
\end{align*}
Here $e_j$ is the $j$th exterior coordinate unit vector and $\mathbf1_H$ is zero at coordinate zero. These are the centered representatives $\sqrt2(x-w/2)$, of square five. The heptads form $S(4,7,23)$ by the octad design verified in Section~\ref{d20:sec:golay}. Distinct heptads intersect in one or three points, since their extensions are distinct octads through zero. Directly,
\[
p_j\cdot p_k=1\ (j\ne k),\qquad p_j\cdot q_H=2\mathbf1_{j\in H}-1,\qquad q_H\cdot q_{H'}=|H\cap H'|-2.
\]
In particular these 276 axes form an equiangular configuration in $W^\perp$ \cite[p.~134]{Atlas}. We write their representatives collectively as $v_i$.

\subsubsection{The evaluation lattice and its moments}
The matrices $Q_i=v_iv_i^{\mathsf T}$ form a basis of $\operatorname{Sym}(W^\perp)$: their Frobenius Gram matrix is $24I+J$, and there are $276=23\cdot24/2$ of them. Testing against this basis gives $\sum_iQ_i=60I$ and $\sum_i\tr(AQ_i)Q_i=24A+12\tr(A)I$ for symmetric $A$. Consequently
\begin{align}
 \sum_i\ip{y}{v_i}^2&=60y^2,\label{d20:second}\\
 \sum_i\ip{y}{v_i}^2\ip{z}{v_i}^2&=12y^2z^2+24\ip{y}{z}^2.\label{d20:fourth}
\end{align}
Let $L=\sum_i\Z v_i$ in $W^\perp$. It is a full integral lattice (its generators have integral pairings and span the space). Its dual has the minimum-norm estimate
\begin{equation}\label{d20:dualmin}
 y\in L^*\setminus\{0\}\quad\Longrightarrow\quad y^2\ge\tfrac53.
\end{equation}
Indeed the integer evaluations $t=(y,v_i)$ satisfy $t^4\ge t^2$, and \eqref{d20:second}--\eqref{d20:fourth} give $36(y^2)^2\ge60y^2$. Equality in \eqref{d20:dualmin} is equivalent to all evaluations belonging to $\{0,\pm1\}$.

Define the binary evaluation code $C$ by the image of
\[
\phi:L^*\longrightarrow\F_2^{276},\qquad y\longmapsto ((y,v_i)\bmod2)_i.
\]
The kernel is $2L^*$, by the definition of $L$, so $C$ has dimension 23. For the binary $[276,23,100]$ code invariant under $\mathrm{Co}_3$ and its weight distribution, see \cite{HPPT}. The next lemma establishes all the weight information required below by incidence counting, without enumerating the $2^{23}$ words of $C$.

\begin{lemma}\label{d20:codeweights}
Every nonzero word of $C$ has weight in
\[
\{100,112,128,132,144,148,164,176,276\}.
\]
Furthermore every word of weight 100 is the image under $\phi$ of a vector in $L^*$ of square $5/3$.
\end{lemma}
\begin{proof}
For $y\in W^\perp$ put $t_j=(y,p_j)$ and $S=\sum_{j=1}^{23}t_j$. The orthogonality to $W$ gives $t_j=2(y_0+y_j)$, and therefore
\begin{equation}\label{d20:tcoords}
 y_0=S/36,\qquad y_j=t_j/2-S/36,\qquad
 y^2=\frac14\sum_jt_j^2-\frac{S^2}{108}.
\end{equation}
Moreover
\begin{equation}\label{d20:hepeval}
 (y,q_H)=\frac{3\sum_{j\in H}t_j-S}{6}.
\end{equation}
Thus $y\in L^*$ precisely when $t_j\in\Z$, $3\mid S$, and $\sum_{j\in H}t_j\equiv S\pmod2$ for every heptad.

The punctured Golay code $G$ has dimension 12, weights in $\{0,7,8,11,12,15,16,23\}$, and contains $\mathbf1$. Its parity extension $G_{24}$ is self-dual because it is doubly even of dimension 12. Hence the even subcode of $G$ equals $G^\perp$ (orthogonality follows from the extension and the dimensions are both 11). The twelve shifts of the generator polynomial used in Section~\ref{d20:sec:golay} are heptads and span $G$. The map $g\mapsto g+\wt(g)\mathbf1$ over $\F_2$ takes $G$ onto its even subcode, with kernel $\langle\mathbf1\rangle$. It follows that the complements of the heptads span $G^\perp$. The above integrality conditions on $y$ are therefore equivalent to $t\in\Z^{23}$, $3\mid S$, and $g=t\bmod2\in G$.

Write $t\equiv g+2b\pmod4$, viewing $g,b$ as zero-one vectors when using integer arithmetic. The binary evaluations on the points are $g$, and on the heptads are
\begin{equation}\label{d20:wordformula}
 \phi(y)_H=\frac{3|H\cap g|-|g|}{2}+|H\cap b|+|b|\pmod2.
\end{equation}
To obtain this, divide the numerator in \eqref{d20:hepeval} by two, and use that division by three does not change its parity. In \eqref{d20:wordformula}, replacing $b$ by a binary translate by an element of $G$ has no effect, by the heptad-complement orthogonality. The radius-three balls of $G$ are disjoint by its minimum distance seven and cover $\F_2^{23}$ since $2^{12}\sum_{i=0}^3\binom{23}{i}=2^{23}$. Thus we may choose $|b|\le3$. Conversely every pair $(g,b)$ so obtained arises from an element of $L^*$: choose $t$ with the prescribed residues modulo four and adjust one entry by a multiple of four to make $S$ a multiple of three. This proves the asserted parameterization of the words. Replacing $g$ by its complement with $b$ unchanged complements the evaluation word: the integer first term in \eqref{d20:wordformula} changes from $A$ to $-1-A$. It suffices to treat $k=|g|\in\{0,7,8,11\}$.

Here are the incidence computations, including the rules for reproducing the table. Set $\lambda_j=\binom{23-j}{4-j}/\binom{7-j}{4-j}$, so $(\lambda_0,\ldots,\lambda_4)=(253,77,21,5,1)$. For $J\subseteq b$ let $u=|J|$, $v=|J\cap g|$, and write $n_i(J)$ for the number of heptads containing $J$ and meeting $g$ in $i$ points. The allowed $i$ are, respectively,
\[
\begin{array}{c|cccc} k&0&7&8&11\\\hline i&0&1,3,7&0,2,4&1,3,5\end{array}
\]
by the possible weights of $H+g$ and $0\le i\le\min(7,k)$. For $0\le j\le4-u$, double counting pairs consisting of a heptad and a $j$-subset of $(H\cap g)\setminus J$ gives
\begin{equation}\label{d20:incidences}
 \sum_i\binom{i-v}{j}n_i(J)=\binom{k-v}{j}\lambda_{u+j}.
\end{equation}
Terms with $i<v$ or $7-i<u-v$ are zero. Define $A_g(J)=\sum_i(-1)^{(3i-k)/2}n_i(J)$. For $k=0$ it equals $\lambda_u$. For $k=7$ the unique heptad meeting $g$ in seven points is $g$ itself, so $n_7(J)=\mathbf1_{u=v}$. The first moment gives
\[
 A_g(J)=2\lambda_u-\big((k-v)\lambda_{u+1}+v\lambda_u\big)+4\mathbf1_{u=v}.
\]
For $k=8,11$, write the three possible intersections as $\alpha,\alpha+2,\alpha+4$ (with $\alpha=0,1$ respectively). The middle value alone has negative sign. If $u\le2$, the first three equations in \eqref{d20:incidences} determine the three counts, and equivalently
\begin{align*}
 M_1&=(k-v)\lambda_{u+1}+v\lambda_u,\\
 M_2&=2\binom{k-v}{2}\lambda_{u+2}+(2v+1)(k-v)\lambda_{u+1}+v^2\lambda_u,\\
 A_g(J)&=\lambda_u-\tfrac12\big(-M_2+(2\alpha+4)M_1-\alpha(\alpha+4)\lambda_u\big).
\end{align*}
This last formula uses the identity $\mathbf1_{i=\alpha+2}=(i-\alpha)(\alpha+4-i)/4$ on the three allowed values. If $u=3$, there are only five heptads containing $J$; retain all nonnegative integer triples $(n_\alpha,n_{\alpha+2},n_{\alpha+4})$ of sum five satisfying $\sum i n_i=M_1$ and the zero restrictions above. In that case the possible values are $A_g(J)=5-2n_{\alpha+2}$. These rules may overcount the possibilities, which is harmless for an upper restriction on weights.

Finally expand $(-1)^{|H\cap b|}=\prod_{j\in b}(1-2\mathbf1_{j\in H})$ in \eqref{d20:wordformula} to obtain
\begin{equation}\label{d20:weighttableformula}
 \wt(\phi(y))=k+\tfrac12\left(253-(-1)^t\sum_{J\subseteq b}(-2)^{|J|}A_g(J)\right),\qquad t=|b|.
\end{equation}
In this sum there are $\binom{r}{v}\binom{t-r}{u-v}$ subsets $J$ of type $(u,v)$, where $r=|g\cap b|$. Substitution gives the allowed weights in Table~\ref{d20:table:incidences} (a dash denotes an impossible value of $r$). All arithmetic uses only the five incidence numbers and $t\le3$.
\begin{table}[htbp]
\centering\small
\caption{Allowed evaluation-code weights grouped by the incidence parameters.}\label{d20:table:incidences}
\begin{tabular}{ccrrrr}\toprule
$k$&$t$&$r=0$&$r=1$&$r=2$&$r=3$\\\midrule
0&0&0&--&--&--\\
0&1&176&--&--&--\\
0&2&112&--&--&--\\
0&3&128&--&--&--\\\midrule
7&0&148&--&--&--\\
7&1&112&164&--&--\\
7&2&148&128&100&--\\
7&3&128&132&144&148\\\midrule
8&0&176&--&--&--\\
8&1&128&100&--&--\\
8&2&128&148&144&--\\
8&3&112,128,144&132&128&148\\\midrule
11&0&176&--&--&--\\
11&1&132&112&--&--\\
11&2&128&148&144&--\\
11&3&148&112,128,144&132&144\\\bottomrule
\end{tabular}
\end{table}
For example, with $k=8,u=3,v=0$, the counts satisfy $n_0+n_2+n_4=5$ and $2n_2+4n_4=8$, so $n_2\in\{0,2,4\}$ and $A_g(J)\in\{-3,1,5\}$; the first two moments for the proper subsets of $J$ then give the indicated three weights. The accompanying short incidence-table checker repeats these elementary substitutions, not an enumeration of evaluation words.

Taking the table together with complement symmetry proves the weight restriction. It also shows that weight 100 occurs only for
\[
 (k,t)=(7,2),(8,1),(12,0),(15,0),(23,1),\qquad b\subseteq g.
\]
For any such word take the integer vector $t_j=\mathbf1_g(j)-2\mathbf1_b(j)$. Its sum $S=k-2t$ belongs to $\{3,6,12,15,21\}$, so the vector $y$ from \eqref{d20:tcoords} lies in $L^*$. Its residues modulo four are the prescribed ones, and its square is
\[
 \frac{k}{4}-\frac{(k-2t)^2}{108}=\frac53
\]
in each of the five cases. This proves the lifting assertion.
\end{proof}

\subsubsection{Three-dimensional subcodes}
If $a,b$ are minimum words with lifts $y_a,y_b$ of square $5/3$, their integral evaluations are in $\{0,\pm1\}$ by \eqref{d20:dualmin}. The support intersection count from \eqref{d20:fourth} implies
\begin{equation}\label{d20:inner}
36\ip{y_a}{y_b}^2=100-\tfrac34\wt(a+b).
\end{equation}
For $a=b$ this gives uniqueness of minimum lifts up to sign. If $a,b,a+b$ are minimum words, choose signs so $\ip{y_a}{y_b}=-5/6$. Then $y_a+y_b$ has square $5/3$ and residue $a+b$, so it is a minimum lift of that word.

\begin{proposition}\label{d20:d3}
Every three-dimensional subcode of $C$ has support size at least 186.
\end{proposition}
\begin{proof}
If a subcode $D$ has support size at most 185, the sum of its seven nonzero weights is $4|\supp D|\le740$. Let $k$ be the number of weight-100 words. Then $k\ge4$. If there is a Fano line consisting of minimum words $a,b,c=a+b$, choose another minimum word $d$ outside it, and take lifts with $y_c=y_a+y_b$. The three positive magnitudes $|6\ip{y_d}{y_a}|$, $|6\ip{y_d}{y_b}|$, $|6\ip{y_d}{y_c}|$ must have one equal to the sum of the other two. By the weight sum restriction the remaining three words have weights $(112,112,112)$ if $k=4$, $(100,112,112)$ or $(100,112,128)$ if $k=5$, and at least two 100s with third weight at most 132 if $k\ge6$. Formula \eqref{d20:inner} thus gives $(4,4,4)$, one 5 and two positive even integers, or at least two 5s with the third between 1 and 5. None satisfies the required addition.

Hence the minimum-word set contains no Fano line. Translating it by any of its elements gives a disjoint set, so it has size at most four. Thus $k=4$. Three distinct minimum words are independent, and the fourth must be their sum (it cannot be a pair sum). Hence $a_1+a_2+a_3+a_4=0$. The three nonminimum weights must all equal 112, so each pair of distinct minimum lifts has inner product $\pm2/3$. Three times their Gram matrix is integral and congruent to $I_4$ modulo 2; consequently the four lifts are linearly independent. Every signed sum $\sum_i\varepsilon_i y_{a_i}$ lies in $2L^*\setminus\{0\}$ and has square at least $20/3$. Choose the first two signs with cross term $-4/3$ and average over the other two independent uniform signs. The average square is $20/3-4/3=16/3$, a contradiction.
\end{proof}

The following bound works in any rank. For a linear code \(C_p\) over \(\F_p\), let \(d_j(C_p)\) be the least size of the union of the supports of a \(j\)-dimensional subcode~\cite{Wei}.

\begin{lemma}[Evaluation codes and sections]\label{lem:section-weights}
Let \(v_1,\ldots,v_n\) generate a lattice \(L\) of rank \(r\), and let \(C_p\) be the image of \(L^*\) under evaluation on the \(v_i\) modulo a prime \(p\). For \(1\le j\le r\), a real subspace of dimension at most \(r-j\) contains at most \(n-d_j(C_p)\) of the indexed vectors.
\end{lemma}
\begin{proof}
Let \(I\) consist of the indices in the subspace and put \(U=\operatorname{span}_{\R}\{v_i:i\in I\}\). Since \(U\) is rational relative to \(L\), its annihilator \(K=L^*\cap U^\perp\) has rank \(r-\dim U\ge j\). The inclusion induces an injection \(K/pK\hookrightarrow L^*/pL^*\), since \(py\in K\) for \(y\in L^*\) implies \(y\in K\). The evaluation map identifies \(L^*/pL^*\) with \(C_p\) because the \(v_i\) generate \(L\). The image of \(K/pK\) vanishes on \(I\), so a \(j\)-dimensional subcode has support of size at most \(n-|I|\).
\end{proof}

For the 276 Leech axes and \(p=2\), Proposition~\ref{d20:d3} and Lemma~\ref{lem:section-weights} give the upper bound \(276-186=90\) in dimension twenty.

\subsection{A sharp construction and completion of the proof}\label{d20:sec:lower}
For completeness we give a direct incidence count for the classical lower bound $N(20)\ge90$ (see Taylor \cite[p.~123]{Taylor} and Lin--Yu \cite[Example~2 and Table~2]{LY}). Use the point and heptad representatives in Section~\ref{d20:sec:leech}. Choose a weight-eight word $O\in G$ and distinct points $a,b\in O$. Such a word exists, since the octads not containing zero puncture to weight-eight words. Take the three normals $y_a,y_b,y_O$ given by \eqref{d20:tcoords} with respective integer evaluation columns
\[
t^{(a)}=\mathbf1-2e_a,\qquad t^{(b)}=\mathbf1-2e_b,\qquad t^{(O)}=-\mathbf1_{O^c}.
\]
These three columns are independent: a linear relation evaluated on $O\setminus\{a,b\}$ and then at $a$ forces the coefficients of the first two to vanish, and $O^c$ forces the last coefficient to vanish. The invertible map \eqref{d20:tcoords} gives independent normals in $W^\perp$. By \eqref{d20:hepeval},
\[
 (y_a,q_H)=-\mathbf1_{a\in H},\qquad (y_b,q_H)=-\mathbf1_{b\in H},\qquad (y_O,q_H)=|H\cap O|/2-1.
\]
Order the normals as \(y_a,y_O,y_b\) and let \(V_{23-j}\) be the common kernel of the first \(j\) in \(W^\perp\). For \(j\ge1\), every \(p_i\) has a nonzero evaluation on \(y_a\). The heptads in the first kernel avoid \(a\), giving \(253-77=176\) axes. The second kernel imposes \(|H\cap O|=2\), and the third also excludes \(b\). The possible weights of \(H+O\) restrict \(|H\cap O|\) to \(0,2,4\). Using \eqref{d20:incidences} with \(g=O\) and \(J=\varnothing,\{a\},\{a,b\}\) gives \(168,42,6\) heptads meeting \(O\) in two points and containing \(J\), respectively. The count for \(J=\{b\}\) is also \(42\). Hence the numbers of axes in the nested sections are
\[
\begin{array}{c|rrrr}
\dim V & 23 & 22 & 21 & 20\\
\#\text{ axes in }V & 276 & 176 & 126 & 90
\end{array}
\]
with \(126=168-42\) and \(90=168-2\cdot42+6\). The rank inequality \(\rk G\ge25m/(m+24)\) for the Gram matrix of \(m\) axes shows that the axes in each section span it.

In each case the section is largest among all sections of the corresponding dimension of these 276 axes. Indeed, Lemma~\ref{d20:codeweights} gives \(d_1(C)\ge100\). The three nonzero words of a binary two-dimensional subcode each have weight at least \(100\), and the sum of their weights is twice the size of its support. Thus \(d_2(C)\ge150\), and Proposition~\ref{d20:d3} gives \(d_3(C)\ge186\). Lemma~\ref{lem:section-weights} and the constructed sections give the reverse bounds, so the first three generalized Hamming weights are \(100,150,186\). In particular the twenty-dimensional section supplies the required lower bound.

\begin{proof}[Proof of Theorem~\ref{d20:main}]
The preceding construction gives the lower bound. For the upper bound, a 91-line system in $\R^{20}$ has common cosine $1/5$ by Lemma~\ref{d20:angle}, and Lemma~\ref{lem:paired-sections} realizes it in a norm-six section of an even unimodular lattice of rank 24. The Niemeier classification leaves exactly the twenty-two nonbinary rooted types in Section~\ref{d20:sec:rooted}, the binary rooted type in Section~\ref{d20:sec:golay}, and the Leech lattice in Section~\ref{d20:sec:leech}. Sections~\ref{d20:sec:rooted}--\ref{d20:sec:leech} exclude all these possibilities, using the finite verifications specified there for the rooted types. Hence no 91-line system exists. The lower and upper bounds coincide.
\end{proof}

\subsection{Verification of the finite computations}\label{d20:sec:verification}
The finite checks for the nonbinary rooted lattices are specified by the local Cartan-coordinate and glue-code rules in Section~\ref{d20:sec:rooted}. The repository supplement contains the root-data, energy-bound and glue-model routines, the scalar and small-model verifiers, and the two Golay and incidence checkers. The file \texttt{README.md} lists their individual roles. From the repository root run
\begin{verbatim}
cd dimension-20/reproduction
unset PYTHONOPTIMIZE
python3 -B audit/referee/replay.py
\end{verbatim}
The replay checks the SHA-256 digests, reproduces the tables, and checks the elementary Golay inputs with Python 3.10 or later, using only the standard library and integer or rational arithmetic. Assertions must be enabled. The 1885 scalar profiles in the twenty general types and the 37 local alternatives for the two small types are verified without constructing compatibility graphs. The Leech argument is analytical: the incidence-table program is a check of the short table in Lemma~\ref{d20:codeweights}, not a search of the evaluation code. No numerical tolerances or unverified clique-search outputs enter the proof.

\section{Sharp fixed-angle bounds}\label{fa:section}
The 276 lines of Section~\ref{d20:sec:leech} also give fixed-angle results in a different range of dimensions. We use their Seidel spectrum, and the spectrum of the 28-line system in \(\R^7\), to construct intervals on which \(N_{1/q}(d)-d\) is constant. These arguments are independent of the exclusions in dimensions eighteen, nineteen, and twenty. For each fixed odd reciprocal cosine, Jiang, Tidor, Yao, Zhang and Zhao~\cite[Theorem~1.2]{JTYZZ} determine the exact value in all sufficiently large dimensions. Our intervals concern dimensions that grow linearly with the varying odd denominator. Here the angle is prescribed, and \(N_\alpha(d)\) allows the lines to span any subspace of \(\R^d\).

A real Seidel matrix \(S\) is symmetric with zero diagonal and off-diagonal entries in \(\{-1,1\}\). If \(q>1\) and \(qI+S\) is positive semidefinite, then \(q^{-1}(qI+S)\) is the Gram matrix of distinct lines of common cosine \(1/q\). Their span has dimension \(\rk(qI+S)\). We write \(J_s\) for the all-ones matrix of order \(s\).

\subsection{Intervals for odd reciprocal cosines}\label{fa:odd}
Use the following parameters, indexed by the excess \(k\) of the line count over the dimension.
\[
\begin{array}{c|rrrr}
k&M&b&a&L_k\\ \hline
7&28&10&3&72/(7\sqrt{10}+2)\\
23&276&56&5&600/(23\sqrt{26}+2)
\end{array}
\]
The parameters \(M\) and \(1/a\) are the size and cosine of the seed system in \(\R^k\). The constant \(b\) satisfies \(b-1=aM/k-a\).

\begin{samepage}
\begin{theorem}\label{fa:intervals}
Let \(q>1\) be an odd integer and use either row of the table. Set
\[
Q=q+1,\qquad t=2\left\lceil\frac{Q}{2b}\right\rceil,\qquad
r=\frac{bt-Q}{2},
\]
and define
\[
d_-=Mt-k,\qquad d_+=\lfloor L_kQ\rfloor-k-1.
\]
If \(r<t\) and \(d_-\le d_+\), then
\[
N_{1/q}(d)=d+k\qquad(d_-\le d\le d_+).
\]
Both conditions hold for all odd \(q\ge319\) when \(k=7\), and for all odd \(q\ge5599\) when \(k=23\). For each row, the number of integers in the interval is \((L_k-M/b)Q+O(1)\).
\end{theorem}
\end{samepage}

\begin{proof}
\emph{Upper bound.}
For \(h\ge1\), write
\[
\beta_h=\frac{(h-1)\sqrt{h+2}+2}{h(h+1)}.
\]
The real coherence bound of Bukh and Cox~\cite[Theorem~2(a)]{BukhCox} says that \(n=d+h\) unit vectors in \(\R^d\) have maximum absolute inner product at least \(1/(\beta_h n-1)\) in the range used here. The same inequality holds for a configuration of smaller span by approximation with \(n\) unit vectors spanning \(\R^d\). In particular, an equiangular set of \(d+h\) lines of cosine \(1/q\) must satisfy \(d+h\ge Q/\beta_h\). With \(h=k+1\) we have \(L_k=1/\beta_h\). If a larger set existed for \(d\le d_+\), we could retain just \(d+k+1\) lines, with
\[
d+k+1\le\lfloor L_kQ\rfloor<L_kQ.
\]
Here the strict inequality follows from the irrationality of \(L_kQ\). This proves \(N_{1/q}(d)\le d+k\).

\emph{A matrix of rank \(Mt-k\).}
For \(k=7\), take the 28 lines obtained from the centered incidence vectors of the pairs of an eight-element set. Their common cosine is \(1/3\), since the vectors \(\mathbf 1_{\{i,j\}}-\tfrac14\mathbf 1_8\) have squared norm \(3/2\) and pairings \(\pm1/2\) for distinct pairs. For \(k=23\), use the Leech system in Section~\ref{d20:sec:leech}. Both systems attain equality in the Gram trace bound \(\rk G\ge(\tr G)^2/\tr(G^2)\), so their nonzero Gram eigenvalues are equal. If \(S\) is the seed Seidel matrix, its eigenvalues are consequently \(b-1\) with multiplicity \(k\) and \(-a\) with multiplicity \(M-k\).

The ceiling in the definition of \(t\) gives \(0\le r<b\), and \(r\) is an integer since \(bt\) and \(Q\) are even. As \(r<t\) and \(t\) is even, we can choose the \(r\)-regular circulant graph on \(\Z/t\Z\) with differences \(\pm1,\ldots,\pm\lfloor r/2\rfloor\) and, if \(r\) is odd, \(t/2\). Let \(A\) be its adjacency matrix. Define
\begin{equation}\label{fa:block-gram}
G=Q I_{Mt}-(S+I_M)\otimes J_t+2I_M\otimes A.
\end{equation}
The diagonal entries are \(q\). The off-diagonal entries are those of \(-S\) between blocks, and \(-1+2A_{ij}\) within blocks. When \(r=0\), the Seidel matrix in \eqref{fa:block-gram} is the tensor construction of Greaves, Koolen, Munemasa and Sz\H{o}ll\H{o}si~\cite[Proposition~5.11, applied to \(-S\)]{GKMS}. The corresponding coherence values for the two seeds are also in Bukh and Cox~\cite[Corollary~17]{BukhCox}.

On vectors that are constant in every block, \(G\) acts as
\[
(Q+2r)I_M-t(S+I_M)=t((b-1)I_M-S).
\]
This is positive semidefinite with a \(k\)-dimensional kernel. On vectors with sum zero in every block, \(G\) acts as \(QI_{Mt}+2I_M\otimes A\). The latter has all eigenvalues at least
\[
Q-2r=bt-4r>(b-4)t>0,
\]
using the eigenvalue bound \(\|A\|\le r\). Hence \(G\) has rank \(Mt-k\).

\emph{Extending the construction to consecutive dimensions.}
Let \(z\in\{-1,1\}^{Mt}\) be the vector whose signs alternate along \(\Z/t\Z\) in each block. Alternation is a character of the cyclic group, so \(z\) is a \(G\)-eigenvector with eigenvalue \(\lambda\ge Q-2r>0\). Put \(\rho=z^{\mathsf T}z/\lambda=Mt/\lambda<M/(b-4)\). For a positive integer \(u\), consider
\[
\widehat G=
\begin{pmatrix}
 G&z\mathbf 1_u^{\mathsf T}\\
 \mathbf 1_u z^{\mathsf T}&(q-1)I_u+J_u
\end{pmatrix}.
\]
Since \(G(z/\lambda)=z\), a block congruence gives the direct sum of \(G\) and \((q-1)I_u+(1-\rho)J_u\). The second summand is positive definite when \(q-1>(\rho-1)u\). In that case \(\widehat G\) realizes \(Mt+u\) lines of cosine \(1/q\) and rank \(Mt-k+u\).

For \(d=Mt-k+u\le d_+\) we have \(u<L_kQ-Mt\le(bL_k-M)t\). When \(k=7\), the inequalities \(bL_k-M<2\) and \(\rho-1<4\) give \((\rho-1)u<8t\). When \(k=23\), \(bL_k-M<6\) and \(\rho-1<5\) give \((\rho-1)u<30t\). These estimates follow by using the displayed values of \(L_k\) and the bound on \(\rho\). In both cases \(q-1=bt-2r-2\ge(b-2)t\). This gives the required positive definiteness for \(u>0\), and for \(u=0\) the matrix \(G\) suffices.

Finally \(t<Q/b+2\) gives \(L_kQ-Mt>(L_k-M/b)Q-2M\). For \(k=7\) we have \(L_k-M/b>0.18\), so for \(Q\ge320\) the right side exceeds \(1\), and \(t\ge32>b-1\ge r\). For \(k=23\) we have \(L_k-M/b>0.10\), so for \(Q\ge5600\) the right side again exceeds \(1\), and \(t\ge100>b-1\ge r\). The inequalities for \(L_k\) in this proof follow, for example, by squaring the rational bounds \(3.162<\sqrt{10}<3.163\) and \(5.099<\sqrt{26}<5.100\). Thus \(d_-\le d_+\) at and above the stated thresholds. The interval length formula follows from \(t=Q/b+O(1)\).
\end{proof}

The slopes are approximately \(0.18310\) for \(k=7\) and \(0.10172\) for \(k=23\). The conditions in Theorem~\ref{fa:intervals} also hold for some smaller odd denominators. Direct substitution gives the following examples.
\begin{center}
\begin{tabular}{rrrrrr}
\toprule
\(q\)&\(k\)&\(t\)&\(r\)&\(d_-\)&\(d_+\)\\
\midrule
19&7&2&0&49&51\\
111&23&2&0&529&539\\
321&7&34&9&945&952\\
5601&7&562&9&15729&16703\\
5601&23&102&55&28129&28155\\
6001&23&108&23&29785&30167\\
\bottomrule
\end{tabular}
\end{center}

In the first row, the same bound \(N_{1/19}(d)\le d+7\) applies in all three dimensions. The construction starts with 56 lines in dimension 49 and adds one line and one dimension at each step, giving exact values 56, 57, and 58 in dimensions 49, 50, and 51, respectively. Theorem~\ref{fa:intervals} gives the same phenomenon for every sufficiently large odd denominator in each row of the parameter table.

\subsection{A family at quadratic irrational cosines}\label{fa:conference}
\begin{theorem}\label{fa:irrational}
Let \(C\) be a real symmetric conference matrix of order \(m\) and suppose \(m-1\) is not a square. For integers \(t\ge2\), set
\[
q=t(1+\sqrt{m-1})-1,\qquad d=mt-m/2.
\]
Then \(N_{1/q}(d)=mt\). In particular, for primes \(p\equiv1\pmod4\) one has
\[
N_{1/(1+2\sqrt p)}\bigl(3(p+1)/2\bigr)=2(p+1).
\]
\end{theorem}
\begin{proof}
By definition \(C\) has zero diagonal, off-diagonal entries \(\pm1\), and \(C^2=(m-1)I_m\). Its eigenvalues are \(\pm\sqrt{m-1}\), each with multiplicity \(m/2\) since its trace is zero. The tensor construction~\cite[Proposition~5.11]{GKMS} gives the Seidel matrix \(I_{mt}+(C-I_m)\otimes J_t\). On block constants it has eigenvalues \(1+t(\pm\sqrt{m-1}-1)\), and on the block sum-zero spaces it has eigenvalue \(1\). Its least eigenvalue is \(-q\) with multiplicity \(m/2\). This gives \(mt\) lines of the required cosine and rank \(d\).

Suppose a system of \(n=mt+1\) lines of that cosine existed in \(\R^d\). Its Seidel matrix \(S'\) would have eigenvalue \(-q\) with multiplicity at least \(m/2+1\), and its integer characteristic polynomial would force the same multiplicity for the algebraic conjugate. Related conjugate-multiplicity arguments for Seidel spectra appear in Lin and Yu~\cite[Theorem~5.3]{LYLemmens}. For \(S'-I_n\) the two eigenvalues in question are \(-t(1\pm\sqrt{m-1})\), whose squares sum to \(2mt^2\). Every entry of \(S'-I_n\) is \(\pm1\), so
\[
(mt+1)^2=\tr\big((S'-I_n)^2\big)\ge (m/2+1)2mt^2.
\]
The right side exceeds the left by \(2mt(t-1)-1>0\), a contradiction. For the stated specialization take the symmetric Paley conference matrix of order \(p+1\)~\cite{Paley} and \(t=2\).
\end{proof}

For \(p=5,13,17\) these formulas give dimension and cardinality \((9,12),(21,28),(27,36)\), at the three corresponding cosines. The 12-line construction also appears as the order-12 Seidel spectrum in \cite[\S5.1]{GKMS}.

\subsection{Two values at cosine \texorpdfstring{\(1/7\)}{1/7}}
\begin{proposition}\label{fa:seven}
One has \(N_{1/7}(10)=12\) and \(N_{1/7}(11)=13\).
\end{proposition}
\begin{proof}
The relative bound~\cite[Lemma~6.1]{vanLintSeidel} gives \(N_{1/7}(d)\le \lfloor48d/(49-d)\rfloor\). The bounds are \(12\) and \(13\) for \(d=10\) and \(d=11\), respectively.

For the lower bounds, the following Gram matrices are instances of the graph construction in \cite[Proposition~3.2]{JTYZZ}. For \(n=12\) take \(B=\operatorname{diag}(J_4,J_4,J_4)\), and for \(n=13\) add a singleton block \(J_1\) to \(B\). In each case \(8I_n-2B+J_n\) has diagonal \(7\) and off-diagonal entries \(\pm1\). The matrix \(8I_n-2B\) is positive semidefinite with kernel spanned by the three full-block constants. Adding \(J_n\) cuts the kernel to the two-dimensional subspace of those constants with total sum zero. The two Gram matrices therefore have ranks \(10\) and \(11\).
\end{proof}

\section{Computational supplement}\label{sec:reproduction}
The supplementary material records the finite inputs, verification scripts and expected outputs for each dimension. Section~\ref{d18:sec:verification} describes the root, height, and binary stages in dimension eighteen. Section~\ref{d19:sec:verification} describes the rational certificates and Niemeier computations for dimension nineteen. Section~\ref{d20:sec:verification} describes the scalar and small-model checks and the two Golay-based verifiers in dimension twenty. Each dimension-specific reproduction directory includes a manifest of file digests.

The \href{https://github.com/Wideep886/Equiangular_Lines_Maxima/tree/d7ffd626aca5398c83e840dd670a58bf339643b1/combined-reproduction}{combined verification guide} lists the dependencies and setup instructions. From the repository root, with Python assertions enabled, run
\begin{verbatim}
python3 -B combined-reproduction/replay.py
\end{verbatim}
Use \verb|--dimension| followed by \texttt{18}, \texttt{19}, or \texttt{20} to replay one dimension. The runner checks the file digests and calls the dimension-specific verification programs.

The \href{https://github.com/Wideep886/Equiangular_Lines_Maxima/blob/d7ffd626aca5398c83e840dd670a58bf339643b1/combined-reproduction/README.md\#dimension-nineteen-certificate-dependencies}{certificate guide} connects the geometry in Section~\ref{d19:sec:certs} to the rows used in the dimension-nineteen verifier.

The guide also gives commands for two standalone checks. The \href{https://github.com/Wideep886/Equiangular_Lines_Maxima/blob/d7ffd626aca5398c83e840dd670a58bf339643b1/combined-reproduction/check_leech_sections.py}{Leech section checker} verifies the sizes, pairings, and rational ranks of the nested sections in Section~\ref{d20:sec:lower}. The \href{https://github.com/Wideep886/Equiangular_Lines_Maxima/blob/d7ffd626aca5398c83e840dd670a58bf339643b1/combined-reproduction/check_fixed_angle_examples.py}{fixed-angle example checker} recomputes the six interval examples in Section~\ref{fa:odd}, checks the arithmetic for the sufficient thresholds, and verifies the two Gram matrices in Proposition~\ref{fa:seven}. Both checkers use Python 3.10 or newer with assertions enabled and require only the standard library. The sharp Leech section bounds follow from the generalized Hamming weights in the text. The general fixed-angle results are proved in Section~\ref{fa:section}.

\section*{Statements and Declarations}
\paragraph{Data and code availability.}
\mbox{\url{https://github.com/Wideep886/Equiangular_Lines_Maxima}}.
This manuscript uses the computational supplement at Git commit
\href{https://github.com/Wideep886/Equiangular_Lines_Maxima/commit/d7ffd626aca5398c83e840dd670a58bf339643b1}{\texttt{d7ffd62}}.
The materials are in the
\href{https://github.com/Wideep886/Equiangular_Lines_Maxima/tree/d7ffd626aca5398c83e840dd670a58bf339643b1/dimension-18/reproduction}{dimension-eighteen},
\href{https://github.com/Wideep886/Equiangular_Lines_Maxima/tree/d7ffd626aca5398c83e840dd670a58bf339643b1/dimension-19/reproduction}{dimension-nineteen},
and \href{https://github.com/Wideep886/Equiangular_Lines_Maxima/tree/d7ffd626aca5398c83e840dd670a58bf339643b1/dimension-20/reproduction}{dimension-twenty}
reproduction directories.

\section*{AI methodology}
The authors used OpenAI's GPT-6 in exploring mathematical approaches and developing the low-dimensional and fixed-angle arguments. GPT-6 also assisted with preparing and checking verification scripts, checking literature and citations, and editing and consolidating the manuscript. Readers can replay the exact finite computations directly from the published scripts and data. Responsibility for the mathematical claims, the code, and the manuscript rests with the authors.

\end{document}